\documentclass[11pt,letterpaper]{article}

\usepackage{amsmath,amsfonts,amsthm,amssymb}
\usepackage{dsfont}
\usepackage{hyperref}
\usepackage[english]{babel}
\usepackage{verbatim}
\usepackage{mathrsfs}
\usepackage{bm}
\usepackage{color}
\usepackage{graphicx}
\usepackage{url}
\usepackage[left=2.5cm, top=2.5cm, bottom=2.5cm, right=2.5cm]{geometry}
\usepackage[active]{srcltx}

\newtheorem{theorem}{Theorem}[section]
\newtheorem{assumption}[theorem]{Assumption}
\newtheorem{corollary}[theorem]{Corollary}

\newtheorem{lemma}[theorem]{Lemma}

\newtheorem{remark}[theorem]{Remark}

\numberwithin{equation}{section}

\definecolor{Red}{rgb}{1.00, 0.00, 0.00}

\definecolor{Green}{rgb}{0.2, 0.5, 0.2}

\definecolor{Blue}{rgb}{0.00, 0.00, 1.00}

\title{Feynman-Kac Formulas for Solutions to Degenerate Elliptic and Parabolic Boundary-Value and Obstacle Problems with Dirichlet Boundary Conditions}

\author{Paul M.N. Feehan\thanks{Department of Mathematics, Rutgers, The State University of New Jersey, Piscataway, NJ 08854-8019, U.~S.~A. {\tt feehan@math.rutgers.edu}} \and Ruoting Gong\thanks{Department of Applied Mathematics, Illinois Institute of Technology, Chicago, IL 60616, U.~S.~A. {\tt rgong2@iit.edu}} \and Jian Song\thanks{Department of Mathematics, The University of Hong Kong, Pokfulam, Hong Kong, P.~R. China. {\tt txjsong@hku.hk}}}

\date{}

\begin{document}

\maketitle

\begin{abstract}

We prove Feynman-Kac formulas for solutions to elliptic and parabolic boundary value and obstacle problems associated with a general Markov diffusion process. Our diffusion model covers several popular stochastic volatility models, such as the Heston model, the CEV model and the SABR model, which are widely used as asset pricing models in mathematical finance. The generator of this Markov process with killing is a second-order, degenerate, elliptic partial differential operator, where the degeneracy in the operator symbol is proportional to the $2\alpha$-power of the distance to the boundary of the half-plane, with $\alpha\in(0,1]$. Our stochastic representation formulas provide the unique solutions to the elliptic boundary value and obstacle problems, when we seek solutions which are suitably smooth up to the boundary portion $\Gamma_{0}$ contained in the boundary of the upper half-plane. In the case when the full Dirichlet condition is given, our stochastic representation formulas provide the unique solutions which are not guaranteed to be any more than continuous up to the boundary portion $\Gamma_{0}$.

\vspace{0.3cm}


\vspace{0.3cm}

\end{abstract}

\section{Introduction}

The \emph{Feynman-Kac} formula, discovered by Mark Kac~\cite{Kac:1949} who was inspired in turn by the doctoral dissertation of Richard Feynman~\cite{Feynman:2005}, establishes a link between stochastic differential equations (SDEs) and parabolic partial differential equations (PDEs). It offers a method of solving certain PDEs by simulating random paths of a stochastic process. Conversely, expectations of an important class of random processes can be computed by deterministic methods.The Feynman-Kac formulas have been well established for strictly elliptic (in the sense of Gilbarg and Trudinger~\cite[pp.~31]{GilbargTrudinger:1983} or Karatzas and Shreve~\cite[Definition 5.7.1]{KaratzasShreve:1991}) Markov generator $\mathscr{A}$ (cf.~\cite{Friedman:1976}, \cite{BensoussanLions:1982}, \cite{KaratzasShreve:1991} and~\cite{Oksendal:2003}). However, the literature is rather incomplete when $\mathscr{A}$ is degenerate elliptic, that is, only has a non-negative definite characteristic form, and its coefficients are unbounded. One such work is Feehan and Pop~\cite{FeehanPop:2013(2)}, which obtains Feynman-Kac formulas for solutions to elliptic and parabolic boundary value and obstacle problems associated with the two-dimensional Heston stochastic volatility process. For certain classes of degenerate elliptic (parabolic) partial differential equations with boundary (terminal-boundary respectively) condition, the Feynman-Kac formula was explored in~\cite{StroockVaradhan:1972} and~\cite{Friedman:1976}, and more recently in~\cite{EkstromTysk:2010}, \cite{EkstromTysk:2011}, \cite{BayraktarKardarasXing:2012} and~\cite{Zhou:2013}. A thorough survey of the results in the above references was given in~\cite{FeehanPop:2013(2)}, where the authors provided detailed comparisons between their results and the previous results.

In this article, we extend the results of~\cite{FeehanPop:2013(2)} to the case when $\mathscr{A}$ is a general elliptic differential operator with possibly unbounded  coefficients:
\begin{align*}
\mathscr{A}u(x):=-\frac{1}{2}tr\left(a(x)D^{2}u(x)\right)-\langle b(x),Dv(x)\rangle+c(x)v(x),\quad u\in C^{2}(\mathscr{O}),\,\,\,x\in\mathscr{O},
\end{align*}
where $\mathscr{O}$ is a possibly unbounded, connected and open subset of the open upper half-space $\mathbb{H}:=\mathbb{R}^{d-1}\times(0,\infty)$ ($d\geq 2$). Under mild conditions (see Section \ref{sec:MainResults}) below, we prove stochastic representation formulas for solutions to an elliptic boundary value problem
\begin{align}\label{eq:EllipBound}
\mathscr{A}u=f\quad\text{on }\,\mathscr{O},
\end{align}
and to an elliptic obstacle problem
\begin{align}\label{eq:EllipObs}
\min\left\{\mathscr{A}u-f,\,u-\psi\right\}=0\quad\text{on }\,\mathscr{O},
\end{align}
respectively, subject to a \emph{partial} Dirichlet boundary condition
\begin{align}\label{eq:BoundCondGamma1}
u=g\quad\text{on }\,\Gamma_{1}.
\end{align}
Here, $\Gamma_{1}:=\partial\mathscr{O}\cap\mathbb{H}$ is the portion of the boundary, $\partial\mathscr{O}$, of $\mathscr{O}$ which lies in $\mathbb{H}$, $f:\mathscr{O}\rightarrow\mathbb{R}$ is a source function, $g:\Gamma_{1}\rightarrow\mathbb{R}$ prescribes a Dirichlet boundary condition along $\Gamma_{1}$, and $\psi:\mathscr{O}\cup\Gamma_{1}\rightarrow\mathbb{R}$ is an obstacle function which is compatible with $g$ in the sense that
\begin{align}\label{eq:CompatObsGamm1}
\psi\leq g\quad\text{on }\,\Gamma_{1},
\end{align}
while $\mathscr{A}$ is an elliptic differential operator on $\mathscr{O}$ which is degenerate along the interior, $\Gamma_{0}$, of $\partial\mathbb{H}\cap\partial\mathscr{O}$ and may have unbounded coefficients. We require $\Gamma_{0}$ to be non-empty throughout this article as, otherwise, $\mathscr{A}$ is non-degenerate, and standard results apply (cf.~\cite{Friedman:1976}, \cite{BensoussanLions:1982}, \cite{KaratzasShreve:1991} and~\cite{Oksendal:2003}). However, an additional boundary condition is \emph{not} necessarily prescribed along $\Gamma_{0}$. Rather, we will see that our stochastic representation formulas provide the unique solutions to (\ref{eq:EllipBound}) or (\ref{eq:EllipObs}), together with (\ref{eq:BoundCondGamma1}), when we seek solutions which are suitably \emph{smooth} up to the boundary portion $\Gamma_{0}$, a property which is guaranteed when the solutions lie in certain weighted H\"{o}lder spaces (by analogy with~\cite{DaskalopoulosHamilton:1998}), or replace the boundary condition (\ref{eq:BoundCondGamma1}) with the \emph{full} Dirichlet condition
\begin{align}\label{eq:BoundCondFull}
u=g\quad\text{on }\,\partial\mathscr{O},
\end{align}
where $g:\partial\mathscr{O}\rightarrow\mathbb{R}$, in which case the solutions are not guaranteed to be any more than \emph{continuous} on the full boundary $\partial\mathscr{O}$, and $\psi:\bar{\mathscr{O}}\rightarrow\mathbb{R}$ is now required to be compatible with $g$ in the sense that
\begin{align}\label{eq:CompatObsFull}
\psi\leq g\quad\text{on }\,\partial\mathscr{O}.
\end{align}

Moreover, we also prove stochastic representation formulas for solutions for a parabolic terminal/boundary value problem
\begin{align}\label{eq:ParaTermBound}
-u_{t}+\mathscr{A}u=f\quad\text{on }\,Q,
\end{align}
and to a parabolic obstacle problem
\begin{align}\label{eq:ParaTermObs}
\min\left\{-u_{t}+\mathscr{A}u-f,\,u-\psi\right\}=0\quad\text{on }\,Q,
\end{align}
respectively, subject to the \emph{partial} terminal/boundary condition
\begin{align}\label{eq:TermBoundPart}
u=g\quad\text{on }\,\eth^{1}Q.
\end{align}
Here, we define $Q:=(0,T)\times\mathscr{O}$, where $T\in(0,\infty)$, and
\begin{align*}
\eth^{1}Q:=(0,T)\times\Gamma_{1}\cup\{T\}\times(\mathscr{O}\cup\Gamma_{1}),
\end{align*}
and assume given a source function $f:Q\rightarrow\mathbb{R}$, a boundary data function $g:\eth^{1}Q\rightarrow\mathbb{R}$, and an obstacle function $\psi:Q\cup\eth^{1}Q\rightarrow\mathbb{R}$ which is compatible with $g$ in the sense that
\begin{align}\label{eq:CompatObsPartPara}
\psi\leq g\quad\text{on }\,\eth^{1}Q.
\end{align}
Just as in the elliptic case, we will either consider solutions which are suitably \emph{smooth} up to $(0,T)\times\Gamma_{0}$, but impose no explicit boundary condition along $(0,T)\times\Gamma_{0}$, or replace the boundary condition in (\ref{eq:TermBoundPart}) with the \emph{full} Dirichlet condition
\begin{align}\label{eq:TermBoundFull}
u=g\quad\text{on }\,\eth Q,
\end{align}
where
\begin{align*}
\eth Q:=(0,T)\times\partial\mathscr{O}\cup\{T\}\times\bar{\mathscr{O}},
\end{align*}
is the full parabolic boundary of $Q$, in which case the solutions are not guaranteed to be any more than \emph{continuous} on $\eth Q$ and $\psi:Q\cup\eth Q\rightarrow\mathbb{R}$ is compatible with $g$ in the sense that
\begin{align}\label{eq:CompatObsFullPara}
\psi\leq g\quad\text{on }\,\eth Q.
\end{align}

The major applications of the Feynman-Kac formulas presented in this article are in the area of mathematical finance. A solution $u$ to the elliptic obstacle problem (\ref{eq:EllipObs}) and (\ref{eq:BoundCondGamma1}) with $f=0$, can be interpreted as the value function for a perpetual American-type option with the payoff function given by the obstacle function $\psi$, and a solution $u$ to the parabolic obstacle problem (\ref{eq:ParaTermObs})  and (\ref{eq:TermBoundPart}) with $f=0$, can be regarded as the value function for a finite-expiration American-type option with the payoff function given by a terminal condition function $h=g(T,\cdot):\mathscr{O}\rightarrow\mathbb{R}$, which typically coincides on ${T}\times\mathscr{O}$ with the obstacle function $\psi$. A  solution $u$ to the  parabolic terminal/boundary value problem (\ref{eq:ParaTermBound}) and (\ref{eq:TermBoundPart}) with $f=0$ can be interpreted as the value function for a European-type option with the payoff function given by the function $h$ as above.

The present article is organized as follows. Section \ref{sec:MainResults} presents the main results of this article, which establish Feynman-Kac formulas for the unique solutions to the elliptic boundary-value and obstacle problems with partial/full boundary conditions. Section \ref{sec:Prelim} contains some basic estimates and boundary classifications for degenerate diffusion processes. Section \ref{sec:EllipBoundValProb} and Section \ref{sec:EllipObsProb} present the proofs of the main results of Section \ref{sec:MainResults} for the elliptic boundary-value and obstacle problems, respectively. Finally, Section \ref{sec:ParProb} establishes Feynman-Kac formulas for the unique solutions to the parabolic terminal/boundary-value and obstacle problems with partial/full terminal/boundary conditions.

\medskip
\noindent
\textbf{Acknowledgement:} The authors gratefully thank Daniel Ocone and Camelia Pop for the constructive and insightful comments, which significantly contributed to improve the quality of the manuscript.

\section{Main Results}\label{sec:MainResults}

\subsection{Basic Settings}\label{subsec:Setting}

Let $\mathscr{S}^{+}(d)\subset\mathbb{R}^{d\times d}$ be the collection of positive-definite symmetric matrices. Let $\mathscr{A}$ be the following elliptic operator
\begin{align}\label{eq:GenAkillc}
\mathscr{A}v(x):=-\frac{1}{2}\text{tr}\!\left(a(x)D^{2}v(x)\right)-\langle b(x),Dv(x)\rangle+c(x)v(x),\quad v\in C^{2}(\mathscr{O}),\,\,\,x\in\mathscr{O},
\end{align}
where
\begin{itemize}
\item $b:\mathscr{O}\rightarrow\mathbb{R}^{d}$ is a continuous vector field with
    \begin{align}\label{eq:OperatorDrift}
    b(x)=(b_{1}(x),\ldots,b_{d-1}(x),b_{d}(x_{d}))^{T},\quad x\in\mathscr{O},
    \end{align}
     where the last component $b_{d}$ only depends on $x_{d}$;
\item $c:\mathscr{O}\rightarrow[0,\infty)$ is a non-negative continuous function;
\item $a:\mathscr{O}\rightarrow\mathscr{S}^{+}(d)$ is a continuous matrix-valued function, such that
    \begin{align}\label{eq:Deftildea}
    a(x)=x_{d}^{\beta}\tilde{a}(x),\quad x\in\mathscr{O},
    \end{align}
    for some continuous $\tilde{a}:\mathscr{O}\rightarrow\mathscr{S}^{+}(d)$ and $\beta\in(0,2]$.
\end{itemize}

To connect the operator $-\mathscr{A}$ with a diffusion process on $\overline{\mathbb{H}}$, we extend the definitions of $a$ (and thus $\tilde{a}$), $b$ and $c$ to the half space. Throughout, $a$ (and $\tilde{a}$), $b$ and $c$ will denote continuous extensions on $\overline{\mathbb{H}}$ such that (\ref{eq:OperatorDrift}) and (\ref{eq:Deftildea}) hold true. Let $m\in\mathbb{N}$, and let $\mathcal{M}^{d\times m}$ be the collection of $d\times m$ matrices. Assume that there exists continuous matrix-valued functions
\begin{align}\label{eq:tildeSigma}
\tilde{\sigma}:\overline{\mathbb{H}}\rightarrow\mathcal{M}^{d\times m}\quad\text{and}\quad\sigma(x):=x_{d}^{\beta/2}\tilde{\sigma}(x),\quad x\in\overline{\mathbb{H}},
\end{align}
such that
\begin{align}\label{eq:Relasigma}
a(x)=\sigma(x)\sigma^{T}(x)\quad\text{and}\quad\tilde{a}(x)=\tilde{\sigma}(x)\tilde{\sigma}^{T}(x),\quad x\in\overline{\mathbb{H}},
\end{align}
and that
\begin{align}\label{eq:LastRowtildeSigma}
\tilde{\sigma}_{dj}(x):=\rho_{j}\tilde{\sigma}_{0}(x_{d}),\quad x\in\mathbb{H},\quad j=1,\ldots,m,
\end{align}
with $\rho_{j}>0$, $j=1,\ldots,m$, $\sum_{j=1}^{m}\rho_{j}^{2}=1$, and $\tilde{\sigma}_{0}:\overline{\mathbb{H}}\rightarrow\mathbb{R}$ continuous.  The operator $-\mathscr{A}$ defined by (\ref{eq:GenAkillc})-(\ref{eq:LastRowtildeSigma}) is the infinitesimal generator, with killing rate $c$, of the diffusion process which satisfies following SDEs:
\begin{align}
\label{eq:MainDegenSDE1} dX_{s}^{(i,t)}&=b_{i}(X_{s}^{(t)})\,ds+\left(X_{s}^{(d,t)}\right)^{\beta/2}\left(\sum_{j=1}^{m}\tilde{\sigma}_{ij}(X_{s}^{t})\,dW^{(j,t)}_{s}\right),\quad i=1,\ldots,d-1,\quad s\geq t,\\
\label{eq:MainDegenSDE2} dX_{s}^{(d,t)}&=b_{d}(X_{s}^{d,t})\,ds+\left(X_{s}^{(d,t)}\right)^{\beta/2}\tilde{\sigma}_{0}\left(X_{s}^{(d,t)}\right)\left(\sum_{j=1}^{m}\rho_{j}\,dW^{(j,t)}_{s}\right),\quad s\geq t,
\end{align}
subject to the initial condition
\begin{align}\label{eq:MainDegenSDEIniCond}
X_{t}^{(t)}=(X_{t}^{(1,t)},\ldots,X_{t}^{(d,t)})^{T}=x=\left(x_{1},\ldots,x_{d}\right)^{T}\in\overline{\mathbb{H}},
\end{align}
We will use the notations  $\mathbb{P}^{t,x}$ and $\mathbb{E}^{t,x}$ to indicate the probability and the expectation corresponding to the initial state (\ref{eq:MainDegenSDEIniCond}), respectively.  Also, we will omit the superscript  $t$ when $t=0$.
\begin{remark} \label{remark1.1}
The diffusion model (\ref{eq:MainDegenSDE1})-(\ref{eq:MainDegenSDEIniCond}) considered here covers several popular stochastic volatility model, e.g., the Heston model~\cite{Heston:1993}, the CEV model~\cite{Cox:1975}, and the SABR model~\cite{HaganKumarLesniewskiWoodward:2003}, which are widely used as asset pricing models in mathematical finance. It also covers a simple version of the continuous affine diffusion process introduced by Duffie, Filipovi\'{c} and Schachermayer~\cite{DuffieFilipovicSchachermayer:2003} (see also~\cite{FilipovicMayerhofer:2009}).
\end{remark}
Throughout this article, we make the following standard assumptions.
\begin{assumption}\label{assupt:Cont} \emph{(Continuity of coefficients)} $\,$ The functions $b$, $\sigma$ (and thus $\tilde{\sigma}$), and $c$, are continuous on $\overline{\mathbb{H}}$.
\end{assumption}
\begin{assumption}\label{assupt:LinearGrowth} \emph{(Linear growth condition)} $\,$
$b$ and $\sigma$ satisfies the following linear growth condition: let $u$ be any vector-valued or matrix-valued function on $\overline{\mathbb{H}}$, then there exists $K>0$, such that
\begin{align}\label{eq:lineargrowth}
\|u(x)\|\leq K(1+\|x\|),\quad x\in\overline{\mathbb{H}},
\end{align}
where $\|\cdot\|$ denotes the Euclidean norm of a vector or a matrix.
\end{assumption}
\begin{assumption}\label{assupt:NonNegbd}
The function $b_{d}$ satisfies
\begin{align*}
b_{d}(x_{d})\geq 0,\quad\text{when }\,x_{d}=0.
\end{align*}
\end{assumption}
\begin{assumption}\label{assupt:UnifPostc} \emph{(Uniform positivity)}
The Borel measurable function $c:\mathbb{R}^{d}\rightarrow\mathbb{R}$ is strictly positive, and moreover, there exists $c_{0}>0$, such that
\begin{align*}
c(x)\geq c_{0},\quad\text{for any }\,x\in\mathbb{R}^{d}.
\end{align*}
\end{assumption}
\begin{assumption}\label{assupt:Elliptildesigma} \emph{(Uniform Ellipticity)}
The function $\tilde{a}$ is \emph{uniformly elliptic} in $\overline{\mathbb{H}}$ (cf.~\cite[Chapter 3]{GilbargTrudinger:1983} and~\cite[Definition 5.7.1]{KaratzasShreve:1991}). That is, there exists a universal positive constant $\delta>0$, such that
\begin{align*}
\sum_{i=1}^{d}\sum_{k=1}^{d}\tilde{a}_{ik}(x)\xi_{i}\xi_{k}\geq\delta\|\xi\|^{2},\quad\text{for any }\,\,x\in\overline{\mathbb{H}},\quad \xi\in\mathbb{R}^{d}.
\end{align*}
\end{assumption}
Assumption \ref{assupt:Cont} and  Assumption \ref{assupt:LinearGrowth} ensures that, for any $t\geq 0$ and $x\in\overline{\mathbb{H}}$, there exists a weak solution\footnote{For the definition of weak solution, see e.g., \cite[Definition 5.3.1]{KaratzasShreve:1991}, \cite[Definitions IX.1.2 \& IX.1.5]{RevusYor:1999} and~\cite[pp.~115]{Yamada:1978}.} $(\Omega,\mathscr{F},(\mathscr{F}_{s})_{s\geq t},\mathbb{P}^{t,x},W^{(t)},X^{(t)})$ to (\ref{eq:MainDegenSDE1})-(\ref{eq:MainDegenSDE2}) with the initial condition (\ref{eq:MainDegenSDEIniCond}) (cf.~\cite[Theorem 3.3]{FeehanPop:2013(1)}),  where
\begin{itemize}
\item $(\Omega,\mathscr{F},\mathbb{P})$ is a probability space, and $(\mathscr{F}_{s})_{s\geq t}$ is a filtration of sub-$\sigma$-fields of $\mathscr{F}$ satisfying the usual conditions;
\item $W^{(t)}=(W_{s}^{(t)})_{s\geq t}$ is an $m$-dimensional Brownian motion with respect to $(\mathscr{F}_{s})_{s\geq t}$, $X^{(t)}=(X_{s}^{(t)})_{s\geq t}$ is a continuous, $(\mathscr{F}_{s})_{s\geq t}$-adapted $\mathbb{R}^{d}$-valued process;
\item for every $i=1,\ldots,d$, $j=1,\ldots,m$ and $s\geq t$,
    \begin{align*}
    \mathbb{P}\left(\int_{t}^{s}\left(\left|b_{i}(X_{u}^{(t)})\right|+\sigma_{ij}^{2}(X_{u}^{(t)})\right)du<\infty\right)=1;
    \end{align*}
\item the integral version of (\ref{eq:MainDegenSDE1})-(\ref{eq:MainDegenSDEIniCond}) holds true $\mathbb{P}^{t,x}$-a.$\,$e..
\end{itemize}
Notice that we do \emph{not} have the uniqueness in law of weak solutions in our SDE system. Indeed, when $\beta\in(0,1)$, the last SDE (\ref{eq:MainDegenSDE2}) may not have a weak solution which is unique in law (cf.~\cite[Example 5.2.15, Theorem 5.5.4 \& Remark 5.5.6]{KaratzasShreve:1991}). Moreover, Assumption \ref{assupt:NonNegbd} ensures that any weak solution  $(\Omega,\mathscr{F},(\mathscr{F}_{s})_{s\geq t},\mathbb{P}^{t,x},W^{(t)},X^{(t)})$ started in $\overline{\mathbb{H}}$ remains in the upper half space (cf.~\cite[Proposition 3.1]{FeehanPop:2013(1)}), i.e.,
\begin{align*}
X_{t}^{(t)}=x\in\overline{\mathbb{H}}\quad\Rightarrow\quad\mathbb{P}^{t,x}\left(X_{s}^{(t)}\in\overline{\mathbb{H}}\right)=1,\quad\text{for any }\,s\geq t.
\end{align*}

Let $U\subseteq\mathbb{H}$ be an open set, and for any $x\in U$ and $t\geq 0$, let $(\Omega,\mathscr{F},(\mathscr{F}_{s})_{s\geq t},\mathbb{P}^{t,x},W^{(t)},X^{(t)})$ be a weak solution to (\ref{eq:MainDegenSDE1})-(\ref{eq:MainDegenSDE2}) with the initial condition (\ref{eq:MainDegenSDEIniCond}). Define
\begin{align}\label{eq:StopTimeTau}
\tau_{U}^{t,x,X}&=\inf\left\{s\geq t:\,X_{s}^{(t)}\notin U,\,X_{t}^{(t)}=x\right\},\\
\label{eq:StopTimeLambda} \lambda_{U}^{t,x,X}&=\inf\left\{s\geq t:\,X_{s}^{(t)}\notin U\cup\left(\overline{U}\cap\partial\mathbb{H}\right),\,X_{t}^{(t)}=x\right\}.
\end{align}
It is clear that $\tau_{U}^{t,x,X}=\lambda_{U}^{t,x,X}$ if $\overline{U}\cap\partial\mathbb{H}=\emptyset$. From the boundary classification lemma (see Lemma \ref{lem:BoundClass} below), we also have $\tau_{\mathscr{O}}^{t,x,X}=\lambda_{\mathscr{O}}^{t,x,X}$, when $\Sigma^{(d)}(0)=\infty$. Indeed, when $\Sigma^{(d)}(0)=\infty$, $X^{(t)}$ can never reach the boundary portion $\Gamma_{0}$ from the interior of $\mathscr{O}$. Moreover, both $\tau_{U}^{t,x,X}$ and $\lambda_{U}^{t,x,X}$ are stopping times with respect to $(\mathscr{F}_{s})_{s\geq t}$, since $(\mathscr{F}_{s})_{s\geq t}$ is assumed to satisfy the usual conditions (cf.~\cite[pp.~117]{Oksendal:2003}). In the sequel, when the initial condition $(t,x)$ is clear from the context, we will omit the superscripts in the preceding definitions. Also, when $t=0$, we will omit the superscript $t$ in the preceding definitions.

In the following two subsections, we present Feynman-Kac formulas for solutions to elliptic boundary-value and obstacle problems. The proofs are deferred to Section \ref{sec:EllipBoundValProb} and Section \ref{sec:EllipObsProb}, respectively. Similar results can be obtained parabolic terminal/boundary-value and obstacle problems, and will be presented in Section \ref{sec:ParProb}.

\subsection{Feynman-Kac Formulas for Solutions to Elliptic Boundary-Value Problems}\label{sec:UniqueExistEllipBound}

For any $x\in\bar{\mathscr{O}}$ and any weak solution $(\Omega,\mathscr{F},(\mathscr{F}_{s})_{s\geq 0},\mathbb{P}^{x},W,X)$ to (\ref{eq:MainDegenSDE1})-(\ref{eq:MainDegenSDE2}) with the initial condition (\ref{eq:MainDegenSDEIniCond}) at $t=0$, let
\begin{align}\label{eq:EllipBoundFeynmanKacTau}
u_{*}^{(X)}(x)\!&:=\!\mathbb{E}^{x}\!\!\left[\exp\!\left(\!-\!\!\int_{0}^{\tau_{\mathscr{O}}}\!\!c(X_{s})\,ds\!\right)\!g\!\left(X_{\tau_{\mathscr{O}}}\right)\!{\bf 1}_{\{\tau_{\mathscr{O}}<\infty\}}\right]\!\!+\!\mathbb{E}^{x}\!\!\left[\int_{0}^{\tau_{\mathscr{O}}}\!\!\exp\!\left(\!-\!\!\int_{0}^{s}\!c(X_{u})\,du\!\right)\!f(X_{s})\,ds\right],\\
\label{eq:EllipBoundFeynmanKacLambda} u_{**}^{(X)}(x)\!&:=\!\mathbb{E}^{x}\!\!\left[\exp\!\left(\!-\!\!\int_{0}^{\lambda_{\mathscr{O}}}\!\!c(X_{s})\,ds\!\right)\!g\!\left(X_{\lambda_{\mathscr{O}}}\right)\!{\bf 1}_{\{\lambda_{\mathscr{O}}<\infty\}}\!\right]\!\!+\!\mathbb{E}^{x}\!\!\left[\int_{0}^{\lambda_{\mathscr{O}}}\!\!\!\exp\!\left(\!-\!\!\int_{0}^{s}\!c(X_{u})\,du\!\right)\!f(X_{s})ds\right].
\end{align}
For any integer $k\geq 0$, let $C^{k}(\mathscr{O})$ be the vector space of functions whose derivatives up to order $k$ are continuous on $\mathscr{O}$, and let $C^{k}(\bar{\mathscr{O}})$ be the Banach space of functions whose derivatives up to order $k$ are {\it uniformly continuous and bounded} on $\bar{\mathscr{O}}$ (cf.~\cite[\S 1.25 \& \S 1.26]{Adams:1975}). For any integer $k\geq 0$ and $\alpha\in(0,1)$, let $C^{k,\alpha}(\mathscr{O})$ be the subspace of $C^{k}(\mathscr{O})$ consisting of functions whose derivatives up to order $k$ are \emph{locally} $\alpha$-H\"{o}lder continuous on $\mathscr{O}$ (cf.~\cite[pp.~52]{GilbargTrudinger:1983}) and let $C^{k,\alpha}(\bar{\mathscr{O}})$ be the subspace of $C^{k}(\bar{\mathscr{O}})$ consisting of functions whose derivatives up to order $k$ are \emph{uniformly} $\alpha$-H\"{o}lder continuous on $\mathscr{O}$ (cf.~\cite[pp.~52]{GilbargTrudinger:1983} or~\cite[\S 1.27]{Adams:1975}). If $T\subsetneq\partial\mathscr{O}$ is a relatively open set, let $C_{\text{loc}}^{k}(\mathscr{O}\cup T)$ (respectively, $C^{k,\alpha}_{\text{loc}}(\mathscr{O}\cup T)$) denote the space of functions of $u$ such that for any pre-compact open subset $U\Subset\mathscr{O}\cup T$, $u\in C^{k}(\overline{U})$ (respectively, $u\in C^{k,\alpha}(\overline{U})$).

Our first result shows that (\ref{eq:EllipBoundFeynmanKacTau}) is the unique solution (assuming existence), which is $C^{2}$ inside the domain $\mathscr{O}$ and is continuous up to an appropriate  portion of $\partial\mathscr{O}$, to the elliptic boundary-value problem (\ref{eq:EllipBound}) with either the partial boundary condition (\ref{eq:BoundCondGamma1}) or the full boundary condition (\ref{eq:BoundCondFull}). In particular, $u_{*}^{(X)}$ is {\it independent} of the choice of the weak solution.
\begin{theorem}\label{thm:UniqueEllipBound}
Let $b$, $\sigma$ and $c$ satisfy (\ref{eq:SupMartCon}), and let $f\in C(\mathscr{O})$ obey the linear growth condition (\ref{eq:lineargrowth}) on $\mathscr{O}$. Moreover, assume that $b_{d}(0)>0$, and that $b_{d}$ is locally H\"{o}lder continuous at the origin, i.e., there exists constants $\gamma\in(0,1]$, $L>0$ and $\kappa\in(0,1)$, such that
\begin{align}\label{eq:LocalHolderbd}
\left|b_{d}(x_{d})-b_{d}(0)\right|\leq L\,|x_{d}|^{\gamma},\quad\text{for any }\,x_{d}\in[-\kappa,\kappa].
\end{align}
\begin{itemize}
\item [(1)] Suppose that one of the following three cases occurs,
    \begin{itemize}
    \item [(a)] $\beta\in(1,2]$;
    \item [(b)] $\beta=1$, and $\,2b_{d}(0)>\tilde{\sigma}_{0}^{2}(0)$;
    \item [(c)] $\beta=1$, $2b_{d}(0)=\tilde{\sigma}_{0}^{2}(0)$, and $\tilde{\sigma}_{0}^{2}$ is constant in a neighborhood of the origin.
    \end{itemize}
    Assume that $g\in C_{\emph{loc}}(\Gamma_{1})$ obeys (\ref{eq:lineargrowth}) on $\Gamma_{1}$. Let
    \begin{align*}
    u\in C_{\emph{loc}}(\mathscr{O}\cup\Gamma_{1})\cap C^{2}(\mathscr{O})
    \end{align*}
    be a solution to the elliptic boundary value problem (\ref{eq:EllipBound}) and (\ref{eq:BoundCondGamma1}), and which obeys (\ref{eq:lineargrowth}) on $\mathscr{O}$. Then  for any $x\in\mathscr{O}\cup\Gamma_{1}$, $u(x)=u_{*}^{(X)}(x)$, for any weak solution  $(\Omega,\mathscr{F},(\mathscr{F}_{s})_{s\geq 0},\mathbb{P}^{x},W,X)$ to (\ref{eq:MainDegenSDE1})-(\ref{eq:MainDegenSDE2}) with the initial condition (\ref{eq:MainDegenSDEIniCond}) at  $t=0$, where $u_{*}^{(X)}$ is given by (\ref{eq:EllipBoundFeynmanKacTau}).
\item [(2)] Suppose that one of the following two cases occurs,
    \begin{itemize}
    \item [(d)] $\beta\in(0,1)$;
    \item [(e)] $\beta=1$, and $\,2b_{d}(0)<\tilde{\sigma}_{0}^{2}(0)$.
    \end{itemize}
    Assume that $g\in C_{\emph{loc}}(\partial\mathscr{O})$ obeys (\ref{eq:lineargrowth}) on $\partial\mathscr{O}$. Let
    \begin{align*}
    u\in C_{\emph{loc}}(\bar{\mathscr{O}})\cap C^{2}(\mathscr{O})
    \end{align*}
    be a solution to the elliptic boundary value problem (\ref{eq:EllipBound}) and (\ref{eq:BoundCondFull}), and which obeys (\ref{eq:lineargrowth}) on $\mathscr{O}$. Then  for any $x\in\bar{\mathscr{O}}$, $u(x)=u_{*}^{(X)}(x)$, for any weak solution  $(\Omega,\mathscr{F},(\mathscr{F}_{s})_{s\geq 0},\mathbb{P}^{x},W,X)$ to (\ref{eq:MainDegenSDE1})-(\ref{eq:MainDegenSDE2}) with the initial condition (\ref{eq:MainDegenSDEIniCond}) at $t=0$, where $u_{*}^{(X)}$ is given by (\ref{eq:EllipBoundFeynmanKacTau}).
\end{itemize}
\end{theorem}
\begin{remark}\label{rem:TwoSenarios}
In the above theorem, we discuss two scenarios which depend on whether the diffusion process $(X_{t})_{t\geq 0}$ can reach the boundary portion $\Gamma_{0}$ or not (see Section \ref{subsec:BoundClass} and the beginning of Section \ref{sec:EllipBoundValProb} for more details). In the first scenario, if one of the cases (a), (b) and (c) occurs, $X$ never  hits $\Gamma_{0}$, while in the second scenario, if either (d) or (e) occurs, $X$ may reach $\Gamma_{0}$. Consequently, to obtain the uniqueness, we assume the partial boundary condition $g\in C_{\emph{loc}}(\Gamma_{1})$ for scenario (1) and the full boundary condition $g\in C_{\emph{loc}}(\partial\mathscr{O})$ for scenario (2), respectively.
\end{remark}
The next result shows that, for the scenario (2) in Theorem \ref{thm:UniqueEllipBound}, the uniqueness of the elliptic boundary-value problem can be obtained under the partial boundary condition (\ref{eq:BoundCondGamma1}), when a compensational regularity condition on $\Gamma_{0}$ is imposed, and the unique solution is given by $u_{**}^{(X)}$ as in (\ref{eq:EllipBoundFeynmanKacLambda}). In particular, $u_{**}^{(X)}$ is independent of the choice of the weak solution.

Following~\cite[Definition 2.2]{DaskalopoulosFeehan:2012} and~\cite[Remark 1.4]{FeehanPop:2013(2)}, for any $\beta\in(0,2]$, let $C_{\text{s,loc}}^{1,1,\beta}(\mathscr{O}\cup\Gamma_{0})$ be the linear subspace of $C^{1,1}(\mathscr{O})\cap C_{\text{loc}}^{2}(\mathscr{O}\cup\Gamma_{0})$ consisting of functions, $\varphi$, such that, for any pre-compact open subset $U\Subset\mathscr{O}\cup\Gamma_{0}$,
\begin{align*}
\sup_{x\in U}\left(|\varphi(x)|+\|D\varphi(x)\|+\left\|x_{d}^{\beta}D^{2}\varphi(x)\right\|\right)<\infty,
\end{align*}
where $D\varphi$ and $D^{2}\varphi$ denote the gradient and the Hessian matrix of $\varphi$, respectively.
\begin{theorem}\label{thm:UniqueEllipBoundPartial}
Suppose that either case \emph{(d)} or \emph{(e)} in Theorem \ref{thm:UniqueEllipBound} occurs. In addition to the hypotheses of Theorem \ref{thm:UniqueEllipBound}, let $c\in C_{\emph{loc}}(\mathscr{O}\cup\Gamma_{0})$, $f\in C_{\emph{loc}}(\mathscr{O}\cup\Gamma_{0})$, and $g\in C_{\emph{loc}}(\Gamma_{1})$ which obeys the linear growth condition (\ref{eq:lineargrowth}) on $\Gamma_{1}$. Let
\begin{align*}
u\in C_{\emph{loc}}(\mathscr{O}\cup\Gamma_{1})\cap C^{2}(\mathscr{O})\cap C_{\emph{s,loc}}^{1,1,\beta}(\mathscr{O}\cup\Gamma_{0})
\end{align*}
be a solution to the elliptic boundary value problem (\ref{eq:EllipBound}) and (\ref{eq:BoundCondGamma1}), and which obeys (\ref{eq:lineargrowth}) on $\mathscr{O}$. Then for any $x\in\mathscr{O}\cup\Gamma_{1}$, $u(x)=u_{**}^{(X)}(x)$, for any weak solution $(\Omega,\mathscr{F},(\mathscr{F}_{s})_{s\geq 0},\mathbb{P}^{x},W,X)$ to (\ref{eq:MainDegenSDE1})-(\ref{eq:MainDegenSDE2}) with the initial condition (\ref{eq:MainDegenSDEIniCond}) at $t=0$, where $u_{**}^{(X)}$ is given by (\ref{eq:EllipBoundFeynmanKacLambda}).
\end{theorem}
\begin{remark}\label{rem:UniqueinLaw}
In both Theorem \ref{thm:UniqueEllipBound} and Theorem \ref{thm:UniqueEllipBoundPartial}, we obtain the unique Feynman-Kac formula for the elliptic boundary value problem (\ref{eq:EllipBound}) with partial/full boundary condition, which is, in particular, independent of the choice of the weak solution. Indeed, the difference in law between weak solutions (possibly defined on different probability spaces) lies only on the time they spend in zeros of the volatility (cf.~\cite[Section 5.5]{KaratzasShreve:1991}, \cite{EngelbertSchmidt:1984} and~\cite{EngelbertSchmidt:1985}), i.e., the boundary portion $\Gamma_{0}$.
\begin{itemize}
\item For $u_{*}^{(X)}$ given as in (\ref{eq:EllipBoundFeynmanKacTau}), the stochastic representation formula only depends on $X$ before it hits $\Gamma_{0}$, as well as the corresponding hitting time $\tau_{\mathscr{O}}$, both of which are unique in law for any weak solution $(\Omega,\mathscr{F},(\mathscr{F}_{s})_{s\geq 0},\mathbb{P}^{x},W,X)$ for each fixed $x$. Hence, $u_{*}^{(X)}$ is independent of the choice of the weak solution.
\item For $u_{**}^{(X)}$ given as in (\ref{eq:EllipBoundFeynmanKacLambda}), the uniqueness of the Feynman-Kac formula does \emph{not} conflict with the non- uniqueness in law of weak solutions. If the upper limit of the integral in (\ref{eq:EllipBoundFeynmanKacLambda}) was a fixed $T>0$, then the uniqueness of (\ref{eq:EllipBoundFeynmanKacLambda}) for a sufficiently large class of $f$ and $g$ would result in the uniqueness in law of $X$ for each fixed $x$. However, the terminal time in (\ref{eq:EllipBoundFeynmanKacLambda}) is the hitting time $\lambda_{\mathscr{O}}$ of the boundary portion $\Gamma_{1}$. Prior to $\lambda_{\mathscr{O}}$, $X$ can hit and stay in $\Gamma_{0}$ for multiple times, which results in different (in law) weak solutions, and in particular, different (in law) hitting times $\lambda_{\mathscr{O}}$. Therefore, the Feynman-Kac formula (\ref{eq:EllipBoundFeynmanKacLambda}) can still be unique without the uniqueness in law of weak solutions.
\end{itemize}
\end{remark}
Next, we show that the Feynman-Kac formula (\ref{eq:EllipBoundFeynmanKacTau}) is indeed the solution to the elliptic boundary-value problem. We first obtain the following theorem concerning the existence of solutions to the elliptic boundary value problem with traditional regularity on $\mathscr{O}$.
\begin{theorem}\label{thm:ExistenceEllipBound}
In addition to the hypothesis of Theorem \ref{thm:UniqueEllipBound}, assume that the boundary portion $\Gamma_{1}$ is of class $C^{2,\alpha}$, and that $b,\sigma,f\in C^{0,\alpha}(\mathscr{O})$, for some $\alpha\in(0,1)$.
\begin{itemize}
\item [(1)] Suppose that either case \emph{(a)}, \emph{(b)} or \emph{(c)} in Theorem \ref{thm:UniqueEllipBound} occurs. Assume that $g\in C_{\emph{loc}}(\overline{\Gamma}_{1})$ which obeys the linear growth condition (\ref{eq:lineargrowth}) on $\Gamma_{1}$. For any $x\in\mathscr{O}\cup\Gamma_{1}$, let $(\Omega,\mathscr{F},(\mathscr{F}_{s})_{s\geq 0},\mathbb{P}^{x},W,X)$ be a weak solution to (\ref{eq:MainDegenSDE1})-(\ref{eq:MainDegenSDE2}) with the initial condition (\ref{eq:MainDegenSDEIniCond}) at $t=0$, and let $u_{*}^{(X)}(x)$ be defined as in (\ref{eq:EllipBoundFeynmanKacTau}). Then, $u_{*}^{(X)}$ is a solution to (\ref{eq:EllipBound}) with the boundary condition (\ref{eq:BoundCondGamma1}) along $\Gamma_{1}$. In particular, $u_{*}^{(X)}\in C_{\emph{loc}}(\mathscr{O}\cup\Gamma_{1})\cap C^{2,\alpha}(\mathscr{O})$ which obeys (\ref{eq:lineargrowth}) on $\mathscr{O}\cup\Gamma_{1}$.
\item [(2)] Suppose that either case \emph{(d)} or \emph{(e)} in Theorem \ref{thm:UniqueEllipBound} occurs. Assume that $g\in C_{\emph{loc}}(\partial\mathscr{O})$ which obeys the linear growth condition (\ref{eq:lineargrowth}) on $\partial\mathscr{O}$. For any $x\in\bar{\mathscr{O}}$, let $(\Omega,\mathscr{F},(\mathscr{F}_{s})_{s\geq 0},\mathbb{P}^{x},W,X)$ be a weak solution to (\ref{eq:MainDegenSDE1})-(\ref{eq:MainDegenSDE2}) with the initial condition (\ref{eq:MainDegenSDEIniCond}) at $t=0$, and let $u_{*}^{(X)}(x)$ be defined as in (\ref{eq:EllipBoundFeynmanKacTau}). Then, $u_{*}^{(X)}$ is a solution to (\ref{eq:EllipBound}) with the boundary condition (\ref{eq:BoundCondFull}) along $\partial\mathscr{O}$. In particular, $u_{*}^{(X)}\in C_{\emph{loc}}(\bar{\mathscr{O}})\cap C^{2,\alpha}(\mathscr{O})$ which obeys (\ref{eq:lineargrowth}) on $\bar{\mathscr{O}}$.
\end{itemize}
\end{theorem}
Moreover, we obtain the existence of solutions to the elliptic boundary value problem when the boundary data $g$ is H\"{o}lder continuous on suitable portions of $\partial\mathscr{O}$.
\begin{theorem}\label{thm:ExistenceHolderCont}
In addition to the hypotheses of Theorem \ref{thm:UniqueEllipBound}, assume that the boundary portion $\Gamma_{1}$ be of class $C^{2,\alpha}$, that the coefficients $b,\sigma\in C^{0,\alpha}(\mathscr{O})$, and that $f\in C_{\emph{loc}}^{0,\alpha}(\mathscr{O}\cup\Gamma_{1})$, for some $\alpha\in(0,1)$.
\begin{itemize}
\item [(1)] Suppose that either case \emph{(a)}, \emph{(b)} or \emph{(c)} in Theorem \ref{thm:UniqueEllipBound} occurs. Assume that $g\in C_{\emph{loc}}^{2,\alpha}(\mathscr{O}\cup\Gamma_{1})$ which obeys the linear growth condition (\ref{eq:lineargrowth}) on $\mathscr{O}\cup\Gamma_{1}$. For any $x\in\mathscr{O}\cup\Gamma_{1}$, let $(\Omega,\mathscr{F},(\mathscr{F}_{s})_{s\geq 0},\mathbb{P}^{x},W,X)$ be a weak solution to (\ref{eq:MainDegenSDE1})-(\ref{eq:MainDegenSDE2}) with the initial condition (\ref{eq:MainDegenSDEIniCond}) at $t=0$, and let $u_{*}^{(X)}(x)$ be defined as in (\ref{eq:EllipBoundFeynmanKacTau}). Then, $u_{*}^{(X)}$ is a solution to (\ref{eq:EllipBound}) with the boundary condition (\ref{eq:BoundCondGamma1}) along $\Gamma_{1}$. In particular, $u_{*}^{(X)}\in C_{\emph{loc}}^{2,\alpha}(\mathscr{O}\cup\Gamma_{1})$ which obeys (\ref{eq:lineargrowth}) on $\mathscr{O}\cup\Gamma_{1}$.
\item [(2)] Suppose that either case \emph{(d)} or \emph{(e)} in Theorem \ref{thm:UniqueEllipBound} occurs. Assume that $g\in C_{\emph{loc}}^{2,\alpha}(\mathscr{O}\cup\Gamma_{1})\cup C_{\emph{loc}}(\partial\mathscr{O})$ which obeys the linear growth condition (\ref{eq:lineargrowth}) on $\bar{\mathscr{O}}$. For any $x\in\bar{\mathscr{O}}$, let $(\Omega,\mathscr{F},(\mathscr{F}_{s})_{s\geq 0},\mathbb{P}^{x},W,X)$ be a weak solution to (\ref{eq:MainDegenSDE1})-(\ref{eq:MainDegenSDE2}) with the initial condition (\ref{eq:MainDegenSDEIniCond}) at $t=0$, and let $u_{*}^{(X)}(x)$ be defined as in (\ref{eq:EllipBoundFeynmanKacTau}). Then, $u_{*}^{(X)}$ is a solution to (\ref{eq:EllipBound}) with the boundary condition (\ref{eq:BoundCondFull}) on $\partial\mathscr{O}$. In particular, $u_{*}^{(X)}\!\in C_{\emph{loc}}(\bar{\mathscr{O}})\cap C^{2,\alpha}_{\emph{loc}}(\mathscr{O}\cup\Gamma_{1})$ which obeys (\ref{eq:lineargrowth}) on $\bar{\mathscr{O}}$.
\end{itemize}
\end{theorem}
\begin{remark}\label{rem:ExistUnique}
Let $u_{*}$ be the \emph{unique} function (regardless of the choice of the weak solution) defined by (\ref{eq:EllipBoundFeynmanKacTau}). Theorem \ref{thm:UniqueEllipBound} and Theorem \ref{thm:ExistenceEllipBound} (or Theorem \ref{thm:ExistenceHolderCont}) implies that, in all cases \emph{(a)}-\emph{(e)}, $u_{*}$ is the \emph{unique} solution, which is $C^{2}$ (or $C^{2,\alpha}$) inside the domain and continuous up to the (partial or full) boundary, to the elliptic boundary value problem (\ref{eq:EllipBound}) with the partial boundary condition (\ref{eq:BoundCondGamma1}), or with the full boundary condition (\ref{eq:BoundCondFull}).
\end{remark}
\begin{remark}\label{rem:ExistWeightSobolev}
In this article, we only consider the stochastic interpretations for {\it classical solutions} to the equations, for which we assume regularity conditions, such as the continuity and the H\"older continuity, on $f, g$ and the coefficient functions in the operator $\mathscr{A}$. Existence of classical solutions to (\ref{eq:EllipBound}), with partial boundary condition (\ref{eq:BoundCondGamma1}), in weighted Sobolev spaces, will be studied elsewhere. For the special case when $\mathscr A$ is the Heston operator, the existence of such solution was proved in~\cite[Theorem 1.18]{DaskalopoulosFeehan:2011}, and the weighted Sobolev and H\"{o}lder regularity was shown in~\cite[Theorem 1.11]{FeehanPop:2013} and~\cite{FeehanPop:2015}. We also refer to~\cite{FeehanPop:2015} for introduction to the literature on related degenerate elliptic and parabolic equations. The existence of solutions in weighted Sobolev spaces to the parabolic problem (\ref{eq:ParaTermBound}) for the Heston operator, with partial terminal/boundary condition (\ref{eq:TermBoundPart}), was studied in~\cite{DaskalopoulosFeehan:2011}.
\end{remark}

\subsection{Feynman-Kac Formulas for Solutions to Elliptic Obstacle Problems}

For any $x\in\bar{\mathscr{O}}$, let $(\Omega,\mathscr{F},(\mathscr{F}_{s})_{s\geq 0},\mathbb{P}^{x},W,X)$ be a weak solution to (\ref{eq:MainDegenSDE1})-(\ref{eq:MainDegenSDE2}) with the initial condition (\ref{eq:MainDegenSDEIniCond}) at $t=0$, and let $\mathscr{T}^{x,X}$ be the collection of all $(\mathscr{F}_{s})_{s\geq 0}$-stopping times. For any $\theta_{1},\,\theta_{2}\in\mathscr{T}^{x,X}$, we define
\begin{align}
J^{\theta_{1},\theta_{2}}_{X}(x)&:=\mathbb{E}^{x}\!\left[\exp\!\left(-\!\int_{0}^{\theta_{1}}c(X_{s})\,ds\right)\!\psi(X_{\theta_{2}}){\bf 1}_{\{\theta_{1}>\theta_{2}\}}\right]+\mathbb{E}^{x}\!\left[\exp\!\left(-\!\int_{0}^{\theta_{1}}c(X_{s})\,ds\right)\!g(X_{\theta_{1}}){\bf 1}_{\{\theta_{1}\leq\theta_{2}\}}\right]\nonumber\\
\label{eq:FunJTheta12} &\quad\,\,+\mathbb{E}^{x}\!\left[\int_{0}^{\theta_{1}\wedge\theta_{2}}\exp\left(-\!\int_{0}^{s}c(X_{u})\,du\right)f(X_{s})\,ds\right],
\end{align}
and
\begin{align}\label{eq:EllipObsFeynmanKacTau}
v_{*}^{(X)}(x)&:=\sup_{\theta\in{ \mathscr{T}^{x,X}}}J^{\tau_{\mathscr{O}},\theta}_{X}(x),\\
\label{eq:EllipObsFeynmanKacLambda} v_{**}^{(X)}(x)&:=\sup_{\theta\in{ \mathscr{T}^{x,X}}}J^{\lambda_{\mathscr{O}},\theta}_{X}(x).
\end{align}
where $\tau_{\mathscr{O}}=\tau_{\mathscr{O}}^{x,X}$ and $\lambda_{\mathscr{O}}=\lambda_{\mathscr{O}}^{x,X}$ are defined by (\ref{eq:EllipBoundFeynmanKacTau}) and (\ref{eq:EllipBoundFeynmanKacLambda}), respectively.

We then have the following uniqueness results for the elliptic obstacle problem. Similar to the results for the elliptic boundary-value problem, if the diffusion cannot reach the boundary portion $\Gamma_{0}$, (\ref{eq:EllipObsFeynmanKacTau}) provides the unique solution to (\ref{eq:EllipObs}) with partial boundary condition, while if the diffusion can hit $\Gamma_{0}$, (\ref{eq:EllipObsFeynmanKacTau}) and (\ref{eq:EllipObsFeynmanKacLambda}) provide the unique solutions (with distinct regularities) to (\ref{eq:EllipObs}) with full and partial boundary conditions, respectively.
\begin{theorem}\label{thm:UniqueEllipObs}
Let $b$, $\sigma$, $c$ and $f$ be as in Theorem \ref{thm:UniqueEllipBound}, and let $\psi\in C(\mathscr{O})$ which obeys the linear growth condition (\ref{eq:lineargrowth}).
\begin{itemize}
\item [(1)] Suppose that either case \emph{(a)}, \emph{(b)} or \emph{(c)} in Theorem \ref{thm:UniqueEllipBound} occurs. Assume that $\psi\in C_{\emph{loc}}(\mathscr{O}\cup\Gamma_{1})$, and that $g\in C_{\emph{loc}}(\Gamma_{1})$ which obeys (\ref{eq:CompatObsGamm1}) and (\ref{eq:lineargrowth}) on $\Gamma_{1}$. Let
    \begin{align*}
    u\in C_{\emph{loc}}(\mathscr{O}\cup\Gamma_{1})\cap C^{2}(\mathscr{O})
    \end{align*}
    be a solution to the elliptic obstacle problem (\ref{eq:EllipObs}) and (\ref{eq:BoundCondGamma1}), such that both $u$ and $\mathscr{A}u$ obey (\ref{eq:lineargrowth}) on $\mathscr{O}$. Then, for any $x\in\mathscr{O}\cup\Gamma_{1}$, $u(x)=v_{*}^{(X)}(x)$, for any weak solution $(\Omega,\mathscr{F},(\mathscr{F}_{s})_{s\geq 0},\mathbb{P}^{x},W,X)$ to (\ref{eq:MainDegenSDE1})-(\ref{eq:MainDegenSDE2}) with the initial condition (\ref{eq:MainDegenSDEIniCond}) at $t=0$, where $v_{*}^{(X)}$ is given by (\ref{eq:EllipObsFeynmanKacTau}).
\item [(2)] Suppose that either case \emph{(d)} or \emph{(e)} in Theorem \ref{thm:UniqueEllipBound} occurs. Assume that $\psi\in C_{\emph{loc}}(\bar{\mathscr{O}})$, and that $g\in C_{\emph{loc}}(\partial\mathscr{O})$ which obeys (\ref{eq:CompatObsFull}) and (\ref{eq:lineargrowth}) on $\partial\mathscr{O}$. Let
    \begin{align*}
    u\in C_{\emph{loc}}(\bar{\mathscr{O}})\cap C^{2}(\mathscr{O})
    \end{align*}
    be a solution to the elliptic obstacle problem (\ref{eq:EllipObs}) and (\ref{eq:BoundCondFull}), such that both $u$ and $\mathscr{A}u$ obey (\ref{eq:lineargrowth}) on $\mathscr{O}$. Then, for any $x\in\bar{\mathscr{O}}$, $u(x)=v_{*}^{(X)}(x)$, for any arbitrary weak solution $(\Omega,\mathscr{F},(\mathscr{F}_{s})_{s\geq 0},\mathbb{P}^{x},W,X)$ to (\ref{eq:MainDegenSDE1})-(\ref{eq:MainDegenSDE2}) with the initial condition (\ref{eq:MainDegenSDEIniCond}) at $t=0$, where $v_{*}^{(X)}$ is given by (\ref{eq:EllipObsFeynmanKacTau}).
\end{itemize}
\end{theorem}
\begin{theorem}\label{thm:UniqueEllipObsPartial}
Suppose that either case \emph{(d)} or \emph{(e)} in Theorem \ref{thm:UniqueEllipBound} occurs. Let $b$, $\sigma$, $c$ and $f$ be as in Theorem \ref{thm:UniqueEllipBoundPartial}. Assume that $\psi\in C_{\emph{loc}}(\bar{\mathscr{O}})$ which obeys the linear growth condition (\ref{eq:lineargrowth}) on $\mathscr{O}$, and that $g\in C_{\emph{loc}}(\Gamma_{1})$ which obeys (\ref{eq:CompatObsGamm1}) and (\ref{eq:lineargrowth}) on $\Gamma_{1}$. Let
\begin{align*}
u\in C_{\emph{loc}}(\mathscr{O}\cup\Gamma_{1})\cap C^{2}(\mathscr{O})\cap C_{\emph{s,loc}}^{1,1,\beta}(\mathscr{O}\cup\Gamma_{0})
\end{align*}
be a solution to the elliptic obstacle problem (\ref{eq:EllipObs}) and (\ref{eq:BoundCondGamma1}), such that both $u$ and $\mathscr{A}u$ obey (\ref{eq:lineargrowth}). Then, for any $x\in\mathscr{O}\cup\Gamma_{1}$, $u(x)=v_{**}^{(X)}(x)$, for any weak solution $(\Omega,\mathscr{F},(\mathscr{F}_{s})_{s\geq 0},\mathbb{P}^{x},W,X)$ to (\ref{eq:MainDegenSDE1})-(\ref{eq:MainDegenSDE2}) with the initial condition (\ref{eq:MainDegenSDEIniCond}) at $t=0$, where $v_{**}^{(X)}$ is given by (\ref{eq:EllipObsFeynmanKacLambda}).
\end{theorem}
\begin{remark}\label{rem:ExistWeightSobolevObs}
Existence of solutions in weighted sobolev spaces to the elliptic obstacle problem (\ref{eq:EllipObs}) for the Heston operator, with the partial boundary condition (\ref{eq:BoundCondGamma1}), was proved in~\cite[Theorem 1.6]{DaskalopoulosFeehan:2011}, and the H\"{o}lder continuity of such solutions up to the boundary portion $\Gamma_{0}$ was proved in~\cite{FeehanPop:2013}.
\end{remark}

\section{Preliminary Results}\label{sec:Prelim}

\subsection{Properties of Diffusion Processes}

In this subsection, we consider a general $d$-dimensional ($d\geq 2$) time-homogeneous SDE system
\begin{align}\label{eq:GeneralSDE}
\left\{\begin{array}{ll} dX_{t}&\!\!\!=\,b(X_{t})\,dt+\sigma(X_{t})\,dW_{t}, \\ X_{0}&\!\!\!=\,x.\end{array} \right.
\end{align}
Throughout this section, we assume that the coefficients $b=(b_{1},\ldots,b_{d})^{T}$ and $\sigma=(\sigma_{ij})$, $1\leq i\leq d$, $1\leq j\leq m$, obey Assumption \ref{assupt:Cont} and Assumption \ref{assupt:LinearGrowth}. Then for any $x\in\mathbb{R}^{d}$, there exists a weak solution $(\Omega,\mathscr{F},(\mathscr{F}_{t})_{t\geq 0},\mathbb{P}^{x},W,X)$ to (\ref{eq:GeneralSDE}) (cf.~\cite[Theorem 4.3]{FeehanPop:2013(1)}). We also assume that $c$ is a non-negative function satisfying Assumption \ref{assupt:UnifPostc}.
\begin{lemma}\label{lem:SupMart}
Suppose there exists $p>0$ such that
\begin{equation}\label{eq:SupMartCon}
\liminf_{\|x\|\rightarrow\infty}\left[c(x)\|x\|^{2}-\left(p\|x\|^{2}+\frac{1}{p}\|b(x)\|^{2}+\|\sigma(x)\|^{2}\right)\right]>-\infty,
\end{equation}
then there exists $M\in[0,\infty)$ such that
\begin{align*}
Z_{t}:=\exp\left\{-\int_{0}^{t}c(X_{s})\,ds\right\}\|X_{t}\|^{2}+Mc_{0}^{-1}e^{-c_{0}t},\quad t\geq 0,
\end{align*}
is a supermartingale with respect to $(\mathscr{F}_{t})_{t\geq 0}$ under $\mathbb{P}^{x}$.
\end{lemma}
\noindent
\textbf{Proof:} For any $t>0$ and $M\in\mathbb{R}$, It\^{o}'s formula implies that
\begin{align}
dZ_{t}&=\exp\left(-\int_{0}^{t}c(X_{s})ds\right)\left(-c(X_{t})\|X_{t}\|^{2}+2\sum_{i=1}^{d}X^{(i)}_{t}b_{i}(X_{t})+\|\sigma(X_{t})\|^{2}\right)dt-Me^{-c_{0}t}\,dt\nonumber\\
\label{eq:SupMartDriftVol} &\quad+\exp\left(-\int_{0}^{t}c(X_{s})\,ds\right)\sum_{j=1}^{m}\left(\sum_{i=1}^{d}X^{(i)}_{t}\sigma_{ij}(X_{t})\right)dW^{(j)}_{t}.
\end{align}
To show that $Z_{t}$ is a supermartingale, since $Z_{t}$ is non-negative and the stochastic integral term in (\ref{eq:SupMartDriftVol}) is a local martingale, it suffices to show that the drift term is non-positive by Fatou's lemma. Indeed, for $p>0$ given in (\ref{eq:SupMartCon}),
\begin{align*}
\left|2\sum_{i=1}^{d}X^{(i)}_{t}b_{i}(X_{t})+\|\sigma(X_{t})\|^{2}\right|&\leq p\|X_{t}\|^{2}+\frac1p\|b(X_{t})\|^{2}+\|\sigma(X_{t})\|^{2}.
\end{align*}
By (\ref{eq:SupMartCon}), we may set
\begin{align}\label{eq:SupMartM}
M:=-\inf\limits_{x\in \mathbb R^d}\left[c(x)\|x\|^2-\left(p\|x\|^2+\dfrac1p\|b(x)\|^2+\|\sigma(x)\|^2\right)\right]\in[0,\infty),
\end{align}
and we obtain that
\begin{align*}
&\quad\exp\left(-\int_{0}^{t}c(X_{s})\,ds\right)\left(-c(X_{t})\|X_{t}\|^{2}+2\sum_{i=1}^{d}X^{(i)}_{t}b_{i}(X_{t})+\|\sigma(X_{t})\|^{2}\right)\\
&\leq\exp\left(-\int_{0}^{t}c(X_{s})\,ds\right)\left(-c(X_{t})\|X_{t}\|^{2}+p\|X_{t}\|^{2}+\frac1p\|b(X_{t})\|^{2}+\|\sigma(X_{t})\|^{2}\right)\leq Me^{-c_{0}t},
\end{align*}
which completes the proof. \hfill $\Box$
\begin{remark}\label{rem:CondSupMart}
\begin{itemize}
\item Condition (\ref{eq:SupMartCon}) holds when $c(x)>c_{0}>0$ for all $x\in\mathbb{R}^{d}$, and $\|b(x)\|+\|\sigma(x)\|\leq K(1+\|x\|^\beta)$ for $\beta\in (0,1)$. Actually, we can choose $0<p<c_{0}$ and consequently
    \begin{align*}
    \liminf_{\|x\|\rightarrow\infty}\left[c(x)\|x\|^{2}-\left(p\|x\|^{2}+\frac{1}{p}\|b(x)\|^{2}+\|\sigma(x)\|^{2}\right)\right]=\infty.
    \end{align*}
\item When $\|b(x)\|+\|\sigma(x)\|\leq K(1+\|x\|)$,  we may require that $c_{0}>0$ is sufficiently large (depending on $K$) to satisfy condition (\ref{eq:SupMartCon}).
\end{itemize}
\end{remark}
\begin{corollary}\label{cor:SupMartLastEle}
For $\mathbb{P}^{x}$-a.e. $\omega\in\Omega$,
\begin{align}
\lim_{t\rightarrow\infty}Z_{t}(\omega)=Z_{\infty}(\omega)\quad\text{exists},
\end{align}
and $(Z_{t})_{t\in[0,\infty]}$ is a supermartingale with respect to $(\mathscr{F}_{t})_{t\geq 0}$ under $\mathbb{P}^{x}$.
\end{corollary}
\noindent
\textbf{Proof:} The proof follows from \cite[Problem 1.3.16]{KaratzasShreve:1991}.\hfill $\Box$

\medskip
As a consequence of Lemma \ref{lem:SupMart} and Corollary \ref{cor:SupMartLastEle}, we obtain the following lemma, which, in particular, implies that the functions $u_{*}^{(X)}(x)$, $u_{**}^{(X)}(x)$, $v_{*}^{(X)}(x)$ and $v_{**}^{(X)}(x)$, given respectively by (\ref{eq:EllipBoundFeynmanKacTau}), (\ref{eq:EllipBoundFeynmanKacLambda}), (\ref{eq:EllipObsFeynmanKacTau}) and (\ref{eq:EllipObsFeynmanKacLambda}), are well defined and obey the linear growth condition (\ref{eq:lineargrowth}), for any $x\in\bar{\mathscr{O}}$ and any weak solution $(\Omega,\mathscr{F},(\mathscr{F}_{s})_{s\geq 0},\mathbb{P}^{x},W,X)$.
\begin{lemma}\label{lem:EstFunJ}
Let $f$, $g$ and $\psi$ are real-valued Borel measurable functions on $\mathbb{R}^{d}$ satisfying the linear growth condition (\ref{eq:lineargrowth}). Assume that the coefficients functions $b$, $\sigma$ and $c$ satisfy Assumption 2.2, Assumption 2.3, Assumption 2.5 and (\ref{eq:SupMartCon}). Then, for any $x\in\bar{\mathscr{O}}$, any weak solution $(\Omega,\mathscr{F},(\mathscr{F}_{s})_{s\geq 0},\mathbb{P}^{x},W,X)$ to (\ref{eq:MainDegenSDE1})-(\ref{eq:MainDegenSDE2}) with the initial condition (\ref{eq:MainDegenSDEIniCond}) at $t=0$, and any $\theta_{1},\theta_{2}\in\mathscr{T}^{x,X}$, the function $J^{\theta_{1},\theta_{2}}_{X}(x)$, given by (\ref{eq:FunJTheta12}), is well defined and satisfies
\begin{align}\label{eq:EstJ12}
\left|J^{\theta_{1},\theta_{2}}_{X}(x)\right|\leq C\left(1+\|x\|\right),
\end{align}
where $C$ is a universal positive constant, depending only on $K$ as in (\ref{eq:lineargrowth}), $c_{0}$ as in Assumption \ref{assupt:UnifPostc}, and $M$ as in (\ref{eq:SupMartM}).
\end{lemma}
\noindent
\textbf{Proof:} To show the well-definedness of { $J^{\theta_{1},\theta_{2}}_{X}(x)$, $x\in\bar{\mathscr{O}}$}, we need to show that the first integral term in (\ref{eq:FunJTheta12}) satisfies the following identity
\begin{align}\label{eq:SecondFunJTheta12FinTheta1}
\mathbb{E}^{x}\!\left[\exp\left(-\!\int_{0}^{\theta_{1}}\!\!c(X_{s})\,ds\right)g(X_{\theta_{1}}){\bf 1}_{\{\theta_{1}\leq\theta_{2}\}}\right]\!=\mathbb{E}^{x}\!\left[\exp\left(-\!\int_{0}^{\theta_{1}}\!\!c(X_{s})\,ds\right)g(X_{\theta_{1}}){\bf 1}_{\{\theta_{1}\leq\theta_{2},\,\theta_{1}<\infty\}}\right],
\end{align}
since $X=(X_{s})_{s\geq 0}$ does not necessarily have a limit at infinity. Indeed, for any $T>0$, the left-hand side of (\ref{eq:SecondFunJTheta12FinTheta1}) can be rewritten as
\begin{align*}
\mathbb{E}^{x}\left[\exp\left(-\int_{0}^{\theta_{1}}c(X_{s})\,ds\right)g(X_{\theta_{1}}){\bf 1}_{\{\theta_{1}\leq\theta_{2}\}}\right]&=\mathbb{E}^{x}\left[\exp\left(-\int_{0}^{\theta_{1}}c(X_{s})\,ds\right)g(X_{\theta_{1}}){\bf 1}_{\{\theta_{1}\leq\theta_{2}\wedge T\}}\right]\\
&\quad\,+\mathbb{E}^{x}\left[\exp\left(-\int_{0}^{\theta_{1}}c(X_{s})\,ds\right)g(X_{\theta_{1}}){\bf 1}_{\{T<\theta_{1}\leq\theta_{2}\}}\right].
\end{align*}
We will show that the second term on the right-hand side above vanishes as $T\rightarrow\infty$, regardless of the choices of random variables for $X_{\infty}$ so that $Z_{\infty}$ in Corollary \ref{cor:SupMartLastEle} holds. Using the linear growth condition (\ref{eq:lineargrowth}) on $g$, we have
\begin{align*}
&\quad\,\,\mathbb{E}^{x}\left[\exp\left(-\int_{0}^{\theta_{1}}c(X_{s})\,ds\right)\left|g(X_{\theta_{1}})\right|{\bf 1}_{\{T<\theta_{1}\leq\theta_{2}\}}\right]\\
&\leq\mathbb{E}^{x}\left[\exp\left(-\int_{0}^{\theta_{1}}c(X_{s})\,ds\right)\cdot K\left(1+\left\|X_{\theta_{1}}\right\|\right){\bf 1}_{\{T<\theta_{1}\}}\right]\\
&\leq Ke^{-c_{0}T}+K\left\{\left[\mathbb{E}^{x}\!\left(\exp\!\left(-\!\int_{0}^{\theta_{1}}\!c(X_{s})\,ds\right){\bf 1}_{\{T<\theta_{1}\}}\right)\right]^{1/2}\!\left[\mathbb{E}^{x}\left(\exp\!\left(-\!\int_{0}^{\theta_{1}}\!c(X_{s})\,ds\right)\left\|X_{\theta_{1}}\right\|^{2}\right)\right]^{1/2}\right\}\\
&\leq Ke^{-c_{0}T}+Ke^{-c_{0}T/2}\left\{\mathbb{E}^{x}\left[\exp\left(-\int_{0}^{\theta_{1}}c(X_{s})\,ds\right)\left\|X_{\theta_{1}}\right\|^{2}\right]\right\}^{1/2}.
\end{align*}
By Corollary \ref{cor:SupMartLastEle} and the Optional Sampling Theorem (cf.~\cite[Theorem 1.3.22]{KaratzasShreve:1991}),
\begin{align}\label{eq:SupMartZOptStop}
\mathbb{E}\left(Z_{\theta_{1}}\right)=\mathbb{E}\left[\exp\left(-\int_{0}^{\theta_{1}}c(X_{s})\,ds\right)\|X_{\theta_{1}}\|^{2}+Mc_{0}^{-1}e^{-c_{0}\theta_{1}}\right]\leq \|x\|^{2}+Mc_{0}^{-1}.
\end{align}
Therefore, we have
\begin{align*}
\mathbb{E}^{x}\left[\exp\left(-\int_{0}^{\theta_{1}}c(X_{s})\,ds\right)\left|g(X_{\theta_{1}})\right|{\bf 1}_{\{T<\theta_{1}\leq\theta_{2}\}}\right]\leq Ke^{-c_{0}T}+Ke^{-c_{0}T/2}\left(\|x\|^{2}+Mc_{0}^{-1}\right)\rightarrow 0,
\end{align*}
as $T\rightarrow\infty$, which justifies the identity (\ref{eq:SecondFunJTheta12FinTheta1}).

To obtain the estimation (\ref{eq:EstJ12}), we first analyze the integral term in (\ref{eq:FunJTheta12}). By the linear growth condition (\ref{eq:lineargrowth}),
\begin{align*}
&\quad\,\,\mathbb{E}^{x}\left[\int_{0}^{\theta_{1}\wedge\theta_{2}}\exp\left(-\int_{0}^{s}c(X_{u})\,du\right)\left|f(X_{s})\right|\,ds\right]\\
&\leq K\,\mathbb{E}^{x}\left[\int_{0}^{\infty}\exp\left(-\int_{0}^{s}c(X_{u})\,du\right)\left(1+\|X_{s}\|\right)ds\right]\\
&\leq K\int_{0}^{\infty}e^{-c_{0}s}\,ds+K\int_{0}^{\infty}e^{-c_{0}s/2}\left\{\mathbb{E}^{x}\left[\exp\left(-\int_{0}^{s}c(X_{u})\,du\right)\|X_{s}\|^{2}\right]\right\}^{1/2}ds\\
&=\frac{K}{c_{0}}+K\int_{0}^{\infty}e^{-c_{0}s/2}\left\{\mathbb{E}^{x}\left[\exp\left(-\int_{0}^{s}c(X_{u})\,du\right)\|X_{s}\|^{2}\right]\right\}^{1/2}ds.
\end{align*}
By Lemma \ref{lem:SupMart}, we have
\begin{align}
\mathbb{E}^{x}\left[\int_{0}^{\theta_{1}\wedge\theta_{2}}\exp\left(-\!\int_{0}^{s}c(X_{u})\,du\right)\left|f(X_{s})\right|ds\right]&\leq\mathbb{E}^{x}\left[\int_{0}^{\infty}\exp\left(-\!\int_{0}^{s}c(X_{u})\,du\right)\left|f(X_{s})\right|ds\right]\nonumber\\
&\leq K\left[\frac{1}{c_{0}}+\int_{0}^{\infty}e^{-c_{0}s/2}\left(\|x\|^{2}\!+\!\frac{M}{c_{0}}(1\!-\!e^{-c_{0}s})\right)^{1/2}\!ds\right]\nonumber\\
&\leq Kc_{0}^{-1}\left[1+2\left(\|x\|^{2}+c_{0}^{-1}M\right)^{1/2}\right]\nonumber\\
\label{eq:EstFirstTermJ12} &\leq 2Kc_{0}^{-1}\|x\|+Kc_{0}^{-1}\left(1+2\sqrt{c_{0}^{-1}M}\right).
\end{align}

Next, we estimate those two non-integral terms in (\ref{eq:FunJTheta12}). We will only provide the proof for the first non-integral term, since both terms have the same form. To do this, using the linear growth condition (\ref{eq:lineargrowth}) on $g$, the identity (\ref{eq:SecondFunJTheta12FinTheta1}), and a similar argument as above, we have
\begin{align*}
&\quad\,\,\mathbb{E}^{x}\left[\exp\left(-\int_{0}^{\theta_{1}}c(X_{s})\,ds\right)\left|g(X_{\theta_{1}})\right|{\bf 1}_{\{\theta_{1}\leq\theta_{2}\}}\right]\\
&=\mathbb{E}^{x}\left[\exp\left(-\int_{0}^{\theta_{1}}c(X_{s})\,ds\right)\left|g(X_{\theta_{1}})\right|{\bf 1}_{\{\theta_{1}\leq\theta_{2},\theta_{1}<\infty\}}\right]\\
&\leq\mathbb{E}^{x}\left[\exp\left(-\int_{0}^{\theta_{1}}c(X_{s})\,ds\right)\cdot K\left(1+\left\|X_{\theta_{1}}\right\|\right){\bf 1}_{\{\theta_{1}<\infty\}}\right]\\
&\leq K+K\left\{\left[\mathbb{E}^{x}\!\left(\exp\left(-\int_{0}^{\theta_{1}}c(X_{s})\,ds\right)\right)\right]^{1/2}\left[\mathbb{E}^{x}\!\left(\exp\left(-\int_{0}^{\theta_{1}}c(X_{s})\,ds\right)\left\|X_{\theta_{1}}\right\|^{2}{\bf 1}_{\{\theta_{1}<\infty\}}\right)\right]^{1/2}\right\}\\
&\leq K+K\left\{\mathbb{E}^{x}\!\left[\exp\left(-\int_{0}^{\theta_{1}}c(X_{s})\,ds\right)\left\|X_{\theta_{1}}\right\|^{2}{\bf 1}_{\{\theta_{1}<\infty\}}\right]\right\}^{1/2}.
\end{align*}
Therefore, by (\ref{eq:SupMartZOptStop}), we see that
\begin{align}\label{eq:EstSecondTermJ12}
\mathbb{E}^{x}\!\!\left[\exp\!\left(\!-\!\int_{0}^{\theta_{1}}\!\!c(X_{s})ds\!\right)\!\left|g(X_{\theta_{1}})\right|\!{\bf 1}_{\{\theta_{1}\leq\theta_{2}\}}\!\right]\!\leq\!K\!+\!K\!\left(\|x\|^{2}\!+\!Mc_{0}^{-1}\right)^{1/2}\!\leq\! K\!\left(\!1\!+\!\sqrt{c_{0}^{-1}M}\!+\!\|x\|\!\right).
\end{align}
Similarly, using the linear growth condition (\ref{eq:lineargrowth}) on $\psi$, we also have
\begin{align}\label{eq:EstThirdTermJ12}
\mathbb{E}^{x}\left[\exp\left(-\int_{0}^{\theta_{2}}c(X_{s})\,ds\right)\left|\psi(X_{\theta_{2}})\right|{\bf 1}_{\{\theta_{1}\geq\theta_{2}\}}\right]\leq K\left(1+\sqrt{c_{0}^{-1}M}+\|x\|\right).
\end{align}
Combining (\ref{eq:EstFirstTermJ12}), (\ref{eq:EstSecondTermJ12}) and (\ref{eq:EstThirdTermJ12}) completes the proof.\hfill $\Box$

\subsection{Preliminary on Boundary Classifications for One-Dimensional Diffusions}\label{subsec:BoundClass}

In this subsection, let us recall some basic results on boundary classifications for one-dimensional diffusion processes (see~\cite[Section 15.6]{KarlinTaylor:1981} for more details). Notice that the last SDE (\ref{eq:MainDegenSDE2}) in our model is independent of the first $d-1$ variables. The results in this subsection can therefore be applied to the one-dimensional process $X^{(d)}$.

Let $-\infty\leq l<r\leq\infty$. Throughout this subsection, let $(Y_{t})_{t\geq 0}$ be a one-dimensional diffusion process, defined on a complete filtered probability space $(\Omega,\mathscr{F},(\mathscr{F}_{t})_{t\geq 0},\mathbb{P})$, which satisfies
\begin{align}\label{eq:Gen1DDifProc}
Y_{t}=y+\int_{0}^{t}\mu(Y_{s})\,ds+\int_{0}^{t}\eta(Y_{s})\,dW_{s},\quad t\geq 0,
\end{align}
where $y\in(l,r)$. Assume that $\mu(\cdot)$ and $\eta(\cdot)$ are continuous on $(l,r)$, and that $\eta^2(\cdot)>0$ on $(l,r)$. Let $T_{z}^{y}$ be the first hitting time of $Y$ to $z$, starting at $Y_{0}=y$, i.e.,
\begin{align*}
T_{z}^{y}:=\inf\left\{t\geq 0:\,Y_{t}=z,\,Y_{0}=y\right\}.
\end{align*}
The upper index $y$ will be omitted when there is no ambiguity. We introduce the following notations, assuming $l<a<y<b<r$:
\begin{align*}
w_{a,b}(y):=\mathbb{P}^{y}\left(T_{b}<T_{a}\right),\quad v_{a,b}(y):=\mathbb{E}^{y}\left(T_{a}\wedge T_{b}\right).
\end{align*}
For any fixed $y_{0}\in(l,r)$, recall the \emph{scale measure} defined as
\begin{align}\label{eq:ScaleMea}
S[a,b]:=\int_{a}^{b}s(y)\,dy,\quad\text{where }\,s(y):=\exp\left\{-\int_{y_{0}}^{y}\frac{2\mu(x)}{\eta^2(x)}\,dx\right\},\quad l<a<b<r,
\end{align}
and the \emph{speed measure} defined as
\begin{align}\label{eq:SpeedMea}
M[a,b]:=\int_{a}^{b}m(y)\,dy,\quad\text{where }\,m(y):=\frac{1}{\eta^2(y)s(y)},\quad l<a<b<r.
\end{align}
At the left endpoint $l$, we define the scale measure and the speed measure as:
\begin{align}\label{eq:ScaleSpeedMeaLeftBound}
S(l,b]:=\lim_{a\downarrow l}S[a,b]\quad\text{ and }\quad M(l,b]:=\lim_{a\downarrow l}M[a,b].
\end{align}
In terms of the scale measure and the speed measure, $w_{a,b}$ and $v_{a,b}$ admit the representations:
\begin{align}\label{eq:ProbTbTaScale}
w_{a,b}(y)=\frac{S[a,y]}{S[a,b]},
\end{align}
and
\begin{align}\label{eq:ExpTaTbScaleSpeed}
v_{a,b}(y)=2\left(w_{a,b}(y)\int_{y}^{b}S[\xi,b]\,M(d\xi)+[1-w_{a,b}(y)]\int_{a}^{y}S[a,\xi]\,M(d\xi)\right).
\end{align}
Finally, for any $y\in(l,r)$, we define
\begin{align}
\label{eq:SigmaLeftBound} \Sigma(l,y)=\Sigma(l)&:=\int_{l}^{y}S(l,\xi]\,M(d\xi)=\int_{l}^{y}M[\xi,y]\,S(d\xi),\\
\label{eq:NLeftBound} N(l,y)=N(l)&:=\int_{l}^{y}S[\xi,y]\,M(d\xi)=\int_{l}^{y}M(l,\xi]\,S(d\xi).
\end{align}
Both definitions above depend on $y$. However, it can be shown that the finiteness of their values is independent of the choices of $y$ (cf.~\cite[Lemma 15.6.2]{KarlinTaylor:1981}). Therefore, we can suppress the dependence on $y$ without ambiguity since only the finiteness of their values is relevant in later arguments.

The following lemma summaries the classifications of the left boundary $l$, which will be useful in the proof of the existence and uniqueness theorems. We refer to~\cite[Table 15.6.2]{KarlinTaylor:1981} for more details.
\begin{lemma}\label{lem:BoundClass}
Let $b\in(l,r)$ be any interior point. The left boundary $l$ of the diffusion process (\ref{eq:Gen1DDifProc}) admits the following classifications.
\begin{itemize}
\item [(1)] If $S(l,b]<\infty$ and $M(l,b]<\infty$ (which implies that both $\Sigma(l)$ and $N(l)$ are finite, see~\cite[Lemma 15.6.3 (v)]{KarlinTaylor:1981}), in which case $l$ is a \emph{regular boundary}; or $M(l,b]=\infty$ and $\Sigma(l)<\infty$, in which case $l$ is an \emph{exit boundary}, we have
    \begin{align}\label{eq:HitTimeLeftBound}
    \lim_{y\downarrow l} T_{l}^{y}=0 \quad\text{a.}\,\text{s.}.
    \end{align}
\item [(2)] If $\Sigma(l)=\infty$ and $S(l,b]<\infty$, in which case $l$ is an \emph{attracting natural (Feller) boundary}, then the diffusion process $(Y(t))_{t\geq 0}$, starting from any interior point $y\in(l,r)$, can never achieve the left boundary $l$.
\item [(3)] If $S(l,b]=\infty$ (which implies that $\Sigma(l)=\infty$, see~\cite[Lemma 15.6.3 (i)]{KarlinTaylor:1981}), in which case $l$ is either a \emph{(non-attracting) natural (Feller) boundary} (when $N(l)=\infty$), or an \emph{entrance boundary} (when $N(l)<\infty$). Assume further that $\lim_{b\uparrow r}T_{b}^{y}=\infty$ almost surely, for any $y\in(l,r)$. Then, the diffusion process $(Y(t))_{t\geq 0}$, starting from any interior point $y\in(l,r)$, can never achieve the left boundary $l$.
\end{itemize}
\end{lemma}
\noindent
\textbf{Proof:} For any $y\in(l,r)$, define $T_{l+}^{y}:=\lim_{a\downarrow l}T_{a}^{y}$. Then
\begin{align}\label{eq:ContHitTimeY}
T_{l+}^{y}=T_{l}^{y},\quad\text{for every }\,y\in(l,r).
\end{align}
To see this, fix any $y\in(l,r)$, and since $T_{a}^{y}\leq T_{l}^{y}$ for any $a\in(l,y)$, we have
\begin{align*}
T_{l+}^{y}=\lim_{a\downarrow l}T_{a}^{y}\leq T_{l}^{y}.
\end{align*}
To show the reverse inequality, we can assume without loss of generality that $T_{l+}^{y}<\infty$ (otherwise $T_{l}^{y}=T_{l+}^{y}=\infty$). By the continuity of the sample paths of $(Y_{t})_{t\geq 0}$,
\begin{align*}
Y_{T_{l+}^{y}}=\lim_{a\downarrow l}Y_{T_{a}^{y}}=\lim_{a\downarrow l}a=l,
\end{align*}
and thus $T_{l+}^{y}\geq T_{l}^{y}$.

\medskip
\noindent
\textbf{(1)} Suppose first that $l$ is a regular boundary. For any fixed $b\in(l,r)$, by (\ref{eq:ScaleSpeedMeaLeftBound}), (\ref{eq:ProbTbTaScale}) and (\ref{eq:ContHitTimeY}),
\begin{align}\label{eq:ProbTleftBound}
\lim_{y\downarrow l}\mathbb{P}^{y}\left(T_{b}<T_{l}\right)=\lim_{y\downarrow l}\lim_{a\downarrow l}\mathbb{P}^{y}\left(T_{b}<T_{a}\right)=\lim_{y\downarrow l}\lim_{a\downarrow l}w_{a,b}(y)=\lim_{y\downarrow l}\lim_{a\downarrow l}\frac{S[a,y]}{S[a,b]}=\lim_{y\downarrow l}\frac{S(l,y]}{S(l,b]}=0.
\end{align}
Moreover, by~\cite[Lemma 15.6.3]{KarlinTaylor:1981}, $S(l,b]<\infty$ and $M(l,b]<\infty$ imply that $\Sigma(l)<\infty$ and $N(l)<\infty$. Hence, by (\ref{eq:ExpTaTbScaleSpeed}), (\ref{eq:ContHitTimeY}) and (\ref{eq:ProbTleftBound}),
\begin{align}
\lim_{y\downarrow l}\mathbb{E}^{y}(T_{l}\wedge T_{b})&=\lim_{y\downarrow l}\lim_{a\downarrow l}\mathbb{E}^{y}(T_{a}\wedge T_{b})=\lim_{y\downarrow l}\lim_{a\downarrow l}v_{a,b}(y)\nonumber\\
&=\lim_{y\downarrow l}\lim_{a\downarrow l}\left(2w_{a,b}(y)\int_{y}^{b}S[\xi,b]\,M(d\xi)+2(1-w_{a,b}(y))\int_{a}^{y}S[a,\xi]\,M(d\xi)\right)\nonumber\\
\label{eq:MeanHitLeftBound} &=2\lim_{y\downarrow l}\lim_{a\downarrow l}w_{a,b}(y)N(l)+2\lim_{y\downarrow l}\lim_{a\downarrow l}(1-w_{a,b}(y))\Sigma(l,y)=0.
\end{align}
Notice that
\begin{align*}
0=\lim_{y\downarrow l}\mathbb{E}^{y}\left(T_{l}\wedge T_{b}\right)\geq\lim_{y\downarrow l}\mathbb{E}^{y}\left(T_{l}{\bf 1}_{\{T_{l}\,\leq\,T_{b}\}}\right)\geq 0.
\end{align*}
Together with (\ref{eq:ProbTleftBound}), we concludes that $T_{l}^{y}$ converges to $0$ in $L^{1}$ (and thus in probability), as $y\downarrow l$. The almost-surely convergence follows immediately since $T_{l}^{y}$ is decreasing in $y$.

Next, assume that $l$ an exit boundary. By \cite[Lemma 15.6.3]{KarlinTaylor:1981}, $\Sigma(l)<\infty$ implies that $S(l,b]<\infty$. Hence, (\ref{eq:ProbTleftBound}) still holds. It suffices to show (\ref{eq:MeanHitLeftBound}) for an exit boundary $l$. As in the case of regular boundary, since $\Sigma(l)<\infty$, we have
\begin{align*}
\lim_{y\downarrow l}\lim_{a\downarrow l}2(1-w_{a,b}(y))\int_{a}^{y}S[a,\xi]\,M(d\xi)=\lim_{y\downarrow l}\lim_{a\downarrow l}(1-w_{a,b}(y))\Sigma(l,y)=0.
\end{align*}
It remains to show that
\begin{align*}
\lim_{y\downarrow l}\lim_{a\downarrow l}w_{a,b}(y)\int_{y}^{b}S[\xi,b]\,M(d\xi)=0.
\end{align*}
Notice that
\begin{align*}
0&\leq\lim_{y\downarrow l}\lim_{a\downarrow l}w_{a,b}(y)\int_{y}^{b}S[\xi,b]\,M(d\xi)=\lim_{y\downarrow l}\lim_{a\downarrow l}\frac{S[a,y]}{S[a,b]}\int_{y}^{b}S[\xi,b]\,M(d\xi)\\
&\leq\lim_{y\downarrow l}\lim_{a\downarrow l}S[a,y]M[y,b]=\lim_{y\downarrow l}S(l,y]M[y,b].
\end{align*}
We only need to show that
\begin{align}\label{eq:LimScaleSpeed}
\lim_{y\downarrow l}S(l,y]M[y,b]=0.
\end{align}
On the one hand, by the integration by parts formula,
\begin{align}
\int_{l}^{b}S(l,\xi]\,M(d\xi)&=\lim_{a\downarrow l}\int_{a}^{b}S(l,\xi]\,d(-M[\xi,b])=\lim_{a\downarrow l}\left[-S(l,\xi]M[\xi,b]\Big|_{a}^{b}+\int_{a}^{b}M[\xi,b]\,d\left(S(l,\xi]\right)\right]\nonumber\\
\label{eq:LimScaleSpeedIntbyPart} &=\lim_{a\downarrow l}S(l,a]M[a,b]+\int_{l}^{b}M[\xi,b]s(\xi)\,d\xi.
\end{align}
On the other hand, by Fubini's Theorem,
\begin{align}\label{eq:LimScaleSpeedFubini}
\int_{l}^{b}S(l,\xi]M(d\xi)=\!\int_{l}^{b}\!\left(\int_{l}^{\xi}s(\zeta)\,d\zeta\right)\!m(\xi)\,d\xi=\!\int_{l}^{b}\!\left(\int_{\zeta}^{b}m(\xi)\,d\xi\right)\!s(\zeta)\,d\zeta=\!\int_{l}^{b}\!M[\zeta,b]s(\zeta)\,d\zeta.
\end{align}
Therefore, (\ref{eq:LimScaleSpeed}) follows immediately from (\ref{eq:LimScaleSpeedIntbyPart}) and (\ref{eq:LimScaleSpeedFubini}), which concludes the proof of part (1).

\medskip
\noindent
\textbf{(2)} By~\cite[Lemma 15.6.2]{KarlinTaylor:1981}, if $l$ is attracting, i.e., $S(l,b]<\infty$, then
\begin{align}\label{eq:InfHitLeftBound}
\Sigma(l)=\infty\quad\Leftrightarrow\quad\mathbb{P}^{y}\left\{T_{l}<\infty\right\}=0,\quad\text{ for any }\,\,y\in(l,r).
\end{align}
Therefore, $T_{l}^{y}=\infty$ almost surely, for any $y\in(l,r)$, which completes the proof of part (2).

\medskip
\noindent
\textbf{(3)}
By~\cite[Lemma 15.6.1(ii)]{KarlinTaylor:1981}, if $S(l,b]=\infty$, then
\begin{align*}
\mathbb{P}^{y}\left\{T_{l+}<T_{b}\right\}=0,\quad\text{ for any }\,\,l<y<b<r.
\end{align*}
Hence, since $\lim_{b\uparrow r}T_{b}^{y}=\infty$, we have
\begin{align*}
\mathbb{P}^{y}\left\{T_{l+}<\infty\right\}=\lim_{b\uparrow r}\mathbb{P}^{y}\left\{T_{l+}<T_{b}\right\}=0,\quad\text{ for any }\,\,y\in(l,r),
\end{align*}
which completes the proof of the lemma.\hfill $\Box$
\begin{remark}\label{rem:BoundClass}
The above classifications can be summarized by the following graph.
\begin{equation*}
\begin{cases}
S(l,b]<\infty
\begin{cases}
M(l,b]<\infty \text{ regular boundary, \bf{case (1)}}\\ M(l,b]=\infty
\begin{cases}
\Sigma(l)=\infty \text{ natural (Feller) boundary, \bf{case (2)}}\\
\Sigma(l)<\infty \text{ exit boundary, \bf{case (1)}}
\end{cases}
\end{cases}\\
S(l,b]=\infty
\begin{cases}
N(l)<\infty \text{ entrance boundary }
\begin{cases}
\,\,\lim_{b\uparrow r}T_{b}^{y}=\infty\,\,\,\,a.s.\text{ \bf{ case (3)}}\\
\text{ otherwise, \bf{ no description}}
\end{cases}\\
N(l)=\infty \text{ natural (Feller) boundary }
\begin{cases}
\,\,\lim_{b\uparrow r}T_{b}^{y}=\infty\,\,\,\,a.s.\text{ \bf{ case (3)}}\\
\text{ otherwise, \bf{ no description}}
\end{cases}
\end{cases}
\end{cases}
\end{equation*}
For the diffusion process $X^{(d)}$ given in (\ref{eq:MainDegenSDE2}), $(l,r)=(0,+\infty)$ and thus $\lim_{b\uparrow \infty}T_{b}^{y}=\infty$ almost surely for any $y\in(0,\infty)$ since $X^{(d)}$ is almost surely finite at any finite time. Therefore, the boundary behavior of $X^{(d)}$ at the origin is completely covered by Lemma \ref{lem:BoundClass}.
\end{remark}

\section{Elliptic Boundary-Value Problems}\label{sec:EllipBoundValProb}

In this section, we will verify the existence and uniqueness of Feynman-Kac representations for solutions to the elliptic boundary value problem with partial/full boundary conditions. Hereafter, we use $S^{(d)}$, $M^{(d)}$, $\Sigma^{(d)}$ and $N^{(d)}$ to denote the quantities introduced in (\ref{eq:ScaleMea}), (\ref{eq:SpeedMea}), (\ref{eq:SigmaLeftBound}) and (\ref{eq:NLeftBound}), for any arbitrarily fixed diffusion process (weak solution) $X^{(d)}$ driven by (\ref{eq:MainDegenSDE2}). We will prove Theorem \ref{thm:UniqueEllipBound}, Theorem \ref{thm:UniqueEllipBoundPartial}, Theorem \ref{thm:ExistenceEllipBound} and Theorem \ref{thm:ExistenceHolderCont} in the following two scenarios, based on the boudary classifications of the left boundary $0$ for $X^{(d)}$.
\begin{itemize}
\item [(A)] $\Sigma^{(d)}(0)=\infty$, in which case the origin is either a natural (Feller) boundary or an entrance boundary of $X^{(d)}$, and thus $X^{(d)}$ can never achieve the origin from any interior point in $(0,\infty)$.
\item [(B)] $S^{(d)}(0,b]<\infty$ and $M^{(d)}(0,b]<\infty$, in which case the origin is a regular boundary for $X^{(d)}$; or $M^{(d)}(0,b]=\infty$ and $\Sigma^{(d)}(0)<\infty$, in which case the origin is an exit boundary for $X^{(d)}$. In both cases, $X^{(d)}$ can reach the origin from the interior $(0,\infty)$.
\end{itemize}
The first scenario covers the cases (2) and (3) in Lemma \ref{lem:BoundClass}, while the second scenario is the case (1) in Lemma \ref{lem:BoundClass}. By Remark \ref{rem:BoundClass}, these two scenarios describe all the possible boundary behavior of $X^{(d)}$ at the origin.

The following lemma shows the relations between the above scenarios and the cases (a)-(e) stated in Theorem \ref{thm:UniqueEllipBound}. In particular, the scenario (A) contains the cases (a), (b) and (c), while the scenario (B) contains the cases (d) and (e).
\begin{lemma}\label{lem:BoundClassXd}
Consider the one-dimensional diffusion process $X^{(d)}$ driven by (\ref{eq:MainDegenSDE2}).
\begin{itemize}
\item [(1)] If $\beta\in(0,1)$, $0$ is a regular boundary for $X^{(d)}$.
\item [(2)] If $\beta\in(1,2]$, assume that $b_{d}(0)>0$ and that $b_{d}$ is locally H\"{o}lder continuous at the origin (see (\ref{eq:LocalHolderbd})). Then, $0$ is either a (non-attracting) natural (Feller) boundary, or an entrance boundary for $X^{(d)}$.
\item [(3)] If $\beta=1$, assume that $b_{d}(0)>0$ and that $b_{d}$ is locally H\"{o}lder continuous at the origin.
    \begin{itemize}
    \item [(i)] If $2b_{d}(0)>\tilde{\sigma}_{0}^{2}(0)$, $0$ is either a (non-attracting) natural (Feller) boundary, or an \emph{entrance} boundary for $X^{(d)}$.
    \item [(ii)] If $2b_{d}(0)<\tilde{\sigma}_{0}^{2}(0)$, $0$ is a regular boundary for $X^{(d)}$.
    \item [(iii)] If $2b_{d}(0)=\tilde{\sigma}_{0}^{2}(0)$, assume further that $\tilde{\sigma}_{0}^{2}$ is constant in a neighborhood of the origin. Then, $0$ is either a (non-attracting) natural (Feller) boundary, or an entrance boundary for $X^{(d)}$.
    \end{itemize}
\end{itemize}
\end{lemma}
\noindent
\textbf{Proof:} Without loss of generality, we assume $y_{0}=1$ and $b\in(0,1)$ throughout this proof.

\medskip
\noindent
\textbf{(1)} By the continuity of $b_{d}$ and $\tilde{\sigma}_{0}^{2}$ (Assumption \ref{assupt:Cont}), there exists a constant $C_{1}>0$, such that
\begin{align*}
\left|b_{d}(y)\right|\leq C_{1},\quad\tilde{\sigma}_{0}^{2}(y)\leq C_{1},\quad\text{for any }\,y\in[0,1].
\end{align*}
Also, by Assumption \ref{assupt:Elliptildesigma},
\begin{align}\label{eq:LowerBoundsigma0}
\tilde{\sigma}_{0}^{2}(y)=\tilde{a}_{dd}(y)\geq\delta,\quad\text{for any }\,y\in[0,+\infty).
\end{align}
Hence, if $\beta\in(0,1)$,
\begin{align*}
S^{(d)}(0,b]&\leq\int_{0}^{b}\exp\left(\int_{y}^{1}\frac{2C_{1}}{\delta}\eta^{-\beta}\,d\eta\right)dy=\exp\left(\frac{2C_{1}}{\delta(1-\beta)}\right)\int_{0}^{b}\exp\left(-\frac{2C_{1}}{\delta(1-\beta)}y^{1-\beta}\right)dy<+\infty,\\
M^{(d)}(0,b]&\leq\int_{0}^{b}\frac{1}{\delta y^{\beta}}\exp\left(\int_{y}^{1}\frac{2C_{1}}{\delta\eta^{\beta}}\,d\eta\right)dy=\frac{1}{\delta}\exp\left(\frac{2C_{1}}{\delta(1\!-\!\beta)}\right)\int_{0}^{b}\frac{1}{y^{\beta}}\exp\left(-\frac{2C_{1}y^{1-\beta}}{\delta(1\!-\!\beta)}\right)dy<+\infty,
\end{align*}
which implies that $0$ is a regular boundary for $X^{(d)}$.

\medskip
\noindent
\textbf{(2)} When $\beta\in(1,2]$, assume without loss of generality that $b\in(0,\kappa]$, where $\kappa$ is given as in (\ref{eq:LocalHolderbd}). By (\ref{eq:LocalHolderbd}), we have
\begin{align*}
S^{(d)}(0,b]&=\int_{0}^{b}\exp\left(\int_{y}^{\kappa}\frac{2b_{d}(\eta)}{\eta^{\beta}\tilde{\sigma}_{0}^{2}(\eta)}\,d\eta+\int_{\kappa}^{1}\frac{2b_{d}(\eta)}{\eta^{\beta}\tilde{\sigma}_{0}^{2}(\eta)}\,d\eta\right)dy\\
&=C_{2}\int_{0}^{b}\exp\left(\int_{y}^{\kappa}\frac{2\left(b_{d}(\eta)-b_{d}(0)\right)}{\eta^{\beta}\tilde{\sigma}_{0}^{2}(\eta)}\,d\eta+\int_{y}^{\kappa}\frac{2b_{d}(0)}{\eta^{\beta}\tilde{\sigma}_{0}^{2}(\eta)}\,d\eta\right)dy\\
&\geq C_{2}\int_{0}^{b}\exp\left(-\frac{2L}{\delta}\int_{y}^{\kappa}\eta^{\gamma-\beta}\,d\eta+\frac{2b_{d}(0)}{C_{1}}\int_{y}^{\kappa}\eta^{-\beta}\,d\eta\right)dy,
\end{align*}
where $C_{2}:=\exp\left(\int_{\kappa}^{1}\frac{2b_{d}(\eta)}{\eta^{\beta}\tilde{\sigma}_{0}^{2}(\eta)}d\eta\right)$. we shall discuss the following three cases for different values $\gamma$, where $\gamma$ is the H\"{o}lder exponent given in (\ref{eq:LocalHolderbd}).

First, if $\gamma-\beta+1>0$,
\begin{align*}
S^{(d)}(0,b]&\geq C_{2}\exp\left(-\frac{2L\kappa^{\gamma-\beta+1}}{\delta(\gamma-\beta+1)}-\frac{2b_{d}(0)\kappa^{1-\beta}}{C_{1}(\beta-1)}\right)\int_{0}^{b}\exp\left(\frac{2Ly^{\gamma-\beta+1}}{\delta(\gamma-\beta+1)}+\frac{2b_{d}(0)y^{1-\beta}}{C_{1}(\beta-1)}\right)dy\\
&\geq C_{2}\exp\left(-\frac{2L\kappa^{\gamma-\beta+1}}{\delta(\gamma-\beta+1)}-\frac{2b_{d}(0)\kappa^{1-\beta}}{C_{1}(\beta-1)}\right)\int_{0}^{b}\exp\left(\frac{2b_{d}(0)}{C_{1}(\beta-1)}y^{1-\beta}\right)dy=+\infty.
\end{align*}
Next, if $\gamma-\beta+1<0$,
\begin{align*}
S^{(d)}(0,b]\geq C_{2}\exp\!\left(\frac{2L\kappa^{\gamma-\beta+1}}{\delta(\beta\!-\!\gamma\!-\!1)}\!-\!\frac{2b_{d}(0)\kappa^{1-\beta}}{C_{1}(\beta\!-\!1)}\right)\!\!\int_{0}^{b}\!\exp\!\left(\!\left[\frac{2b_{d}(0)y^{-\gamma}}{C_{1}(\beta\!-\!1)}\!-\!\frac{2L}{\delta(\beta\!-\!\gamma\!-\!1)}\right]\!y^{\gamma-\beta+1}\!\right)\!dy=+\infty,
\end{align*}
since
\begin{align*}
\frac{2b_{d}(0)y^{-\gamma}}{C_{1}(\beta-1)}>\frac{2L}{\delta(\beta-\gamma-1)}>0,\quad\text{for }\,y>0\,\,\,\text{small enough}.
\end{align*}
Finally, if $\gamma-\beta+1=0$, we also have
\begin{align*}
S^{(d)}(0,b]\geq C_{2}\kappa^{-\frac{2L}{\delta}}\exp\left(-\frac{2b_{d}(0)\kappa^{1-\beta}}{C_{1}(\beta-1)}\right)\int_{0}^{b}y^{\frac{2L}{\delta}}\exp\left(\frac{2b_{d}(0)}{C_{1}(\beta-1)}y^{1-\beta}\right)dy=+\infty.
\end{align*}
Therefore, when $\beta\in(1,2]$, we always have $S^{(d)}(0,b]=+\infty$, which implies that $0$ is either a (non-attracting) natural (Feller) boundary or an entrance boundary for $X^{(d)}$.

\medskip
\noindent
\textbf{(3)} Finally, we assume that $\beta=1$. Again, by the continuity of $b_{d}$ and $\tilde{\sigma}_{0}^{2}$ (Assumption \ref{assupt:Cont}), the assumption $b_{d}(0)>0$, and (\ref{eq:LowerBoundsigma0}), for any $\varepsilon\in\left(0,\min\left\{b_{d}(0),\tilde{\sigma}_{0}^{2}(0)\right\}\right)$, there exists $\kappa_{1}\in(0,\kappa)$ (where $\kappa$ is given as in (\ref{eq:LocalHolderbd})), such that, for any $0\leq y\leq\kappa_{1}$,
\begin{align}\label{eq:ContOribdsigma0}
\left|b_{d}(y)-b_{d}(0)\right|\leq\varepsilon\quad\text{and}\quad\left|\tilde{\sigma}_{0}^{2}(y)-\tilde{\sigma}_{0}^{2}(0)\right|\leq\varepsilon.
\end{align}
Without loss of generality, we can choose $b\in(0,\kappa_{1})$, then we have
\begin{align*}
S^{(d)}(0,b]&=\int_{0}^{b}\exp\left(\int_{y}^{\kappa_{1}}\frac{2b_{d}(\eta)}{\eta\,\tilde{\sigma}_{0}^{2}(\eta)}\,d\eta+\int_{\kappa_{1}}^{1}\frac{2b_{d}(\eta)}{\eta\,\tilde{\sigma}_{0}^{2}(\eta)}\,d\eta\right)dy\\
&=C_{3}\int_{0}^{b}\exp\left(\int_{y}^{\kappa_{1}}\frac{2\left(b_{d}(\eta)-b_{d}(0)\right)}{\eta\,\tilde{\sigma}_{0}^{2}(\eta)}\,d\eta+\int_{y}^{\kappa_{1}}\frac{2b_{d}(0)}{\eta\,\tilde{\sigma}_{0}^{2}(\eta)}\,d\eta\right)dy,
\end{align*}
where $C_{3}:=\exp\left\{\int_{\kappa_{1}}^{1}\frac{2b_{d}(\eta)}{\eta\,\tilde{\sigma}_{0}^{2}(\eta)}\,d\eta\right\}$.

\smallskip
\noindent
\textbf{(i)} If $2b_{d}(0)>\tilde{\sigma}_{0}^{2}(0)$, we may choose $\varepsilon>0$ so that $2b_{d}(0)>\tilde{\sigma}_{0}^{2}(0)+\varepsilon$. By (\ref{eq:LocalHolderbd}) and (\ref{eq:ContOribdsigma0}),
\begin{align*}
S^{(d)}(0,b]&\geq C_{3}\int_{0}^{b}\exp\left(-\frac{2L}{\tilde{\sigma}_{0}^{2}(0)-\varepsilon}\int_{y}^{\kappa_{1}}\eta^{\gamma-1}\,d\eta\right)\exp\left(\frac{2b_{d}(0)}{\tilde{\sigma}_{0}^{2}(0)+\varepsilon}\int_{y}^{\kappa_{1}}\eta^{-1}\,d\eta\right)dy\\
&=C_{3}\exp\left(-\frac{2L\kappa_{1}^{\gamma}}{\gamma\left(\tilde{\sigma}_{0}^{2}(0)-\varepsilon\right)}\right)\kappa_{1}^{\frac{2b_{d}(0)}{\tilde{\sigma}_{0}^{2}(0)+\varepsilon}}\int_{0}^{b}\exp\left(\frac{2Ly^{\gamma}}{\gamma\left(\tilde{\sigma}_{0}^{2}(0)-\varepsilon\right)}\right)y^{-\frac{2b_{d}(0)}{\tilde{\sigma}_{0}^{2}(0)+\varepsilon}}\,dy\\
&\geq C_{3}\exp\left(-\frac{2L\kappa_{1}^{\gamma}}{\gamma\left(\tilde{\sigma}_{0}^{2}(0)-\varepsilon\right)}\right)\kappa_{1}^{\frac{2b_{d}(0)}{\tilde{\sigma}_{0}^{2}(0)+\varepsilon}}\int_{0}^{b}y^{-\frac{2b_{d}(0)}{\tilde{\sigma}_{0}^{2}(0)+\varepsilon}}\,dy=+\infty.
\end{align*}
Therefore, $0$ is either a (non-attracting) natural (Feller) boundary or an entrance boundary for $X^{(d)}$.

\smallskip
\noindent
\textbf{(ii)} If $2b_{d}(0)<\tilde{\sigma}_{0}^{2}(0)$, we may choose $\varepsilon>0$ so that $2b_{d}(0)<\tilde{\sigma}_{0}^{2}(0)-\varepsilon$. Again by (\ref{eq:LocalHolderbd}) and (\ref{eq:ContOribdsigma0}),
\begin{align*}
S^{(d)}(0,b]&\leq C_{3}\int_{0}^{b}\exp\left(\frac{2L}{\tilde{\sigma}_{0}^{2}(0)-\varepsilon}\int_{y}^{\kappa_{1}}\eta^{\gamma-1}\,d\eta\right)\exp\left(\frac{2b_{d}(0)}{\tilde{\sigma}_{0}^{2}(0)-\varepsilon}\int_{y}^{\kappa_{1}}\eta^{-1}\,d\eta\right)dy\\
&=C_{3}\exp\left(\frac{2L\kappa_{1}^{\gamma}}{\gamma\left(\tilde{\sigma}_{0}^{2}(0)-\varepsilon\right)}\right)\kappa_{1}^{\frac{2b_{d}(0)}{\tilde{\sigma}_{0}^{2}(0)-\varepsilon}}\int_{0}^{b}\exp\left(-\frac{2Ly^{\gamma}}{\gamma\left(\tilde{\sigma}_{0}^{2}(0)-\varepsilon\right)}\right)y^{-\frac{2b_{d}(0)}{\tilde{\sigma}_{0}^{2}(0)-\varepsilon}}\,dy\\
&\leq C_{3}\exp\left(\frac{2L\kappa_{1}^{\gamma}}{\gamma\left(\tilde{\sigma}_{0}^{2}(0)-\varepsilon\right)}\right)\kappa_{1}^{\frac{2b_{d}(0)}{\tilde{\sigma}_{0}^{2}(0)-\varepsilon}}\int_{0}^{b}y^{-\frac{2b_{d}(0)}{\tilde{\sigma}_{0}^{2}(0)-\varepsilon}}\,dy<+\infty.
\end{align*}
Moreover, we have
\begin{align*}
M^{(d)}(0,b]&=\int_{0}^{b}\frac{1}{y\,\tilde{\sigma}_{0}^{2}(y)}\exp\left(-\int_{y}^{\kappa_{1}}\frac{2b_{d}(\eta)}{\eta\,\tilde{\sigma}_{0}^{2}(\eta)}\,d\eta-\int_{\kappa_{1}}^{1}\frac{2b_{d}(\eta)}{\eta\,\tilde{\sigma}_{0}^{2}(\eta)}\,d\eta\right)dy\\
&\leq\frac{C_{3}^{-1}}{\tilde{\sigma}_{0}^{2}(0)-\varepsilon}\int_{0}^{b}y^{-1}\exp\left(-\int_{y}^{\kappa_{1}}\frac{2\left(b_{d}(\eta)-b_{d}(0)\right)}{\eta\,\tilde{\sigma}_{0}^{2}(\eta)}\,d\eta-\int_{y}^{\kappa_{1}}\frac{2b_{d}(0)}{\eta\,\tilde{\sigma}_{0}^{2}(\eta)}\,d\eta\right)dy\\
&\leq\frac{C_{3}^{-1}}{\tilde{\sigma}_{0}^{2}(0)-\varepsilon}\int_{0}^{b}y^{-1}\exp\left(\frac{2L}{\tilde{\sigma}_{0}^{2}(0)-\varepsilon}\int_{y}^{\kappa_{1}}\eta^{\gamma-1}\,d\eta\right)\exp\left(-\frac{2b_{d}(0)}{\tilde{\sigma}_{0}^{2}(0)+\varepsilon}\int_{y}^{\kappa_{1}}\eta^{-1}\,d\eta\right)dy\\
&=\frac{C_{3}^{-1}}{\tilde{\sigma}_{0}^{2}(0)-\varepsilon}\exp\left(\frac{2L\kappa_{1}^{\gamma}}{\gamma\left(\tilde{\sigma}_{0}^{2}(0)-\varepsilon\right)}\right)\kappa_{1}^{-\frac{2b_{d}(0)}{\tilde{\sigma}_{0}^{2}(0)+\varepsilon}}\int_{0}^{b}y^{-1}\exp\left(-\frac{2Ly^{\gamma}}{\gamma\left(\tilde{\sigma}_{0}^{2}(0)-\varepsilon\right)}\right)y^{\frac{2b_{d}(0)}{\tilde{\sigma}_{0}^{2}(0)+\varepsilon}}\,dy\\
&\leq\frac{C_{3}^{-1}}{\tilde{\sigma}_{0}^{2}(0)-\varepsilon}\exp\left(\frac{2L\kappa_{1}^{\gamma}}{\gamma\left(\tilde{\sigma}_{0}^{2}(0)-\varepsilon\right)}\right)\kappa_{1}^{-\frac{2b_{d}(0)}{\tilde{\sigma}_{0}^{2}(0)+\varepsilon}}\int_{0}^{b}y^{\frac{2b_{d}(0)}{\tilde{\sigma}_{0}^{2}(0)+\varepsilon}-1}\,dy<+\infty.
\end{align*}
Therefore, $0$ is a regular boundary for $X^{(d)}$.

\smallskip
\noindent
\textbf{(iii)} If $2b_{d}(0)=\tilde{\sigma}_{0}^{2}(0)$, assume further that there exists $\kappa_{2}\in(0,\kappa)$ (where $\kappa>0$ is given in (\ref{eq:LocalHolderbd})), such that $\tilde{\sigma}_{0}^{2}(y)\equiv\tilde{\sigma}_{0}^{2}(0)$ for all $0\leq y\leq\kappa_{2}$. Choose any $b\in(0,\kappa_{2})$, by (\ref{eq:LocalHolderbd}), we have
\begin{align*}
S^{(d)}(0,b]&=\int_{0}^{b}\exp\left(\int_{y}^{\kappa_{2}}\frac{2b_{d}(\eta)}{\eta\,\tilde{\sigma}_{0}^{2}(\eta)}\,d\eta+\int_{\kappa_{2}}^{1}\frac{2b_{d}(\eta)}{\eta\,\tilde{\sigma}_{0}^{2}(\eta)}\,d\eta\right)dy\\
&=C_{4}\int_{0}^{b}\exp\left(\int_{y}^{\kappa_{2}}\frac{2\left(b_{d}(\eta)-b_{d}(0)\right)}{\eta\,\tilde{\sigma}_{0}^{2}(\eta)}\,d\eta+\int_{y}^{\kappa_{2}}\frac{2b_{d}(0)}{\eta\,\tilde{\sigma}_{0}^{2}(\eta)}\,d\eta\right)dy\\
&\geq C_{4}\int_{0}^{b}\exp\left(-\frac{2L}{\tilde{\sigma}_{0}^{2}(0)}\int_{y}^{\kappa_{2}}\eta^{\gamma-1}\,d\eta\right)\exp\left(\int_{y}^{\kappa_{2}}\eta^{-1}\,d\eta\right)dy\\
&=C_{4}\exp\left(-\frac{2L\kappa_{2}^{\gamma}}{\gamma\,\tilde{\sigma}_{0}^{2}(0)}\right)\kappa_{2}\int_{0}^{b}\exp\left(\frac{2Ly^{\gamma}}{\gamma\,\tilde{\sigma}_{0}^{2}(0)}\right)y^{-1}\,dy\\
&\geq C_{4}\exp\left(-\frac{2L\kappa_{2}^{\gamma}}{\gamma\,\tilde{\sigma}_{0}^{2}(0)}\right)\kappa_{2}\int_{0}^{b}y^{-1}\,dy=+\infty,
\end{align*}
where $C_{4}:=\exp\left(\int_{\kappa_{2}}^{1}\frac{2b_{d}(\eta)}{\eta\,\tilde{\sigma}_{0}^{2}(\eta)}\,d\eta\right)$. Therefore, $0$ is either a (non-attracting) natural (Feller) boundary or an entrance boundary for $X^{(d)}$. The proof is now complete. \hfill $\Box$

\medskip
We are now in the position of verifying the existence and uniqueness theorems for the elliptic boundary value problems. Notice that the existence and uniqueness of solutions to (\ref{eq:EllipBound}) with the partial boundary condition (\ref{eq:BoundCondGamma1}) along $\Gamma_{1}$, when the scenario (A) occurs, and with the full boundary condition (\ref{eq:BoundCondFull}) along $\partial\mathscr{O}$, when the scenario (B) occurs, are similar in nature. Therefore, we define
\begin{align}\label{eq:DefPartialOstar}
\partial_{*}\mathscr{O}:=\left\{\begin{array}{ll} \Gamma_{1},\quad &\text{if the scenario (A) occurs},\\
\partial\mathscr{O},\quad &\text{if the scenario (B) occurs},\end{array}\right.
\end{align}
and treat the previous mentioned boundary value problems together as
\begin{align}\label{eq:CombinedEllipBound}
\left\{\begin{array}{ll}\mathscr{A}u=f &\text{on }\,\mathscr{O},\\ u=g\quad
&\text{on }\,\partial_{*}\mathscr{O}.\end{array}\right.
\end{align}

\noindent
\textbf{Proof of Theorem \ref{thm:UniqueEllipBound}:} It suffices to show that if $u\in C_{loc}(\mathscr{O}\cup\partial_{*}\mathscr{O})\cap C^{2}(\mathscr{O})$ is a solution to (\ref{eq:CombinedEllipBound}) and satisfies the linear growth condition (\ref{eq:lineargrowth}), then it agrees with the Feynman-Kac formula (\ref{eq:EllipBoundFeynmanKacTau}) on every $x\in\mathscr{O}\cup\partial_{*}\mathscr{O}$, for any weak solution $(\Omega,\mathscr{F},(\mathscr{F}_{s})_{s\geq 0},\mathbb{P}^{x},W,X)$ to (\ref{eq:MainDegenSDE1})-(\ref{eq:MainDegenSDE2}) with the initial condition (\ref{eq:MainDegenSDEIniCond}) at $t=0$. From the expression (\ref{eq:EllipBoundFeynmanKacTau}) of $u_{*}$ and the boundary condition of $u$ in (\ref{eq:CombinedEllipBound}), we see that
\begin{align}\label{eq:uustarBound}
u(x)=u_{*}^{(X)}(x)=g(x),\quad\text{for any }\,x\in\partial_{*}\mathscr{O}.
\end{align}
Hence, we just need to show that $u=u_{*}^{(X)}$ on $\mathscr{O}$ for any weak solution. In the following proof, for any $x\in\mathscr{O}$, we will fix an arbitrary weak solution $(\Omega,\mathscr{F},(\mathscr{F}_{s})_{s\geq 0},\mathbb{P}^{x},W,X)$.

Let $\left\{\mathscr{O}_{k}\right\}_{k\in\mathbb{N}}$ be an increasing sequence of $C^{2,\alpha}$ open subdomains of $\mathscr{O}$ (cf.~\cite[Definition 6.2]{GilbargTrudinger:1983}) with $\alpha\in(0,1)$, such that $\bar{\mathscr{O}}_{k}\subseteq\mathscr{O}$ for each $k\in\mathbb{N}$, and $\cup_{k\in\mathbb{N}}\mathscr{O}_{k}=\mathscr{O}$. For any $x\in\mathscr{O}$, we have $x\in\mathscr{O}_{k}$ when $k$ is large enough. By It\^{o}'s formula, for any $s\geq 0$,
\begin{align}
&\quad\,\,\mathbb{E}^{x}\left[\exp\left(-\int_{0}^{s\wedge\tau_{\mathscr{O}_{k}}}c(X_{v})\,dv\right)u\left(X_{s\wedge\tau_{\mathscr{O}_{k}}}\right)\right]\nonumber\\ &=u(x)-\mathbb{E}^{x}\left[\int_{0}^{s\wedge\tau_{\mathscr{O}_{k}}}\exp\left(-\int_{0}^{v}c(X_{w})\,dw\right)\mathscr{A}u(X_{v})\,dv\right]\nonumber\\ &\quad\,+\mathbb{E}^{x}\left[\sum_{j=1}^{m}\int_{0}^{s\wedge\tau_{\mathscr{O}_{k}}}\exp\left(-\int_{0}^{v}c(X_{w})\,dw\right)\sum_{i=1}^{d}\frac{\partial u}{\partial x_{i}}(X_{v})\,\sigma_{ij}(X_{v})\,dW_{v}^{(j)}\right]\nonumber\\
\label{eq:DiscuXOkIntForm} &=u(x)-\mathbb{E}^{x}\left[\int_{0}^{s\wedge\tau_{\mathscr{O}_{k}}}\exp\left(-\int_{0}^{v}c(X_{w})\,dw\right)\mathscr{A}f(X_{v})\,dv\right],
\end{align}
since the stochastic integrals in the second sum above are martingales as the subdomain $\mathscr{O}_{k}\subset\mathscr{O}$ is bounded and $u\in C^{2}(\mathscr{O})$. We first take the limit in (\ref{eq:DiscuXOkIntForm}) as $k\rightarrow\infty$. By the growth estimate (\ref{eq:EstFirstTermJ12}), we can apply the dominated convergence theorem to obtain that
\begin{align}\label{eq:LimitktoInfuXOkIntTerm}
\lim_{k\rightarrow\infty}\mathbb{E}^{x}\!\left[\int_{0}^{s\wedge\tau_{\mathscr{O}_{k}}}\!\exp\!\left(\!-\!\int_{0}^{v}\!c(X_{w})\,dw\!\right)\!f(X_{v})\,dv\right]\!=\mathbb{E}^{x}\!\left[\int_{0}^{s\wedge\tau_{\mathscr{O}}}\!\exp\!\left(\!-\!\int_{0}^{v}\!c(X_{w})\,dw\!\right)\!f(X_{v})\,dv\right].
\end{align}
For the non-integral term on the left-hand side of (\ref{eq:DiscuXOkIntForm}), by the continuity of $u$ on $\mathscr{O}\cup\partial_{*}\mathscr{O}$ and the continuity of the sample paths of $(X_{s})_{s\geq 0}$,
\begin{align*}
\lim_{k\rightarrow\infty}\exp\left(-\int_{0}^{s\wedge\tau_{\mathscr{O}_{k}}}c(X_{v})\,dv\right)u\left(X_{s\wedge\tau_{\mathscr{O}_{k}}}\right)=\exp\left(-\int_{0}^{s}c(X_{v})\,dv\right)u\left(X_{s\wedge\tau_{\mathscr{O}}}\right),\quad\text{a.}\,\text{s.}\,.
\end{align*}
Hence, in order to show that
\begin{align}\label{eq:LimitktoInfuXOKNonIntTerm}
\lim_{k\rightarrow\infty}\mathbb{E}^{x}\left[\exp\left(-\int_{0}^{s\wedge\tau_{\mathscr{O}_{k}}}c(X_{v})\,dv\right)u\left(X_{s\wedge\tau_{\mathscr{O}_{k}}}\right)\right]=\mathbb{E}^{x}\left[\exp\left(-\int_{0}^{s}c(X_{v})\,dv\right)u\left(X_{s\wedge\tau_{\mathscr{O}}}\right)\right],
\end{align}
we only need to show that
\begin{align*}
\left\{\exp\left(-\int_{0}^{s\wedge\tau_{\mathscr{O}_{k}}}c(X_{v})\,dv\right)u\left(X_{s\wedge\tau_{\mathscr{O}_{k}}}\right):\,k\in\mathbb{N}\right\}
\end{align*}
is a collection of uniformly integrable random variables, for which it suffices to obtain the uniform boundedness of their second moments. Since $b$, $\sigma$ and $c$ satisfy (\ref{eq:SupMartCon}), by the linear growth condition (\ref{eq:lineargrowth}), Lemma \ref{lem:SupMart}, Corollary \ref{cor:SupMartLastEle} and the Optional Sampling Theorem (cf.~\cite[Theorem 1.3.22]{KaratzasShreve:1991}), we have
\begin{align*}
\mathbb{E}^{x}\!\left[\exp\left(\!-2\!\int_{0}^{s\wedge\tau_{\mathscr{O}_{k}}}\!\!c(X_{v})\,dv\!\right)\!u^{2}\!\left(X_{s\wedge\tau_{\mathscr{O}_{k}}}\right)\right]&\leq K^{2}\mathbb{E}^{x}\!\left[\exp\left(\!-2\!\int_{0}^{s\wedge\tau_{\mathscr{O}_{k}}}\!\!c(X_{v})\,dv\!\right)\!\left(1+\left\|X_{s\wedge\tau_{\mathscr{O}_{k}}}\right\|\right)^{2}\right]\\
&\leq 2K^{2}\!\left\{1+\mathbb{E}^{x}\!\left[\exp\left(\!-\!\int_{0}^{s\wedge\tau_{\mathscr{O}_{k}}}\!c(X_{v})\,dv\!\right)\!\left\|X_{s\wedge\tau_{\mathscr{O}_{k}}}\right\|^{2}\right]\right\}\\
&\leq 2K^{2}\left\{1+\|x\|^{2}+Mc_{0}^{-1}\left[1-\mathbb{E}^{x}\left(e^{-c_{0}(s\wedge\tau_{\mathscr{O}_{k}})}\right)\right]\right\}\\
&\leq 2K^{2}\left(1+\|x\|^{2}+Mc_{0}^{-1}\right).
\end{align*}
Combining (\ref{eq:DiscuXOkIntForm})-(\ref{eq:LimitktoInfuXOKNonIntTerm}), we obtain that
\begin{align}\label{eq:DiscuXOIntForm}
\mathbb{E}^{x}\!\left[\exp\!\left(-\!\int_{0}^{s\wedge\tau_{\mathscr{O}}}\!c(X_{v})\,dv\right)\!u\!\left(X_{s\wedge\tau_{\mathscr{O}}}\right)\right]=u(x)-\mathbb{E}^{x}\!\left[\int_{0}^{s\wedge\tau_{\mathscr{O}}}\!\exp\!\left(-\!\int_{0}^{v}c(X_{w})\,dw\right)\!f(X_{v})\,dv\right].
\end{align}
By the estimate (\ref{eq:EstFirstTermJ12}) and the dominated convergence theorem,
\begin{align}\label{eq:LimitttoInfuXOIntTerm}
\lim_{s\rightarrow\infty}\mathbb{E}^{x}\!\left[\int_{0}^{s\wedge\tau_{\mathscr{O}}}\!\exp\!\left(\!-\!\int_{0}^{v}c(X_{w})\,dw\!\right)\!f(X_{v})\,dv\right]=\mathbb{E}^{x}\!\left[\int_{0}^{\tau_{\mathscr{O}}}\exp\!\left(\!-\!\int_{0}^{v}c(X_{w})\,dw\!\right)\!f(X_{v})\,dv\right].
\end{align}
It remains to consider the left-hand side of (\ref{eq:DiscuXOIntForm}). Since $u\in C(\mathscr{O}\cup\partial_{*}\mathscr{O})$ solves (\ref{eq:CombinedEllipBound}), we can rewrite this term as
\begin{align*}
\mathbb{E}^{x}\left[\exp\left(-\int_{0}^{s\wedge\tau_{\mathscr{O}}}c(X_{v})\,dv\right)u\left(X_{s\wedge\tau_{\mathscr{O}}}\right)\right] &=\mathbb{E}^{x}\left[\exp\left(-\int_{0}^{\tau_{\mathscr{O}}}c(X_{v})\,dv\right)g\left(X_{\tau_{\mathscr{O}}}\right){\bf 1}_{\{\tau_{\mathscr{O}}\leq s\}}\right]\\
&\quad\,+\mathbb{E}^{x}\left[\exp\left(-\int_{0}^{s}c(X_{v})\,dv\right)u(X_{s}){\bf 1}_{\{\tau_{\mathscr{O}}>s\}}\right].
\end{align*}
Using the linear growth condition (\ref{eq:lineargrowth}) on $g$ and $u$, and a similar argument as above, we see that both collections of random variables on the right-hand side of the preceding identity
\begin{align*}
\left\{\exp\left(-\int_{0}^{\tau_{\mathscr{O}}}c(X_{v})\,dv\right)g\left(X_{\tau_{\mathscr{O}}}\right){\bf 1}_{\{\tau_{\mathscr{O}}\leq s\}}:\,\,s\geq 0\right\}
\end{align*}
and
\begin{align*}
\left\{\exp\left(-\int_{0}^{s}c(X_{v})\,dv\right)u(X_{s}){\bf 1}_{\{\tau_{\mathscr{O}}>s\}}:\,\,s\geq 0\right\}
\end{align*}
are uniformly integrable. It is easy to see that
\begin{align*}
\lim_{s\rightarrow\infty}\exp\left(-\int_{0}^{\tau_{\mathscr{O}}}c(X_{v})\,dv\right)g\left(X_{\tau_{\mathscr{O}}}\right){\bf 1}_{\{\tau_{\mathscr{O}}\leq s\}}&=\exp\left(-\int_{0}^{\tau_{\mathscr{O}}}c(X_{v})\,dv\right)g\left(X_{\tau_{\mathscr{O}}}\right){\bf 1}_{\{\tau_{\mathscr{O}}<\infty\}},\quad\text{a.}\,\text{s.}\,,\\
\lim_{s\rightarrow\infty}\exp\left(-\int_{0}^{s}c(X_{v})\,dv\right)u(X_{s}){\bf 1}_{\{\tau_{\mathscr{O}}>s\}}&=0,\quad\text{a.}\,\text{s.}\,.
\end{align*}
Hence, we have
\begin{align}\label{eq:LimitttoInfuXONonIntTerm}
\lim_{s\rightarrow\infty}\mathbb{E}^{x}\!\left[\exp\left(\!-\!\int_{0}^{s\wedge\tau_{\mathscr{O}}}\!c(X_{v})\,dv\!\right)u\!\left(X_{s\wedge\tau_{\mathscr{O}}}\right)\right]\!=\mathbb{E}^{x}\!\left[\exp\left(\!-\!\int_{0}^{\tau_{\mathscr{O}}}\!c(X_{v})\,dv\!\right)g\!\left(X_{\tau_{\mathscr{O}}}\right)\!{\bf 1}_{\{\tau_{\mathscr{O}}<\infty\}}\right].
\end{align}
Combining (\ref{eq:DiscuXOIntForm})-(\ref{eq:LimitttoInfuXONonIntTerm}), we obtain that
\begin{align*}
\mathbb{E}^{x}\left[\exp\left(-\int_{0}^{\tau_{\mathscr{O}}}c(X_{v})\,dv\right)g\left(X_{\tau_{\mathscr{O}}}\right){\bf 1}_{\{\tau_{\mathscr{O}}<\infty\}}\right]=u(x)-\mathbb{E}^{x}\left[\int_{0}^{\tau_{\mathscr{O}}}\exp\left(-\int_{0}^{v}c(X_{w})\,dw\right)f(X_{v})\,dv\right].
\end{align*}
Together with (\ref{eq:uustarBound}), we obtain that $u(x)=u_{*}^{(X)}(x)$ for any $x\in\mathscr{O}\cup\partial_{*}\mathscr{O}$ and any weak solution $(\Omega,\mathscr{F},(\mathscr{F}_{s})_{s\geq 0},\mathbb{P}^{x},W,X)$, where $u_{*}^{(X)}$ is given by (\ref{eq:EllipBoundFeynmanKacTau}).\hfill $\Box$

\medskip
\noindent
\textbf{Proof of Theorem \ref{thm:UniqueEllipBoundPartial}:} We will show that under the scenario (B) (which contains the cases (d) and (e) given as in Theorem \ref{thm:UniqueEllipBound}), and if $u\in C_{\text{loc}}(\mathscr{O}\cup\Gamma_{1})\cap C^{2}(\mathscr{O})\cap C^{1,1,\beta}_{\text{s,loc}}(\mathscr{O}\cup\Gamma_{0})$ is a solution to (\ref{eq:EllipBound}) with the partial boundary condition (\ref{eq:BoundCondGamma1}) along $\Gamma_{1}$, which satisfies the linear growth condition (\ref{eq:lineargrowth}), then it admits the stochastic representation (\ref{eq:EllipBoundFeynmanKacLambda}) for any $x\in\mathscr{O}\cup\Gamma_{1}$, and any weak solution $(\Omega,\mathscr{F},(\mathscr{F}_{s})_{s\geq 0},\mathbb{P}^{x},W,X)$ to (\ref{eq:MainDegenSDE1})-(\ref{eq:MainDegenSDE2}) with the initial condition (\ref{eq:MainDegenSDEIniCond}) at $t=0$. From (\ref{eq:EllipBoundFeynmanKacLambda}) and the boundary condition (\ref{eq:BoundCondGamma1}), we see that
\begin{align}\label{eq:uu2starBound}
u(x)=u_{**}^{(X)}(x)=g(x),\quad\text{for any }\,x\in\Gamma_{1}.
\end{align}
Hence, we just need to show that $u=u_{**}^{(X)}$ on $\mathscr{O}$ for any weak solution. In the following proof, for any $x\in\mathscr{O}$, we will fix an arbitrary weak solution $(\Omega,\mathscr{F},(\mathscr{F}_{s})_{s\geq 0},\mathbb{P}^{x},W,X)$.

Consider the following sequence of increasing subdomains of $\mathscr{O}$,
\begin{align*}
\mathscr{U}_{k}:=\left\{x\in\mathscr{O}:\,\|x\|<k,\,\text{dist}(x,\Gamma_{1})>\frac{1}{k}\right\},\quad k\in\mathbb{N},
\end{align*}
each with non-empty boundary portion $\bar{\Gamma}_{0}\cap\mathscr{U}_{k}$. For any $x\in\mathscr{O}$, we have $x\in\mathscr{U}_{k}$ when $k$ is large enough. For the simplicity of notations, we denote the process $(X^{(d)}_{s})_{s\geq 0}$ by $(Y_{s})_{s\geq 0}$ in the following proof. For any $\varepsilon>0$, define
\begin{align}\label{eq:DefSDEEpsilon}
Y^{(\varepsilon)}:=Y+\varepsilon,\quad X^{(\varepsilon)}:=\left(X^{(1)},X^{(2)},\ldots,Y^{(\varepsilon)}\right)^{T}.
\end{align}
By It\^{o}'s formula, and noticing that the stochastic integral terms are martingales up to $\lambda_{\mathscr{U}_{k}}$ since the subdomain $\mathscr{U}_{k}\subset\mathscr{O}$ is bounded and $u\in C^{2}(\mathscr{O})$, for each $k\in\mathbb{N}$, we have
\begin{align*}
\mathbb{E}^{x}\!\!\left[\exp\!\left(\!-\!\!\int_{0}^{s\wedge\lambda_{\mathscr{U}_{k}}}\!\!c(X^{(\varepsilon)}_{v})\,dv\!\right)\!u\!\left(X^{(\varepsilon)}_{s\wedge\lambda_{\mathscr{U}_{k}}}\right)\!\right]\!=u(x)-\mathbb{E}^{x}\!\!\left[\int_{0}^{s\wedge\lambda_{\mathscr{U}_{k}}}\!\!\exp\!\left(\!-\!\!\int_{0}^{v}c(X^{(\varepsilon)}_{w})dw\!\right)\!\mathscr{A}^{\varepsilon}\!u(X^{(\varepsilon)}_{v})\,dv\right],
\end{align*}
where $\mathscr{A}^{\varepsilon}$ denotes the elliptic differential operator
\begin{align}\label{eq:DefGenEps}
-\mathscr{A}^{\varepsilon}v(x):=-\mathscr{A}v(x)+\langle b(x^{(-\varepsilon)})-b(x),Dv(x)\rangle+\frac{1}{2}\,\text{tr}\left(\left(a(x^{(-\varepsilon)})-a(x)\right)D^{2}v(x)\right),
\end{align}
and where $x^{(-\varepsilon)}=(x_{1},\ldots,x_{d-1},x_{d}-\varepsilon)^{T}$. Since $u$ is a solution to (\ref{eq:EllipBound}), we can rewrite the preceding identity as
\begin{align}
&\quad\,\,\mathbb{E}^{x}\left[\exp\left(-\int_{0}^{s\wedge\lambda_{\mathscr{U}_{k}}}c(X^{(\varepsilon)}_{v})\,dv\right)u\left(X^{(\varepsilon)}_{s\wedge\lambda_{\mathscr{U}_{k}}}\right)\right]\nonumber\\ &=u(x)-\mathbb{E}^{x}\left[\int_{0}^{s\wedge\lambda_{\mathscr{U}_{k}}}\exp\left(-\int_{0}^{v}c(X^{(\varepsilon)}_{w})\,dw\right)f(X^{(\varepsilon)}_{v})\,dv\right]\nonumber\\
\label{eq:ItoXeps} &\quad\,-\mathbb{E}^{x}\left[\int_{0}^{s\wedge\lambda_{\mathscr{U}_{k}}}\exp\left(-\int_{0}^{v}c(X^{(\varepsilon)}_{w})\,dw\right)\left(\mathscr{A}^{\varepsilon}-\mathscr{A}\right)u(X^{(\varepsilon)}_{v})\,dv\right].
\end{align}
We shall first take the limit in (\ref{eq:ItoXeps}) as $\varepsilon\downarrow 0$ (for fixed $k\geq 1$ and $t\geq 0$). Using Assumption \ref{assupt:UnifPostc}, the continuity of $f$, $u$ and $c$ on compact subsets of $\mathscr{O}\cup\Gamma_{0}$, as well as the dominated convergence theorem, we have
\begin{align}\label{eq:LimitEpsLeft}
\lim_{\varepsilon\downarrow 0}\mathbb{E}^{x}\!\left[\exp\!\left(\!-\!\!\int_{0}^{s\wedge\lambda_{\mathscr{U}_{k}}}\!\!c(X^{(\varepsilon)}_{v})\,dv\!\right)\!u\!\left(X^{(\varepsilon)}_{s\wedge\lambda_{\mathscr{U}_{k}}}\right)\right]&=\mathbb{E}^{x}\!\left[\exp\!\left(\!-\!\!\int_{0}^{s\wedge\lambda_{\mathscr{U}_{k}}}\!\!c(X_{v})\,dv\!\right)\!u\!\left(X_{s\wedge\lambda_{\mathscr{U}_{k}}}\right)\right],\\
\label{eq:LimitEpsRight2} \lim_{\varepsilon\downarrow 0}\mathbb{E}^{x}\!\!\left[\int_{0}^{s\wedge\lambda_{\mathscr{U}_{k}}}\!\!\!\exp\!\left(\!-\!\!\int_{0}^{v}\!c(X^{(\varepsilon)}_{w})dw\!\right)\!f(X^{(\varepsilon)}_{v})dv\right]\!&=\mathbb{E}^{x}\!\!\left[\int_{0}^{s\wedge\lambda_{\mathscr{U}_{k}}}\!\!\!\exp\!\left(\!-\!\!\int_{0}^{v}\!c(X_{w})dw\!\right)\!f(X_{v})dv\right].
\end{align}
To estimate the last integral term on the right-hand side of (\ref{eq:ItoXeps}), we assume without loss of generality that $\varepsilon<1/2k$ for any fixed $k\geq 1$. For any $v\in[0,s\wedge\lambda_{\mathscr{U}_{k}}]$,
\begin{align}
\left|\left(\mathscr{A}^{\varepsilon}\!\!-\!\!\mathscr{A}\right)u(X^{(\varepsilon)}_{v})\right|&\leq\left\|b(X_{v})\!-\!b(X_{v}^{(\varepsilon)})\right\|\!\left\|Du(X^{(\varepsilon)}_{v})\right\|\!+\!\frac{1}{2}\left|\text{tr}\!\left[\left(\!Y_{v}^{\beta}\tilde{a}(X_{v})\!-\!\!\left(Y_{v}^{(\varepsilon)}\right)^{\beta}\!\!\tilde{a}(X_{v}^{(\varepsilon)})\!\right)\!D^{2}u(X^{(\varepsilon)}_{v})\right]\right|\nonumber\\
&\leq\left\|b(X_{v})-b(X_{v}^{(\varepsilon)})\right\|\left\|Du\right\|_{C(\bar{\mathscr{U}}_{2k})}\!+\!\frac{1}{2}\left|\text{tr}\left[\left(Y_{v}^{(\varepsilon)}\right)^{\beta}\left(\tilde{a}(X_{v})\!-\!\tilde{a}(X_{v}^{(\varepsilon)})\right)D^{2}u(X^{(\varepsilon)}_{v})\right]\right|\nonumber\\
&\quad\,+\frac{1}{2}\left|\text{tr}\left[\left(Y_{v}^{\beta}-\left(Y_{v}^{(\varepsilon)}\right)^{\beta}\right)\tilde{a}(X_{v})D^{2}u(X^{(\varepsilon)}_{v})\right]\right|\nonumber\\
&\leq\left\|b(X_{v})-b(X_{v}^{(\varepsilon)})\right\|\left\|Du\right\|_{C(\bar{\mathscr{U}}_{2k})}+\frac{1}{2}\left\|\tilde{a}(X_{v})-\tilde{a}(X_{v}^{(\varepsilon)})\right\|\left\|x_{d}^{\beta}D^{2}u\right\|_{C(\bar{\mathscr{U}}_{2k})}\nonumber\\
\label{eq:EstResidAEpsA} &\quad+\frac{1}{2}\left\|\left(Y_{v}^{\beta}-\left(Y_{v}^{(\varepsilon)}\right)^{\beta}\right)D^{2}u(X^{(\varepsilon)}_{v})\right\|\left\|\tilde{a}\right\|_{C(\bar{\mathscr{U}}_{2k})},
\end{align}
where $\|\cdot\|_{C(\bar{\mathscr{U}}_{2k})}$ denotes the uniform norm on $\bar{\mathscr{U}}_{2k}$. Since $b$ and $\tilde{a}$ are uniformly continuous and bounded on $\bar{\mathscr{U}}_{2k}$, the first two terms on the right-hand side of (\ref{eq:EstResidAEpsA}) converge to $0$ as $\varepsilon\downarrow 0$. For the last term in (\ref{eq:EstResidAEpsA}), notice that
\begin{align}
&\quad\,\,\left\|\left(Y_{v}^{\beta}-\left(Y_{v}^{(\varepsilon)}\right)^{\beta}\right)D^{2}u(X^{(\varepsilon)}_{v})\right\|\nonumber\\
&\leq\varepsilon^{-\beta/2}\left|\left(Y_{v}^{(\varepsilon)}\right)^{\beta}-Y_{v}^{\beta}\right|\left\|x_{d}^{\beta}D^{2}u\right\|_{C(\bar{\mathscr{U}}_{2k})}{\bf 1}_{\{Y_{v}^{(\varepsilon)}\geq\sqrt{\varepsilon}\}}+2^{\beta}\left\|x_{d}^{\beta}D^{2}u\right\|_{C(\bar{\mathscr{U}}_{2k})}{\bf 1}_{\{Y_{v}^{(\varepsilon)}<\sqrt{\varepsilon}\}}\nonumber\\
&\leq\beta\max\left\{(k+\varepsilon)^{\beta-1},\varepsilon^{(\beta-1)/2}\right\}\varepsilon^{1-\beta/2}\left\|x_{d}^{\beta}D^{2}u\right\|_{C(\bar{\mathscr{U}}_{2k})}+2^{\beta}\left\|x_{d}^{\beta}D^{2}u\right\|_{C(\bar{\mathscr{U}}_{2k})}{\bf 1}_{\{Y_{v}^{(\varepsilon)}<\sqrt{\varepsilon}\}}\nonumber\\
\label{eq:EstEpsXdXdEpsHessian} &=\beta\max\left\{(k+\varepsilon)^{\beta-1}\varepsilon^{1-\beta/2},\sqrt{\varepsilon}\right\}\left\|x_{d}^{\beta}D^{2}u\right\|_{C(\bar{\mathscr{U}}_{2k})}+2^{\beta}\left\|x_{d}^{\beta}D^{2}u\right\|_{C(\bar{\mathscr{U}}_{2k})}{\bf 1}_{\{Y_{v}^{(\varepsilon)}<\sqrt{\varepsilon}\}},
\end{align}
where in the second inequality we use the convexity or concavity, depending on $\beta\in(0,2)$, of the function $p(x)=x^{\beta}$ on $[\sqrt{\varepsilon},k+\varepsilon]$:
\begin{align*}
(x+\varepsilon)^{\beta}-x^{\beta}&\leq p'(\sqrt{\varepsilon})\varepsilon=\beta\varepsilon^{(\beta-1)/2}\cdot\varepsilon,\quad\beta\in(0,1],\\
(x+\varepsilon)^{\beta}-x^{\beta}&\leq p'(k+\varepsilon)\varepsilon=\beta(k+\varepsilon)^{\beta-1}\cdot\varepsilon,\quad\beta\in[1,2].
\end{align*}
Combining (\ref{eq:EstResidAEpsA})-(\ref{eq:EstEpsXdXdEpsHessian}), using the assumption $u\in C_{\text{s,loc}}^{1,1,\beta}(\mathscr{O}\cup\Gamma_{0})$ and the fact that
\begin{align*}
\lim_{\varepsilon\downarrow 0}{\bf 1}_{\{Y_{v}^{(\varepsilon)}<\sqrt{\varepsilon}\}}=0,\quad\text{a.}\,\text{s.}\,,
\end{align*}
we obtain, by the dominated convergence theorem, that
\begin{align}\label{eq:LimitEpsRight3}
\lim_{\varepsilon\downarrow 0}\mathbb{E}^{x}\left[\int_{0}^{s\wedge\lambda_{\mathscr{U}_{k}}}\exp\left\{-\int_{0}^{v}c(X^{(\varepsilon)}_{w})\,dw\right\}\left(\mathscr{A}^{\varepsilon}-\mathscr{A}\right)u(X^{(\varepsilon)}_{v})\,dv\right]=0.
\end{align}
Therefore, by combining (\ref{eq:ItoXeps}), (\ref{eq:LimitEpsLeft}), (\ref{eq:LimitEpsRight2}) and (\ref{eq:LimitEpsRight3}), we have
\begin{align}\label{eq:ItoForXk}
\mathbb{E}^{x}\!\!\left[\exp\!\left(\!-\!\!\int_{0}^{s\wedge\lambda_{\mathscr{U}_{k}}}\!\!c(X_{v})dv\!\right)\!u\!\left(X_{s\wedge\lambda_{\mathscr{U}_{k}}}\right)\right]\!=u(x)\!-\!\mathbb{E}^{x}\!\!\left[\int_{0}^{s\wedge\lambda_{\mathscr{U}_{k}}}\!\!\!\exp\!\left(\!-\!\!\int_{0}^{v}\!c(X_{w})dw\!\right)\!f(X_{v})\,dv\right].
\end{align}
As $k\rightarrow\infty$ and { $s\rightarrow\infty$}, clearly we have
\begin{align}\label{eq:LimitLambdaktinfty}
s\wedge\lambda_{\mathscr{U}_{k}}\rightarrow\lambda_{\mathscr{O}},\quad\text{a.}\,\text{s.}\,.
\end{align}
By the same argument as in the proof of Theorem \ref{thm:UniqueEllipBound}, we can take the limit in (\ref{eq:ItoForXk}), as $k\rightarrow\infty$ and $s\rightarrow\infty$, to obtain that, for any $x\in\mathscr{O}$,
\begin{align*}
\mathbb{E}^{x}\!\left[\exp\left(-\!\int_{0}^{\lambda_{\mathscr{O}}}\!c(X^{(\varepsilon)}_{v})\,dv\right)g\!\left(X_{\lambda_{\mathscr{O}}}\right){\bf 1}_{\{\lambda_{\mathscr{O}}<\infty\}}\right]=u(x)-\mathbb{E}^{x}\!\left[\int_{0}^{\lambda_{\mathscr{O}}}\exp\left(-\!\int_{0}^{v}\!c(X_{w})\,dw\right)f(X_{v})\,dv\right].
\end{align*}
Together with (\ref{eq:uu2starBound}), we have shown that $u=u_{**}^{(X)}$ on $\mathscr{O}\cup\Gamma_{1}$, where $u_{**}^{(X)}$ is given by (\ref{eq:EllipBoundFeynmanKacLambda}).\hfill $\Box$

\medskip
Next, we prove existence of solutions to (\ref{eq:CombinedEllipBound}) when the boundary condition $g$ is H\"{o}lder continuous on suitable portion of $\mathscr{O}$.

\medskip
\noindent
\textbf{Proof of Theorem \ref{thm:ExistenceHolderCont}:} For each $x\in\mathscr{O}\cup\partial_{*}\mathscr{O}$, let $(\Omega,\mathscr{F},(\mathscr{F}_{s})_{s\geq 0},\mathbb{P}^{x},W,X)$ be any weak solution to (\ref{eq:MainDegenSDE1})-(\ref{eq:MainDegenSDE2}) with the initial condition (\ref{eq:MainDegenSDEIniCond}) at $t=0$, and let $u_{*}^{(X)}(x)$ be defined as in (\ref{eq:EllipBoundFeynmanKacTau}). It suffices to show that, $u_{*}^{(X)}$ is a solution to (\ref{eq:CombinedEllipBound}), such that $u_{*}^{(X)}\in C_{\text{loc}}(\mathscr{O}\cup\partial_{*}\mathscr{O})\cap C_{\text{loc}}^{2,\alpha}(\mathscr{O}\cup\Gamma_{1})$ and $u_{*}^{(X)}$ satisfies the linear growth condition (\ref{eq:lineargrowth}).

The fact that $u_{*}^{(X)}$ satisfies the linear growth condition (\ref{eq:lineargrowth}) follows from Lemma \ref{lem:EstFunJ} with $\theta_{1}=\tau_{\mathscr{O}}^{x,X}$, $\theta_{2}=0$ and $\psi=0$, for every $x\in\mathscr{O}\cup\partial_{*}\mathscr{O}$. It remains to show that $u_{*}^{(X)}\in C_{\text{loc}}(\mathscr{O}\cup\partial_{*}\mathscr{O})\cap C_{\text{loc}}^{2,\alpha}(\mathscr{O}\cup\Gamma_{1})$ and $u_{*}^{(X)}$ is a solution to (\ref{eq:CombinedEllipBound}). In the following proof, without loss of generality, we will fix a single filtered probability space $(\Omega,\mathscr{F},(\mathscr{F}_{s})_{s\geq 0},\mathbb{P})$, on which we fix a weak solution $(X(x),W)$ to (\ref{eq:MainDegenSDE1})-(\ref{eq:MainDegenSDEIniCond}) for each initial data $x\in\mathscr{O}\cup\partial_{*}\mathscr{O}$ at $t=0$. Indeed, if $(\Omega,\mathscr{F},(\mathscr{F}_{s})_{s\geq 0},\mathbb{P},X,W)$ is a weak solution to (\ref{eq:MainDegenSDE1})-(\ref{eq:MainDegenSDEIniCond}) with the initial data $x\in\mathscr{O}\cup\partial_{*}\mathscr{O}$, then $(\Omega,\mathscr{F},(\mathscr{F}_{s})_{s\geq 0},\mathbb{P},X+y-x,W)$ is a weak solution on the same probability space with a different initial data $y\in\mathscr{O}\cup\partial_{*}\mathscr{O}$. We will still use $\mathbb{P}^{x}$ and $\mathbb{E}^{x}$ to denote the probability and the expectation with respect to different initial data, in which the parameter $x$ of the weak solution $X$ as well as all stopping times will be omitted. We will also omit the superscript $X$ of $u_{*}$ and all stopping times for simplicity. The proof is organized in the following two steps.

\medskip
\noindent
\textbf{Step 1} ($u_{*}\in C^{2,\alpha}_{\text{loc}}(\mathscr{O}\cup\Gamma_1)$, and $u_{*}$ is a solution to (\ref{eq:CombinedEllipBound}))$\,$ We first notice that $u_{*}=g$ on $\Gamma_{1}$ follows directly from (\ref{eq:EllipBoundFeynmanKacTau}). Let $(\mathscr{D}_{k})_{k\in\mathbb{N}}$ be an increasing sequence of $C^{2,\alpha}$ subdomains of $\mathscr{O}$ (cf.~\cite[Definition 6.2]{GilbargTrudinger:1983}) such that $\bigcup_{k\in\mathbb{N}}\mathscr{D}_{k}=\mathscr{O}$ and
\begin{align*}
\mathscr{O}\cap(-k,k)^{d-1}\times(1/k,k)\subset\mathscr{D}_k\subset\mathscr{O}\cap(-2k,2k)^{d-1}\times(1/2k,2k),\quad k\in\mathbb{N}.
\end{align*}
Notice that on each domain $\mathscr{D}_k$, the differential operator $-\mathscr{A}$ is uniformly elliptic with bounded and $C^{0,\alpha}(\bar{\mathscr{D}_k})$ coefficients. Moreover, from our hypotheses, we have $f\in C^{0,\alpha}(\bar{\mathscr{D}}_{k})$ and $g\in C^{2,\alpha}(\bar{\mathscr{D}}_{k})$. Therefore, \cite[Theorem 6.13]{GilbargTrudinger:1983} implies that the elliptic boundary value problem
\begin{align}\label{eq:CombinedEllipBoundDk}
\left\{\begin{array}{ll}\mathscr{A}u=f\,\,\,\,\,\,\text{ on }\,\mathscr{D}_{k},\\ u=g\quad\,\,\,\,\,\,\,\text{ on }\,\partial\mathscr{D}_{k}\end{array}\right.
\end{align}
admits a unique solution $u_{k}\in C(\bar{\mathscr{D}}_{k})\cap C^{2,\alpha}(\mathscr{D}_{k})$. Moreover, by the classical theory on stochastic representations of solutions to uniformly elliptic partial differential equations (cf.~\cite[Theorem 6.5.1]{Friedman:1976}, \cite[Proposition 5.7.2]{KaratzasShreve:1991} and~\cite[Theorem 9.1.1, Corollary 9.1.2]{Oksendal:2003}), $u_{k}$ admits the unique stochastic representation on $\bar{\mathscr{D}}_{k}$: for any $x\in\bar{\mathscr{D}}_{k}$,
\begin{align*}
u_{k}(x)=\mathbb{E}^{x}\!\left[\exp\!\left(-\!\int_{0}^{\tau_{\mathscr{D}_{k}}}c(X_{s})\,ds\right)\!g(X_{\tau_{\mathscr{D}_{k}}}){\bf 1}_{\{\tau_{\mathscr{D}_{k}}<\infty\}}\right]+\mathbb{E}^{x}\!\left[\int_{0}^{\tau_{\mathscr{O}_{k}}}\exp\!\left(-\!\int_{0}^{s}c(X_{v})\,dv\right)\!f(X_{s})\,ds\right].
\end{align*}
Here, we notice that, for every $x\in\bar{\mathscr{D}}_{k}$, any weak solution within the time period $[0,\tau_{\mathscr{D}_{k}}^{x}]$ has the unique law (Remark \ref{rem:UniqueinLaw}). Also, by Lemma \ref{lem:EstFunJ}, $u_{k}$ obeys the linear growth condition (\ref{eq:lineargrowth}), for all $k\in\mathbb{N}$. Since, for each $x\in\mathscr{O}$ ($x\in\mathscr{D}_{k}$ for $k$ large enough), $(\tau_{\mathscr{D}_{k}}^{x})_{k\in\mathbb{N}}$ is an increasing sequence of stopping times, which converges to $\tau_{\mathscr{O}}^{x}$ $\mathbb{P}$-a.$\,$s., as $k\rightarrow\infty$, using the hypothesis $g\in C_{\text{loc}}(\mathscr{O}\cup\partial_{*}\mathscr{O})$, the continuity of sample paths of $X$, the linear growth condition (\ref{eq:lineargrowth}) on $(u_{k})_{k\in\mathbb{N}}$ and Lemma \ref{lem:EstFunJ}, we can apply the same argument used in the proof of Theorem \ref{thm:UniqueEllipBound} to obtain that
\begin{align*}
\lim_{k\rightarrow\infty}u_{k}(x)=u_{*}(x),\quad\text{ for any }\,x\in\mathscr{O}.
\end{align*}

Fix $x^{(0)}\in\mathscr{O}$, and let $B:=B(x^{(0)},r_{0})$ be a Euclidean ball centered at $x^{(0)}$ with radius $r_{0}$, such that $\bar{B}\subset\mathscr{O}$. Set $B_{1/2}:=B(x^{(0)},r_{0}/2)$. Since $(u_{k})_{k\in\mathbb{N}}$ obeys the linear growth condition (\ref{eq:lineargrowth}), the same argument as in the proof of Lemma \ref{lem:EstFunJ} shows that $(u_{k})_{k\in\mathbb{N}}$ is uniformly bounded on $\bar{B}$. By the interior Schauder estimates (cf.~\cite[Corollary 6.3]{GilbargTrudinger:1983}), the sequence $(u_{k})_{k\in\mathbb{N}}$ has uniformly bounded $C^{2,\alpha}(\bar{B}_{1/2})$ norms. Compactness of the embedding $C^{2,\alpha}(\bar{B}_{1/2})\hookrightarrow C^{2,\gamma}(\bar{B}_{1/2})$, for any $0\leq\gamma<\alpha$, shows that, after passing to a subsequence, $(u_{k})_{k\in\mathbb{N}}$ converges in $C^{2,\gamma}(\bar{B}_{1/2})$ to $u_{*}\in C^{2,\gamma}(\bar{B}_{1/2})$, and hence $\mathscr{A}u_{*}=f$ on $\bar{B}_{1/2}$. Since this subsequence has uniformly bounded $C^{2,\alpha}(\bar{B}_{1/2})$ norms and it converges strongly in $C^{2}(\bar{B}_{1/2})$ to $u_{*}$, we obtain that $u_{*}\in C^{2,\alpha}(\bar{B}_{1/2})$. Since $x^{(0)}\in\mathscr{O}$ is arbitrarily chosen, we conclude that $u_{*}$ is a solution to (\ref{eq:CombinedEllipBound}) and $u_{*}\in C^{2,\alpha}(\mathscr{O})$.

Next, we fix any $z_{0}\in\Gamma_{1}$ and choose $r_{0}>0$ small enough such that $B(z_{0},r_{0})\cap\Gamma_{0}=\emptyset$. Set
\begin{align*}
D:=B(z_{0},r_{0})\cap\mathscr{O}\quad\text{and}\quad D':=B(z_{0},r_{0}/2)\cap\mathscr{O}.
\end{align*}
From the construction of $(\mathscr{D}_{k})_{k\in\mathbb{N}}$, we can find $k_{0}\in\mathbb{N}$ large enough such that $D\subset\mathscr{D}_{k}$, for all $k\geq k_{0}$. Since $f\in C^{0,\alpha}(\bar{D})$ and $g\in C^{2,\alpha}(\bar{D})$, by~\cite[Corollary 6.7]{GilbargTrudinger:1983} and the fact that $u_{k}$ solves (\ref{eq:CombinedEllipBoundDk}), we have
\begin{align}\label{eq:ukC2alphaEst}
\|u_{k}\|_{C^{2,\alpha}(\bar{D'})}\leq C\left(\|u_{k}\|_{C(\bar D)}+\|g\|_{C^{2,\alpha}(\bar{D})}+\|f\|_{C^{0,\alpha}(\bar{D})}\right),\quad k\geq k_{0},
\end{align}
where $C>0$ is a constant depending only on the coefficients of $\mathscr{A}$, and the domains $D$ and $D'$. Combining (\ref{eq:ukC2alphaEst}) with the uniform bound on the $C(\bar{D})$ norms of $(u_{k})_{k\in\mathbb{N}}$, resulting from the linear growth condition (\ref{eq:lineargrowth}) and a similar argument as in the proof of Lemma \ref{lem:EstFunJ}, the compactness of the embedding of $C^{2,\alpha}(\bar{D'})\hookrightarrow C^{2,\gamma}(\bar{D'})$, where $0\leq\gamma<\alpha$, implies that there exists a subsequence of $(u_{k})_{k\in\mathbb{N}}$ which converges strongly to $u_{*}$ in $C^{2,\gamma}(\bar{D'})$. In particular, this subsequence converges in $C^{2}(\bar{D'})$ to $u_{*}$, and therefore $u_{*}\in C^{2,\alpha}(\bar{D'})$. Combining the above two cases, we conclude that $u_{*}\in C^{2,\alpha}_{\text{loc}}(\mathscr{O}\cup\Gamma_{1})$.
\medskip
\noindent
\textbf{Step 2} ($u_{*}$ is continuous { up to} $\bar{\Gamma}_{0}$ in the scenario (B))$\,$ Let $x^{(0)}=(x_{1}^{(0)},\ldots,x_{d-1}^{(0)},0)\in\bar{\Gamma}_{0}$. For any $x=(x_{1},\ldots,x_{d})\in\mathscr{O}$, let
\begin{align*}
\theta^{x}&:=\inf\left\{s\geq 0:\,X_{s}\in\partial{\mathbb{H}},\,X_{0}=x\right\},\\
T_{0}^{x_{d}}&:=\inf\left\{s\geq 0:\,X^{(d)}_{s}=0,\,X^{(d)}_{0}=x_{d}\right\}.
\end{align*}
Clearly we have
\begin{align}\label{eq:IneqHitTimes}
\tau_{\mathscr{O}}^{x}\leq\theta^{x}\leq T_{0}^{x_{d}},\quad\mathbb{P}-\text{a.}\,\text{s.}\,.
\end{align}
It follows from (\ref{eq:HitTimeLeftBound}) and (\ref{eq:IneqHitTimes}) that
\begin{align}\label{eq:LimitTauTheta}
\lim_{x\rightarrow x^{(0)}}\theta^{x}=\lim_{x\rightarrow x^{(0)}}\tau_{\mathscr{O}}^{x}=0\quad { \text{uniformly in }\,\,(x_{1},\ldots,x_{d-1})\in\mathbb{R}^{d-1},\,\,\,\,\mathbb{P}-\text{a.}\,\text{s.}}.
\end{align}
Therefore, by the dominated convergence theorem,
\begin{align}\label{eq:EstIntTermux}
\lim_{x\rightarrow x^{(0)}}\mathbb{E}^{x}\left[\int_{0}^{\tau_{\mathscr{O}}}\exp\left(-\int_{0}^{s}c(X_{v})\,dv\right)f(X_{s})\,ds\right]=0.
\end{align}

Next, we need to show that
\begin{align}\label{eq:ExistLimitg}
\lim_{x\rightarrow x^{(0)}}\mathbb{E}^{x}\left[\exp\left(-\int_{0}^{\tau_{\mathscr{O}}}c(X_{s})\,ds\right)g(X_{\tau_{\mathscr{O}}}){\bf 1}_{\{\tau_{\mathscr{O}}<\infty\}}\right]=g(x^{(0)}).
\end{align}
Notice that
\begin{align*}
&\quad\,\,\mathbb{E}^{x}\left[\exp\left(-\int_{0}^{\tau_{\mathscr{O}}}c(X_{s})\,ds\right)g(X_{\tau_{\mathscr{O}}}){\bf 1}_{\{\tau_{\mathscr{O}}<\infty\}}\right]-g(x^{(0)})\nonumber\\
&=\mathbb{E}^{x}\!\!\left[\exp\!\left(\!-\!\!\int_{0}^{\tau_{\mathscr{O}}}\!\!c(X_{s})ds\!\right)\!\!\left(\!g(X_{\tau_{\mathscr{O}}})\!-\!g(x^{(0)})\!\right)\!\!{\bf 1}_{\{\tau_{\mathscr{O}}<\infty\}}\!\right]\!\!+\!g(x^{(0)})\!\!\left[\!1\!-\!\mathbb{E}^{x}\!\!\left(\!\exp\!\left(\!-\!\!\int_{0}^{\tau_{\mathscr{O}}}\!\!c(X_{s})ds\!\right)\!\!{\bf 1}_{\{\tau_{\mathscr{O}}<\infty\}}\!\right)\!\right],
\end{align*}
and by Assumption \ref{assupt:UnifPostc}, (\ref{eq:LimitTauTheta}) and the dominated convergence theorem,
\begin{align}\label{eq:LimitLapTauFinite}
\lim_{x\rightarrow x^{(0)}}\mathbb{E}^{x}\left[\exp\left(-\int_{0}^{\tau_{\mathscr{O}}}c(X_{s})\,ds\right){\bf 1}_{\{\tau_{\mathscr{O}}<\infty\}}\right]=1,
\end{align}
we just need to show that
\begin{align*}
\lim_{x\rightarrow x^{(0)}}\mathbb{E}^{x}\left[\exp\left(-\int_{0}^{\tau_{\mathscr{O}}}c(X_{s})\,ds\right)\left(g(X_{\tau_{\mathscr{O}}})-g(x^{(0)})\right){\bf 1}_{\{\tau_{\mathscr{O}}<\infty\}}\right]=0.
\end{align*}
Let $\varepsilon>0$ be fixed. By the continuity of $g$, we may choose $\delta_{1,\varepsilon}>0$ such that
\begin{align}\label{eq:Contg}
\left|g(x)-g(x^{(0)})\right|\leq\varepsilon,\quad\forall x\in B(x^{(0)},\delta_{1,\varepsilon})\cap\partial\mathscr{O}.
\end{align}
For this $\delta_{1,\varepsilon}>0$, by~\cite[Equation (5.3.18) in Problem 5.3.15]{KaratzasShreve:1991}, there exists $C_{1}>0$, depending on $x^{(0)}$ and $\delta_{1,\varepsilon}>0$, such that { for any $s\geq 0$,}
\begin{align*}
\sup_{x\in B(x^{(0)},\delta_{1,\varepsilon})\cap\mathscr{O}}\mathbb{E}^{x}\left(\sup_{0\leq v\leq s}\left|X_{v}-x\right|\right)\leq C_{1}\sqrt{s},
\end{align*}
which implies that, by choosing $s_{\varepsilon}>0$ small enough, for any $s\in(0,s_{\varepsilon}]$,
\begin{align}\label{eq:EstProbXx0delta1}
\sup_{x\in B(x^{(0)},\delta_{1,\varepsilon})\cap\mathscr{O}}\mathbb{P}^{x}\left(\sup_{0\leq v\leq s}\left|X_{v}-x\right|\geq\delta_{1,\varepsilon}/2\right)\leq\frac{2C_{1}\sqrt{s}}{\delta_{1,\varepsilon}}{ \leq\varepsilon}.
\end{align}
Moreover, applying (\ref{eq:HitTimeLeftBound}) to $X^{(d)}(\delta_{2,\varepsilon})$, we may choose $\delta_{2,\varepsilon}>0$ sufficiently small, such that
\begin{align}\label{eq:EstT0delta2}
\mathbb{P}\left(T_{0}^{\delta_{2,\varepsilon}}>s\right)\leq\varepsilon,\quad{ \text{for any }\,s\in(0,s_{\varepsilon}]}.
\end{align}
Let $\delta_{\varepsilon}:=\min\{\delta_{1,\varepsilon}/2,\delta_{2,\varepsilon}\}$, for any $s\in(0,s_{\varepsilon}]$ and any $x\in B(x^{(0)},\delta_{\varepsilon})\cap\mathscr{O}$,
\begin{align}
&\quad\,\,\mathbb{E}^{x}\left[\exp\left(-\int_{0}^{\tau_{\mathscr{O}}}c(X_{v})\,dv\right)\left(g(X_{\tau_{\mathscr{O}}})-g(x^{(0)})\right){\bf 1}_{\{\tau_{\mathscr{O}}<\infty\}}\right]\nonumber\\
&=\mathbb{E}^{x}\left[\exp\left(-\int_{0}^{\tau_{\mathscr{O}}}c(X_{v})\,dv\right)\left(g(X_{\tau_{\mathscr{O}}})-g(x^{(0)})\right){\bf 1}_{\{\tau_{\mathscr{O}}\leq s,\,\sup_{v\in[0,s]}|X_{v}-x|<\delta_{1,\varepsilon}/2\}}\right]\nonumber\\
&\quad\,+\mathbb{E}^{x}\left[\exp\left(-\int_{0}^{\tau_{\mathscr{O}}}c(X_{v})\,dv\right)\left(g(X_{\tau_{\mathscr{O}}})-g(x^{(0)})\right){\bf 1}_{\{\tau_{\mathscr{O}}\leq s,\,\sup_{v\in[0,s]}|X_{v}-x|\geq\delta_{1,\varepsilon}/2\}}\right]\nonumber\\
\label{eq:DecomgXtaux0} &\quad\,+\mathbb{E}^{x}\left[\exp\left(-\int_{0}^{\tau_{\mathscr{O}}}c(X_{v})\,dv\right)\left(g(X_{\tau_{\mathscr{O}}})-g(x^{(0)})\right){\bf 1}_{\{\tau_{\mathscr{O}}\in(s,\infty)\}}\right].
\end{align}
By (\ref{eq:Contg}), for any $x\in B(x^{(0)},\delta_{\varepsilon})\cap\mathscr{O}$,
\begin{align}\label{eq:EstgXtaux01}
\mathbb{E}^{x}\left[\exp\left(-\int_{0}^{\tau_{\mathscr{O}}}c(X_{v})\,dv\right)\left|g(X_{\tau_{\mathscr{O}}})-g(x^{(0)})\right|{\bf 1}_{\{\tau_{\mathscr{O}}\leq s,\,\sup_{v\in[0,s]}|X_{v}-x|<\delta_{1,\varepsilon}/2\}}\right]\leq\varepsilon.
\end{align}
To estimate the second and the third terms on the right-hand side of (\ref{eq:DecomgXtaux0}), we first notice that, by the linear growth condition (\ref{eq:lineargrowth}) on $g$, Lemma \ref{lem:SupMart}, Corollary \ref{cor:SupMartLastEle} and the Optional Sampling Theorem (cf.~\cite[Theorem 1.3.22]{KaratzasShreve:1991}), for any $x\in B(x^{(0)},\delta_{\varepsilon})\cap\mathscr{O}$,
\begin{align*}
\mathbb{E}^{x}\left[\exp\left(-\int_{0}^{\tau_{\mathscr{O}}}c(X_{v})\,dv\right)\left|g(X_{\tau_{\mathscr{O}}})\right|^{2}\right]&\leq 2K\left\{1+\mathbb{E}^{x}\left[\exp\left(-\int_{0}^{\tau_{\mathscr{O}}}c(X_{v})\,dv\right)\left\|X_{\tau_{\mathscr{O}}}\right\|^{2}\right]\right\}\\
&\leq 2K\left[1+\|x\|^{2}+Mc_{0}^{-1}\left(1-\mathbb{E}^{x}\left(e^{-c_{0}\tau_{\mathscr{O}}}\right)\right)\right]\\
&\leq C_{2}:=\sup_{x\in B(x^{(0)},\delta)\cap\mathscr{O}}2K\left(1+\|x\|^{2}+Mc_{0}^{-1}\right).
\end{align*}
By the Cauchy-Schwarz inequality and (\ref{eq:EstProbXx0delta1}), for any $s\in(0,s_{\varepsilon}]$ and any $x\in B(x^{(0)},\delta_{\varepsilon})\cap\mathscr{O}$,
\begin{align}
&\quad\,\,\mathbb{E}^{x}\left[\exp\left(-\int_{0}^{\tau_{\mathscr{O}}}c(X_{v})\,dv\right)\left|g(X_{\tau_{\mathscr{O}}})-g(x^{(0)})\right|{\bf 1}_{\{\tau_{\mathscr{O}}\leq s,\,\sup_{v\in[0,s]}|X_{v}-x|\geq\delta_{1,\varepsilon}/2\}}\right]\nonumber\\
\label{eq:EstgXtaux02} &\leq\left\{\mathbb{E}^{x}\!\left[\exp\!\left(\!-\!\!\int_{0}^{\tau_{\mathscr{O}}}\!\!c(X_{v})dv\!\right)\!\left|g(X_{\tau_{\mathscr{O}}})\!-\!g(x^{(0)})\right|^{2}\right]\right\}^{1/2}\!\!\mathbb{P}^{x}\!\left(\sup_{v\in[0,s]}\!|X_{v}\!-\!x|\geq\frac{\delta_{1,\varepsilon}}{2}\!\right)^{1/2}\!\!\leq 2\sqrt{C_{2}\varepsilon}.
\end{align}
Similarly, using (\ref{eq:IneqHitTimes}) and (\ref{eq:EstT0delta2}), for any $x\in B(x^{(0)},\delta_{\varepsilon})\cap\mathscr{O}$ and any $s\in(0,s_{\varepsilon}]$,
\begin{align}
&\quad\,\,\mathbb{E}^{x}\left[\exp\left(-\int_{0}^{\tau_{\mathscr{O}}}c(X_{v})\,dv\right)\left(g(X_{\tau_{\mathscr{O}}})-g(x^{(0)})\right){\bf 1}_{\{\tau_{\mathscr{O}}\in(s,\infty)\}}\right]\nonumber\\
\label{eq:EstgXtaux03} &\leq\left\{\mathbb{E}^{x}\left[\exp\left(-\int_{0}^{\tau_{\mathscr{O}}}c(X_{v})\,dv\right)\left|g(X_{\tau_{\mathscr{O}}})-g(x^{(0)})\right|^{2}\right]\right\}^{1/2}\mathbb{P}^{x}\left(\tau_{\mathscr{O}}>s\right)^{1/2}\leq 2\sqrt{C_{2}\varepsilon}.
\end{align}
Finally, by combining (\ref{eq:DecomgXtaux0})-(\ref{eq:EstgXtaux03}), for any $x\in B(x^{(0)},\delta_{\varepsilon})\cap\mathscr{O}$, we obtain that
\begin{align*}
\mathbb{E}^{x}\left[\exp\left(-\int_{0}^{\tau_{\mathscr{O}}}c(X_{s})\,ds\right)\left|g(X_{\tau_{\mathscr{O}}})-g(x^{(0)})\right|{\bf 1}_{\{\tau_{\mathscr{O}}<\infty\}}\right]\leq\varepsilon+4\sqrt{C_{2}\varepsilon},
\end{align*}
which, together with (\ref{eq:EstIntTermux})-(\ref{eq:LimitLapTauFinite}), shows that $u_{*}$ is continuous at $x^{(0)}$. Since $x^{(0)}\in\bar{\Gamma}_{0}$ is arbitrarily chosen, we conclude that $u_{*}$ is continuous on $\bar{\Gamma}_{0}$, which completes the proof of the theorem. \hfill $\Box$

\medskip
We now prove existence of solutions to (\ref{eq:CombinedEllipBound}) when the boundary data $g$ is only continuous on suitable portions of the boundary of $\mathscr{O}$.

\medskip
\noindent
\textbf{Proof of Theorem \ref{thm:ExistenceEllipBound}:} For each $x\in\mathscr{O}\cup\partial_{*}\mathscr{O}$, let $(\Omega,\mathscr{F},(\mathscr{F}_{s})_{s\geq 0},\mathbb{P}^{x},W,X)$ be any weak solution to (\ref{eq:MainDegenSDE1})-(\ref{eq:MainDegenSDE2}) with the initial condition (\ref{eq:MainDegenSDEIniCond}) at $t=0$. We need to show that, $u_{*}^{(X)}$, given by (\ref{eq:EllipBoundFeynmanKacTau}), is a solution to (\ref{eq:CombinedEllipBound}), that $u_{*}^{(X)}\in C_{\text{loc}}(\mathscr{O}\cup\partial_{*}\mathscr{O})\cap C^{2,\alpha}(\mathscr{O})$, and that $u_{*}^{(X)}$ satisfies the linear growth condition (\ref{eq:lineargrowth}). Again, for each $x\in\mathscr{O}\cup\partial_{*}\mathscr{O}$, we fix an arbitrary weak solution $(X(x),W)$, with initial data $x$, defined on a single filtered probability space $(\Omega,\mathscr{F},(\mathscr{F}_{s})_{s\geq 0},\mathbb{P})$, and denote $\mathbb{P}^{x}$ and $\mathbb{E}^{x}$ the corresponding probability and expectation. Also, we will omit the superscript $X$ of $u_{*}$ and all stopping times for simplicity.

Since $g\in C_{\text{loc}}(\overline{\partial_{*}\mathscr{O}})$, where $\overline {\partial_{*}\mathscr{O}}$ is a closed set, we may use~\cite[Thoerem 3.1.2]{Friedman:1964} to extend $g$ to $\mathbb{R}^{d}$ such that its extension (denoted by $g$ again) belongs to $C_{\text{loc}}(\mathbb{R}^{d})$. The proof of this theorem is similar to that of Theorem \ref{thm:ExistenceHolderCont}.

\medskip
\noindent
\textbf{Step 1} ($u_{*}\in C^{2,\alpha}(\mathscr{O})$, and $u_{*}$ is a solution to (\ref{eq:CombinedEllipBound}))$\,$ The argument is the same as Step 1 of the proof of Theorem \ref{thm:ExistenceHolderCont}, excluding the part of verifying $u_{*}\in C^{2+\alpha}(\bar{D}')$ at the end.

\medskip
\noindent
\textbf{Step 2} ($u_{*}\in C_{\text{loc}}(\mathscr{O}\cup\partial_{*}\mathscr{O})$)$\,$ For the scenario (B), we may use the same argument as Step 2 of the proof of Theorem \ref{thm:ExistenceHolderCont} to prove that $u_{*}$ is continuous on $\bar{\Gamma}_{0}$. It remains to show that $u_{*}$ is continuous on $\Gamma_{1}$ for both scenarios. That is, for any $x^{(0)}\in\Gamma_{1}$, $\lim_{x\rightarrow x^{(0)}}u_{*}(x)=g(x^{(0)})$.

From the argument in Step 2 of the proof of Theorem \ref{thm:ExistenceHolderCont}, it suffices to show $\lim_{x\rightarrow x^{(0)}}\tau_{\mathscr{O}}^{x}=0$ in probability { (with respect to $\mathbb{P}$)}, for any $x^{(0)}\in\Gamma_{1}$. Now fix any $x^{(0)}\in\Gamma_{1}$, and let $(\mathscr{D}_{k})_{k\in\mathbb{N}}$ be an increasing sequence of $C^{2,\alpha}$ subdomain of $\mathscr{O}$ as in the proof of Theorem \ref{thm:ExistenceHolderCont}. Then, \cite[Theorem 6.13]{GilbargTrudinger:1983} implies that the following elliptic boundary value problem
\begin{align*}
\left\{\begin{array}{ll} \mathscr{A}u=0 &\text{on }\,\mathscr{D}_{k},\\ u=1 &\text{on }\,\partial\mathscr{D}_{k}\end{array}\right.
\end{align*}
admits a unique solution $v_{k}\in C(\bar{\mathscr{D}_{k}})\cap C^{2,\alpha}(\mathscr{D}_{k})$, which moreover admits the unique stochastic representation on $\bar{\mathscr{D}_{k}}$ (cf.~\cite[Theorem 6.5.1]{Friedman:1976}, \cite[Proposition 5.7.2]{KaratzasShreve:1991} and~\cite[Theorem 9.1.1, Corollary 9.1.2]{Oksendal:2003}): for any $x\in\bar{\mathscr{D}_{k}}$,
\begin{align*}
v_{k}(x)=\mathbb{E}^{x}\left[\exp\left(-\int_{0}^{\tau_{\mathscr{D}_{k}}}c(X_{s})\,ds\right)\right].
\end{align*}
Define
\begin{align*}
v_{*}(x)=\mathbb{E}^{x}\left[\exp\left(-\int_{0}^{\tau_{\mathscr{O}}}c(X_{s})\,ds\right)\right],\quad x\in\bar{\mathscr{O}}.
\end{align*}
Similar to the proof of Theorem \ref{thm:ExistenceHolderCont}, and using the continuity of sample paths of $X$ as well as Lemma \ref{lem:EstFunJ}, we have
\begin{align*}
\lim_{k\rightarrow\infty}v_{k}(x)=v_{*}(x),\quad\text{for any }\,x\in\mathscr{O}.
\end{align*}
By Theorem \ref{thm:ExistenceHolderCont}, $v_{*}$ is a solution to
\begin{align*}
\left\{\begin{array}{ll} \mathscr{A}u=0 &\text{on }\,\mathscr{O}, \\ u=1 &\text{on }\,\partial_*\mathscr{O}, \end{array}\right.
\end{align*}
and $v_{*}\in C^{2,\alpha}(\mathscr{O}\cap\Gamma_{1})$. Hence we have
\begin{align*}
\lim_{x\rightarrow x^{(0)}}v_{*}(x)=v_{*}(x^{(0)})=1,
\end{align*}
which implies that $\lim_{x\rightarrow x^{(0)}}\tau_{\mathscr O}^{x}=0$ in probability by Assumption \ref{assupt:UnifPostc}.\hfill $\Box$

\section{Elliptic Obstacle Problems}\label{sec:EllipObsProb}

This section contains the proofs of Theorem \ref{thm:UniqueEllipObs} and \ref{thm:UniqueEllipObsPartial}. By Lemma \ref{lem:BoundClass}, we will prove both theorems in scenarios (A) and (B), as stated at the beginning of Section 3. Also, similar to (\ref{eq:CombinedEllipBound}), the uniqueness of solutions to (\ref{eq:EllipObs}) with the partial Dirichlet boundary condition along $\Gamma_{1}$, when the origin is either a natural (Feller) boundary or an entrance boundary for $X^{(d)}$, and with the full Dirichlet boundary condition along $\partial\mathscr{O}$, when the origin is either a regular boundary or an exit boundary for $X^{(d)}$, are similar in nature. For convenience, we treat them together as
\begin{align}\label{eq:CombinedEllipObs}
\left\{\begin{array}{ll}\min\left\{\mathscr{A}u-f,\,u-\psi\right\}=0\quad\,\text{ on }\,\mathscr{O},\\ u=g\qquad\qquad\qquad\qquad\quad\quad\,\text{on }\,\partial_{*}\mathscr{O},\end{array}\right.
\end{align}
where $\partial_{*}\mathscr{O}$ is given by (\ref{eq:DefPartialOstar}).

\medskip
\noindent
\textbf{Proof of Theorem \ref{thm:UniqueEllipObs}:} We need to show that if $u\in C_{\text{loc}}(\mathscr{O}\cup\partial_{*}\mathscr{O})\cap C^{2}(\mathscr{O})$ is a solution to (\ref{eq:CombinedEllipObs}), which satisfies the linear growth condition (\ref{eq:lineargrowth}), then it admits the stochastic representation (\ref{eq:EllipObsFeynmanKacTau}), for every $x\in\mathscr{O}\cup\partial_{*}\mathscr{O}$, and for any weak solution $(\Omega,\mathscr{F},(\mathscr{F}_{s})_{s\geq 0},\mathbb{P}^{x},W,X)$ to (\ref{eq:MainDegenSDE1})-(\ref{eq:MainDegenSDE2}) with the initial condition (\ref{eq:MainDegenSDEIniCond}) at $t=0$. From (\ref{eq:FunJTheta12}) and (\ref{eq:EllipObsFeynmanKacTau}), we see that
\begin{align}\label{eq:uvstarBound}
u(x)=v_{*}^{(X)}(x)=g(x),\quad\text{for any }\,x\in\partial_{*}\mathscr{O}.
\end{align}
It remains to show that $u=v_{*}$ on $\mathscr{O}$, which we organize in the following two steps. Again, we will fix an arbitrary weak solution for each $x\in\mathscr{O}\cup\partial_{*}\mathscr{O}$, and will omit the superscript $X$ of $v_{*}$.

\medskip
\noindent
\textbf{Step 1} ($u\geq v_{*}$ on $\mathscr{O}$)$\,$ Let $(\mathscr{O}_{k})_{k\in\mathbb{N}}$ be an increasing sequence of $C^{2,\alpha}$ subdomains of $\mathscr{O}$ as in the proof of Theorem \ref{thm:UniqueEllipBound}. For any $x\in\mathscr{O}$, $x\in\mathscr{O}_{k}$ when $k$ is large enough. Since $u\in C^{2}(\mathscr{O})$, by It\^{o}'s formula, for any stopping time $\theta\in\mathscr{T}^{x,X}$ and $s\geq 0$,
\begin{align*}
\mathbb{E}^{x}\!\!\left[\exp\!\left(\!-\!\!\int_{0}^{s\wedge\theta\wedge\tau_{\mathscr{O}_{k}}}\!\!c(X_{v})dv\!\right)\!u\!\left(X_{s\wedge\theta\wedge\tau_{\mathscr{O}_{k}}}\right)\!\right]\!=u(x)-\mathbb{E}^{x}\!\!\left[\int_{0}^{s\wedge\theta\wedge\tau_{\mathscr{O}_{k}}}\!\!\exp\!\left(\!-\!\!\int_{0}^{v}c(X_{w})dw\!\right)\!\mathscr{A}u(X_{v})dv\right].
\end{align*}
By splitting the left-hand side, and using $u\geq\psi$ and $\mathscr{A}u\geq f$ on $\mathscr{O}$, the preceding identity gives
\begin{align}
u(x)&\geq\mathbb{E}^{x}\!\!\left[\exp\!\left(\!-\!\!\int_{0}^{s\wedge\theta}\!\!c(X_{v})dv\!\right)\!\psi(X_{s\wedge\theta}){\bf 1}_{\{\theta<\tau_{\mathscr{O}_{k}}\}}\!\right]\!+\mathbb{E}^{x}\!\!\left[\exp\!\left(\!-\!\!\int_{0}^{s\wedge\tau_{\mathscr{O}_{k}}}\!\!c(X_{v})dv\!\right)\!u\!\left(X_{s\wedge\tau_{\mathscr{O}_{k}}}\right)\!{\bf 1}_{\{\theta\geq\tau_{\mathscr{O}_{k}}\}}\!\right]\nonumber\\
\label{eq:IneqItouOk} &\quad\,+\mathbb{E}^{x}\left[\int_{0}^{s\wedge\theta\wedge\tau_{\mathscr{O}_{k}}}\exp\left(-\int_{0}^{v}c(X_{w})\,dw\right)f(X_{v})\,dv\right].
\end{align}
As in the proof of Theorem \ref{thm:UniqueEllipBound}, the collections of random variables
\begin{align*}
\left\{\exp\left(-\int_{0}^{s\wedge\theta}c(X_{v})\,dv\right)\psi(X_{s\wedge\theta}){\bf 1}_{\{\theta<\tau_{\mathscr{O}_{k}}\}}:\,\,k\in\mathbb{N}\right\}
\end{align*}
and
\begin{align*}
\left\{\exp\left(-\int_{0}^{s\wedge\tau_{\mathscr{O}_{k}}}c(X_{v})\,dv\right)u\left(X_{s\wedge\tau_{\mathscr{O}_{k}}}\right){\bf 1}_{\{\theta\geq\tau_{\mathscr{O}_{k}}\}}:\,\,k\in\mathbb{N}\right\}
\end{align*}
are uniformly integrable because $u$ and $\psi$ satisfy the linear growth condition (\ref{eq:lineargrowth}). By the continuity of $u$ and $\psi$ on $\mathscr{O}\cup\partial_{*}\mathscr{O}$, we also have
\begin{align*}
\lim_{k\rightarrow\infty}\exp\!\left(\!-\!\!\int_{0}^{s\wedge\theta}\!\!c(X_{v})dv\!\right)\!\psi(X_{s\wedge\theta}){\bf 1}_{\{\theta<\tau_{\mathscr{O}_{k}}\}}\!&=\exp\!\left(\!-\!\!\int_{0}^{s\wedge\theta}\!\!c(X_{v})dv\!\right)\!\psi(X_{s\wedge\theta}){\bf 1}_{\{\theta<\tau_{\mathscr{O}}\}},\,\,\,\mathbb{P}^{x}-\text{a.}\,\text{s.}\,,\\
\lim_{k\rightarrow\infty}\!\exp\!\left(\!-\!\!\int_{0}^{s\wedge\tau_{\mathscr{O}_{k}}}\!\!\!c(X_{v})dv\!\right)\!u\!\left(\!X_{s\wedge\tau_{\mathscr{O}_{k}}}\!\right)\!{\bf 1}_{\{\theta\geq\tau_{\mathscr{O}_{k}}\}}\!&=\exp\!\left(\!-\!\!\int_{0}^{s\wedge\tau_{\mathscr{O}}}\!\!\!c(X_{v})dv\!\right)\!u(X_{s\wedge\tau_{\mathscr{O}}}){\bf 1}_{\{\theta\geq\tau_{\mathscr{O}}\}},\,\,\,\mathbb{P}^{x}-\text{a.}\,\text{s.}\,.
\end{align*}
Hence, by~\cite[Theorem 4.5.4]{KailaiChung:2001}, the growth estimate (\ref{eq:EstFirstTermJ12}) and the dominated convergence theorem, we can take the limit in (\ref{eq:IneqItouOk}), as $k\rightarrow\infty$, and obtain that
\begin{align}
u(x)&\geq\mathbb{E}^{x}\!\left[\exp\left(\!-\!\!\int_{0}^{s\wedge\theta}\!\!c(X_{v})dv\!\right)\!\psi(X_{s\wedge\theta}){\bf 1}_{\{\theta<\tau_{\mathscr{O}}\}}\right]\!+\mathbb{E}^{x}\!\left[\exp\left(\!-\!\!\int_{0}^{s\wedge\tau_{\mathscr{O}}}\!\!c(X_{v})dv\!\right)\!u(X_{s\wedge\tau_{\mathscr{O}}}){\bf 1}_{\{\theta\geq\tau_{\mathscr{O}}\}}\right]\nonumber\\
\label{eq:IneqItouO} &\quad\,+\mathbb{E}^{x}\left[\int_{0}^{s\wedge\theta\wedge\tau_{\mathscr{O}}}\exp\left(-\int_{0}^{v}c(X_{w})\,dw\right)f(X_{v})\,dv\right].
\end{align}

Next, we will take $s\rightarrow\infty$ in (\ref{eq:IneqItouO}). Again using a similar argument in the proof of Theorem \ref{thm:UniqueEllipBound}, the collections of random variables
\begin{align*}
\left\{\exp\left(-\int_{0}^{s\wedge\theta}c(X_{v})\,dv\right)\psi(X_{s\wedge\theta}){\bf 1}_{\{\theta<\tau_{\mathscr{O}}\}}:\,\,s\geq 0\right\}
\end{align*}
and
\begin{align*}
\left\{\exp\left(-\int_{0}^{s\wedge\tau_{\mathscr{O}}}c(X_{v})\,dv\right)u(X_{s\wedge\tau_{\mathscr{O}}}){\bf 1}_{\{\theta\geq\tau_{\mathscr{O}}\}}:\,\,s\geq 0\right\}
\end{align*}
are uniformly integrable since $u$ and $\psi$ satisfy the linear growth condition (\ref{eq:lineargrowth}). Also, by the continuity of $u$ and $\psi$ on $\mathscr{O}\cup\partial_{*}\mathscr{O}$,
\begin{align*}
\lim_{s\rightarrow\infty}\exp\!\left(\!-\!\!\int_{0}^{s\wedge\theta}\!c(X_{v})\,dv\!\right)\!\psi(X_{s\wedge\theta}){\bf 1}_{\{\theta<\tau_{\mathscr{O}}\}}&=\exp\!\left(\!-\!\!\int_{0}^{\theta}\!c(X_{v})\,dv\!\right)\!\psi(X_{\theta}){\bf 1}_{\{\theta<\tau_{\mathscr{O}}\}},\quad\mathbb{P}^{x}-\text{a.}\,\text{s.}\,,\\
\lim_{s\rightarrow\infty}\exp\!\left(\!-\!\!\int_{0}^{s\wedge\tau_{\mathscr{O}}}\!\!c(X_{v})dv\!\right)\!u(X_{s\wedge\tau_{\mathscr{O}}}){\bf 1}_{\{\tau_{\mathscr{O}}\leq\theta<\infty\}}&=\exp\!\left(\!-\!\!\int_{0}^{\tau_{\mathscr{O}}}\!\!c(X_{v})dv\!\right)\!u(X_{\tau_{\mathscr{O}}}){\bf 1}_{\{\theta\geq\tau_{\mathscr{O}}\}},\quad\mathbb{P}^{x}-\text{a.}\,\text{s.}\,.
\end{align*}
Therefore, by~\cite[Theorem 4.5.4]{KailaiChung:2001}, the boundary condition (\ref{eq:BoundCondGamma1}), the identity (\ref{eq:SecondFunJTheta12FinTheta1}), the growth estimate (\ref{eq:EstFirstTermJ12}), and the dominated convergence theorem, we can take the limit in (\ref{eq:IneqItouO}), as $s\rightarrow\infty$, and obtain that
\begin{align*}
u(x)&\geq\mathbb{E}^{x}\!\left[\exp\left(-\!\int_{0}^{\theta}\!c(X_{v})\,dv\right)\psi(X_{\theta}){\bf 1}_{\{\theta<\tau_{\mathscr{O}}\}}\right]\!+\mathbb{E}^{x}\!\left[\exp\left(-\!\int_{0}^{\tau_{\mathscr{O}}}\!c(X_{v})\,dv\right)g(X_{\tau_{\mathscr{O}}}){\bf 1}_{\{\theta\geq\tau_{\mathscr{O}},\,\tau_{\mathscr{O}}<\infty\}}\right]\\ &\quad\,+\mathbb{E}^{x}\left[\int_{0}^{\theta\wedge\tau_{\mathscr{O}}}\exp\left(-\int_{0}^{v}c(X_{w})\,dw\right)f(X_{v})\,dv\right]\\
&=\mathbb{E}^{x}\left[\exp\left(-\int_{0}^{\theta}c(X_{v})\,dv\right)\psi(X_{\theta}){\bf 1}_{\{\theta<\tau_{\mathscr{O}}\}}\right]+\mathbb{E}^{x}\left[\exp\left(-\int_{0}^{\tau_{\mathscr{O}}}c(X_{v})\,dv\right)g(X_{\tau_{\mathscr{O}}}){\bf 1}_{\{\theta\geq\tau_{\mathscr{O}}\}}\right]\\
&\quad\,+\mathbb{E}^{x}\left[\int_{0}^{\theta\wedge\tau_{\mathscr{O}}}\exp\left(-\int_{0}^{v}c(X_{w})\,dw\right)f(X_{v})\,dv\right],
\end{align*}
for any $\theta\in\mathscr{T}^{x,X}$ and $x\in\mathscr{O}$, which yields $u\geq v_{*}$ on $\mathscr{O}$.

\medskip
\noindent
\textbf{Step 2} ($u\leq v_{*}$ on $\mathscr{O}$)$\,$ The continuation region
\begin{align}\label{eq:DefContRegion}
\mathscr{C}:=\left\{x\in\mathscr{O}:\,u(x)>\psi(x)\right\}
\end{align}
is an open set of $\mathbb{R}^{d}$ by the continuity of $u$ and $\psi$. We denote
\begin{align}\label{eq:DefTildeTau}
\tilde{\tau}^{t,x,X}:=\inf\left\{s\geq t:\,X_{s}\notin\mathscr{C},\,X_{t}=x\right\},
\end{align}
and write $\tilde{\tau}=\tilde{\tau}^{t,x,X}$ for brevity when $t=0$. $\tilde{\tau}^{t,x,X}$ is indeed a stopping time with respect to $(\mathscr{F}_{s})_{s\geq t}$, since $(X_{s})_{s\geq t}$ is continuous and $\mathscr{O}$ is open. Using the same argument as in Step 1 with $\theta$ replaced by $\tilde{\tau}$, and since $u(X_{\tilde{\tau}})=\psi(X_{\tilde{\tau}})$ and $\mathscr{A}u=f$ on the continuation region $\mathscr{C}$, we obtain that
\begin{align*}
u(x)&=\mathbb{E}^{x}\left[\exp\left(-\int_{0}^{\tilde{\tau}}c(X_{s})\,ds\right)\psi(X_{\tilde{\tau}}){\bf 1}_{\{\tilde{\tau}<\tau_{\mathscr{O}}\}}\right]+\mathbb{E}^{x}\left[\exp\left(-\int_{0}^{\tau_{\mathscr{O}}}c(X_{s})\,ds\right)g(X_{\tau_{\mathscr{O}}}){\bf 1}_{\{\tilde{\tau}\geq\tau_{\mathscr{O}}\}}\right]\\
&\quad\,+\mathbb{E}^{x}\left[\int_{0}^{\tilde{\tau}\wedge\tau_{\mathscr{O}}}\exp\left(-\int_{0}^{s}c(X_{v})\,dv\right)f(X_{s})\,ds\right],
\end{align*}
for any $x\in\mathscr{O}$, which implies that $u\leq v_{*}$ on $\mathscr{O}$. The proof is now complete.\hfill $\Box$

\medskip
\noindent
\textbf{Proof of Theorem \ref{thm:UniqueEllipObsPartial}:}$\,$ We need to show that under the scenario (B), and if $u\in C_{\text{loc}}(\mathscr{O}\cup\Gamma_{1})\cap C^{2}(\mathscr{O})\cap C^{1,1,\beta}_{\text{s,loc}}(\mathscr{O}\cup\Gamma_{0})$ is a solution to (\ref{eq:EllipObs}) with the partial boundary condition (\ref{eq:BoundCondGamma1}), which satisfies the linear growth condition (\ref{eq:lineargrowth}), then it admits the stochastic representation (\ref{eq:EllipObsFeynmanKacLambda}), for any { $x\in\mathscr{O}\cup\Gamma_{1}$, and any} weak solution $(\Omega,\mathscr{F},(\mathscr{F}_{s})_{s\geq 0},\mathbb{P}^{x},W,X)$ to (\ref{eq:MainDegenSDE1})-(\ref{eq:MainDegenSDE2}) with the initial condition (\ref{eq:MainDegenSDEIniCond}) at $t=0$. Again we will omit the superscript $X$ of $v_{**}$, when we fix an arbitrary weak solution for each $x\in\mathscr{O}\cup\Gamma_{1}$. From (\ref{eq:FunJTheta12}) and (\ref{eq:EllipObsFeynmanKacLambda}), we see that
\begin{align}\label{eq:uv2starBound}
u(x)=v_{**}(x)=g(x),\quad\text{for any}\,x\in\Gamma_{1}.
\end{align}
It remains to show $u=v_{**}$ on $\mathscr{O}$. As in the proof of Theorem \ref{thm:UniqueEllipObs}, we organize the proof in the following two steps.

\medskip
\noindent
\textbf{Step 1} ($u\geq v_{**}$ on $\mathscr{O}$)$\,$ Let $\varepsilon>0$ and let $(\mathscr{U}_{k})_{k\in\mathbb{N}}$ be the collection of increasing subdomains of $\mathscr{O}$ as in the proof of Theorem \ref{thm:UniqueEllipBoundPartial}. For any $x\in\mathscr{O}$, $x\in\mathscr{U}_{k}$ when $k$ is large enough. By It\^{o}'s formula, for any $s\geq 0$ and $\theta\in\mathscr{T}^{x,X}$,
\begin{align*}
u(x)\!=\!\mathbb{E}^{x}\!\!\left[\exp\!\left(\!-\!\!\int_{0}^{s\wedge\theta\wedge\lambda_{\mathscr{U}_{k}}}\!\!\!c(X_{v}^{(\varepsilon)})dv\!\right)\!u\!\left(\!X^{(\varepsilon)}_{s\wedge\theta\wedge\lambda_{\mathscr{U}_{k}}}\!\right)\!\right]\!\!+\!\mathbb{E}^{x}\!\!\left[\int_{0}^{s\wedge\theta\wedge\lambda_{\mathscr{U}_{k}}}\!\!\!\!\exp\!\left(\!-\!\!\int_{0}^{v}\!c(X_{w})dw\!\right)\!\mathscr{A}^{\varepsilon}u(X_{v}^{(\varepsilon)})dv\right],
\end{align*}
where $X^{(\varepsilon)}$ is defined by (\ref{eq:DefSDEEpsilon}), and where $\mathscr{A}^{\varepsilon}$ is defined by (\ref{eq:DefGenEps}). By (\ref{eq:DefGenEps}) and using $\mathscr{A}u\geq f$ on $\mathscr{O}$, the preceding identity gives
\begin{align}
u(x)&\geq\mathbb{E}^{x}\!\!\left[\exp\!\left(\!-\!\!\int_{0}^{s\wedge\theta\wedge\lambda_{\mathscr{U}_{k}}}\!\!c(X_{v}^{(\varepsilon)})dv\!\right)\!u\!\left(\!X^{(\varepsilon)}_{s\wedge\theta\wedge\lambda_{\mathscr{U}_{k}}}\!\right)\!\right]\!+\mathbb{E}^{x}\!\!\left[\int_{0}^{s\wedge\theta\wedge\lambda_{\mathscr{U}_{k}}}\!\!\exp\!\left(\!-\!\!\int_{0}^{v}\!c(X_{w}^{(\varepsilon)})dw\!\right)\!f(X^{(\varepsilon)}_{v})dv\right]\nonumber\\
\label{eq:IneqItouUkEps} &\quad\,+\mathbb{E}^{x}\left[\int_{0}^{s\wedge\theta\wedge\lambda_{\mathscr{U}_{k}}}\exp\left(-\int_{0}^{v}c(X_{w}^{(\varepsilon)})\,dw\right)(\mathscr{A}^{\varepsilon}-\mathscr{A})u(X_{v}^{(\varepsilon)})\,dv\right].
\end{align}
Without loss of generality, we assume that $\varepsilon<1/k$, for any fixed large $k\in\mathbb{N}$. By the continuity of $f$ and $u$ on compact subsets of $\mathscr{O}\cup\Gamma_{0}$, as well as the dominated convergence theorem, we see that (\ref{eq:LimitEpsLeft}) and (\ref{eq:LimitEpsRight2}) hold. Also, since the residual term $(\mathscr{A}^{\varepsilon}-\mathscr{A})u$ obeys the estimate (\ref{eq:EstResidAEpsA}), (\ref{eq:LimitEpsRight3}) also holds in the present case. Therefore, by taking $\varepsilon\downarrow 0$ in (\ref{eq:IneqItouUkEps}),
\begin{align}\label{eq:IneqItoUkObs}
u(x)\!\geq\!\mathbb{E}^{x}\!\!\left[\exp\!\left(\!-\!\!\int_{0}^{s\wedge\theta\wedge\lambda_{\mathscr{U}_{k}}}\!\!\!\!c(X_{v})dv\!\right)\!u\!\left(\!X_{s\wedge\theta\wedge\lambda_{\mathscr{U}_{k}}}\!\right)\!\right]\!\!+\!\mathbb{E}^{x}\!\!\left[\!\int_{0}^{s\wedge\theta\wedge\lambda_{\mathscr{U}_{k}}}\!\!\!\!\exp\!\left(\!-\!\!\int_{0}^{v}\!\!c(X_{w})dw\!\right)\!f(X_{v})dv\!\right].
\end{align}
Finally, applying the same argument employed in the proof of Theorem \ref{thm:UniqueEllipBound} and using (\ref{eq:LimitLambdaktinfty}), we can take $k\rightarrow\infty$ and { $s\rightarrow\infty$} in the preceding inequality and obtain that
\begin{align}
u(x)&\geq\mathbb{E}^{x}\left[\exp\left(-\int_{0}^{\theta\wedge\lambda_{\mathscr{O}}}c(X_{v})\,dv\right)u(X_{\theta\wedge\lambda_{\mathscr{O}}})\right]\!+\mathbb{E}^{x}\left[\int_{0}^{\theta\wedge\lambda_{\mathscr{O}}}\exp\left(-\int_{0}^{v}c(X_{w})\,dw\right)f(X_{v})\,dv\right]\nonumber\\
&\geq\mathbb{E}^{x}\!\left[\exp\left(\!-\!\!\int_{0}^{\theta}c(X_{v})dv\!\right)\!\psi(X_{\theta}){\bf 1}_{\{\theta<\lambda_{\mathscr{O}}\}}\!\right]\!+\mathbb{E}^{x}\!\left[\exp\left(\!-\!\!\int_{0}^{\lambda_{\mathscr{O}}}c(X_{v})dv\!\right)\!g(X_{\lambda_{\mathscr{O}}}){\bf 1}_{\{\theta\geq\lambda_{\mathscr{O}},\,\lambda_{\mathscr{O}}<\infty\}}\!\right]\nonumber\\
&\quad\,+\mathbb{E}^{x}\left[\int_{0}^{\theta\wedge\lambda_{\mathscr{O}}}\exp\left(-\int_{0}^{v}c(X_{w})\,dw\right)f(X_{v})\,dv\right]\nonumber\\
&=\mathbb{E}^{x}\left[\exp\left(-\int_{0}^{\theta}c(X_{v})\,dv\right)\psi(X_{\theta}){\bf 1}_{\{\theta<\lambda_{\mathscr{O}}\}}\right]+\mathbb{E}^{x}\left[\exp\left(-\int_{0}^{\lambda_{\mathscr{O}}}c(X_{v})\,dv\right)g(X_{\lambda_{\mathscr{O}}}){\bf 1}_{\{\theta\geq\lambda_{\mathscr{O}}\}}\right]\nonumber\\
\label{eq:IneqItoObs} &\quad\,+\mathbb{E}^{x}\left[\int_{0}^{\theta\wedge\lambda_{\mathscr{O}}}\exp\left(-\int_{0}^{v}c(X_{w})\,dw\right)f(X_{v})\,dv\right],
\end{align}
for any $\theta\in\mathscr{T}^{x,X}$ and $x\in\mathscr{O}$, where in the second inequality we have used $u\geq\psi$ on $\mathscr{O}\cup\Gamma_{0}$, which follows from (\ref{eq:EllipObs}) and the continuity of $u$ and $\psi$ up to $\Gamma_{0}$, and where in the third equation we have used the identity (\ref{eq:SecondFunJTheta12FinTheta1}). This concludes that $u\geq v_{**}$ on $\mathscr{O}$.

\medskip
\noindent
\textbf{Step 2} ($u\leq v_{**}$ on $\mathscr{O}$)$\,$ We choose $\theta=\tilde{\tau}$ in the preceding step, where $\tilde{\tau}$ is defined by (\ref{eq:DefTildeTau}) with $t=0$. By the definition of the continuation region $\mathscr{C}$ given as in (\ref{eq:DefContRegion}), and the obstacle problem (\ref{eq:EllipObs}), the inequalities (\ref{eq:IneqItouUkEps}), (\ref{eq:IneqItoUkObs}) and (\ref{eq:IneqItoObs}) hold with equality. Therefore, we conclude that $u\leq v_{**}$ on $\mathscr{O}$, which completes the proof. \hfill $\Box$

\section{Parabolic Terminal/Boundary-Value and Obstacle Problems}\label{sec:ParProb}

In this section, we will derive Feynman-Kac formulas for solutions to parabolic terminal/boundary-value problem (\ref{eq:ParaTermBound}) and obstacle problem (\ref{eq:ParaTermObs}) with partial/full boundary conditions. Recall that $Q=(0,T)\times\mathscr{O}$, where $T\in(0,\infty)$ is fixed, and where $\mathscr{O}$ is a (possibly unbounded) connected, open subset of the upper half-space $\mathbb{H}=\mathbb{R}^{d-1}\times(0,\infty)$ such that $\Gamma_{0}=\partial\mathscr{O}\cap\partial\mathbb{H}\neq\emptyset$. We will need to appeal to the following analogue of Assumption \ref{assupt:LinearGrowth}:
\begin{assumption}\label{assupt:GrowthPara} \emph{(Linear growth condition)}
If $u$ is a vector-valued or matrix-valued function on a subset of $[0,\infty)\times\overline{\mathbb{H}}$, there exists a universal constant $K>0$, such that
\begin{align}\label{eq:GrowthPara}
\|u(t,x)\|\leq K\left(1+\|x\|\right)
\end{align}
on its domain.
\end{assumption}
Let $C(Q)$ denote the vector space of continuous function on $Q$, while $C(\overline{Q})$ denotes the Banach space of functions which are uniformly continuous and bounded on $\overline{Q}$. Let $Du$ and $D^{2}u$ denote the gradient and the Hessian matrix, respectively, of a function $u$ on $Q$ with respect to spatial variables. Let $C^{1}(Q)$ denote the vector space of functions $u$ such that, $u$, $u_{t}$ and $Du$ are continuous on $Q$, while $C^{1}(\overline{Q})$ denotes the Banach space of functions $u$ such that, $u$, $u_{t}$ and $Du$ are uniformly continuous and bounded on $\bar{Q}$. Finally, let $C^{2}(Q)$ denote the vector space of functions $u$ such that, $u$, $u_{t}$, $Du$ and $D^{2}u$ are continuous on $Q$, while $C^{2}(\overline{Q})$ denotes the Banach space of functions $u$ such that, $u$, $u_{t}$, $Du$ and $D^{2}u$ are uniformly continuous and bounded on $\overline{Q}$. If $T\subsetneqq\partial Q$ is a relatively open set, let $C_{\text{loc}}(Q\cup T)$ denote the vector space of functions $u$ such that, for any precompact open subset $V\Subset Q\cup T$, $u\in C(\overline{V})$.

For any $(t,x)\in\overline{Q}$, let $(\Omega,\mathscr{F},(\mathscr{F}_{s})_{s\geq t},\mathbb{P}^{t,x},W^{(t)},X^{(t)})$ be an arbitrary weak solution to (\ref{eq:MainDegenSDE1})-(\ref{eq:MainDegenSDE2}) with the initial condition (\ref{eq:MainDegenSDEIniCond}) (recall that the existence of such a solution is guaranteed by Assumption \ref{assupt:Cont} and Assumption \ref{assupt:LinearGrowth}). We define
\begin{align}
u_{*}^{(X)}(t,x)&=\mathbb{E}^{t,x}\left[\exp\left(-\int_{t}^{\tau_{\mathscr{O}}\wedge T}c(X_{s})\,ds\right)g\left(\tau_{\mathscr{O}}\wedge T,X_{\tau_{\mathscr{O}}\wedge T}\right)\right]\nonumber\\
\label{eq:ParaBoundFeynmanKacTau} &\quad\,+\mathbb{E}^{t,x}\left[\int_{t}^{\tau_{\mathscr{O}}\wedge T}\exp\left(-\int_{t}^{s}c(X_{v})\,dv\right)f(s,X_{s})\,ds\right],\\
u_{**}^{(X)}(t,x)&=\mathbb{E}^{t,x}\left[\exp\left(-\int_{t}^{\lambda_{\mathscr{O}}\wedge T}c(X_{s})\,ds\right)g\left(\lambda_{\mathscr{O}}\wedge T,X_{\lambda_{\mathscr{O}}\wedge T}\right)\right]\nonumber\\
\label{eq:ParaBoundFeynmanKacLambda} &\quad\,+\mathbb{E}^{t,x}\left[\int_{t}^{\lambda_{\mathscr{O}}\wedge T}\exp\left(-\int_{t}^{s}c(X_{v})\,dv\right)f(s,X_{s})\,ds\right],
\end{align}
where $\tau_{\mathscr{O}}=\tau_{\mathscr{O}}^{t,x,X}$ and $\lambda_{\mathscr{O}}=\lambda_{\mathscr{O}}^{t,x,X}$ are defined as in (\ref{eq:StopTimeTau}) and (\ref{eq:StopTimeLambda}), respectively. Let $\mathscr{T}_{t,T}^{x,X}$ be the collection of all $(\mathscr{F}_{s})_{s\in[t,T]}$-stopping times taking values in $[t,T]$. For any $\theta_{1},\,\theta_{2}\in\mathscr{T}_{t,T}^{x,X}$, define
\begin{align}
\mathcal{J}^{\theta_{1},\theta_{2}}_{X}(t,x)\!&:=\!\mathbb{E}^{t,x}\!\!\left[\int_{t}^{\theta_{1}\wedge\theta_{2}}\!\!\!\exp\!\left(\!-\!\!\int_{t}^{s}\!c(X_{v})dv\!\right)\!f(s,\!X_{s})ds\right]\!\!+\!\mathbb{E}^{t,x}\!\!\left[\exp\!\left(\!-\!\!\int_{t}^{\theta_{1}}\!c(X_{s})ds\!\right)\!g(\theta_{1},\!X_{\theta_{1}}){\bf 1}_{\{\theta_{1}\leq\theta_{2}\}}\!\right]\nonumber\\
\label{eq:ParaFunJTheta12} &\quad+\mathbb{E}^{t,x}\left[\exp\left(-\int_{t}^{\theta_{2}}c(X_{s})\,ds\right)\psi(\theta_{2},X_{\theta_{2}}){\bf 1}_{\{\theta_{1}>\theta_{2}\}}\right],
\end{align}
and
\begin{align}\label{eq:ParaObsTau}
v_{*}^{(X)}(t,x)&:=\sup_{\theta\in\mathscr{T}_{t,T}^{x,X}}\mathcal J^{\tau_{\mathscr{O}}\wedge T,\theta}_{X}(t,x),\\
\label{eq:ParaObsLambda} v_{**}^{(X)}(t,x)&:=\sup_{\theta\in\mathscr{T}_{t,T}^{x,X}}\mathcal J^{\lambda_{\mathscr{O}}\wedge T,\theta}_{X}(t,x).
\end{align}
Above, and in the sequel, we omit the superscripts $t$ and $x$ for all random variables inside the probability $\mathbb{P}^{t,x}$ and the expectation $\mathbb{E}^{t,x}$.

As an analog of Lemma \ref{lem:EstFunJ} in the parabolic case, we have the following estimate on the function $\mathcal{J}^{\theta_{1},\theta_{2}}_{X}$. In particular, the functions $u_{*}^{(X)}$, $u_{**}^{(X)}$, $v_{*}^{(X)}$ and $v_{**}^{(X)}$, given respectively by (\ref{eq:ParaBoundFeynmanKacTau}), (\ref{eq:ParaBoundFeynmanKacLambda}), (\ref{eq:ParaObsTau}) and (\ref{eq:ParaObsLambda}), are well defined and satisfy the linear growth condition (\ref{eq:GrowthPara}). The proof is similar to that of Lemma \ref{lem:EstFunJ}, and is thus omitted.
\begin{lemma}\label{lem:EstParaFunJ}
Fix $T>0$. Let $f$, $g$ and $\psi$ are real-valued Borel measurable functions on $[0,T]\times\mathbb{R}^{d}$ satisfying the linear growth condition (\ref{assupt:GrowthPara}). Assume that the coefficients functions $b$, $\sigma$ and $c$ satisfy (\ref{eq:SupMartCon}). Then, for any $(t,x)\in\overline{Q}$, any weak solution $(\Omega,\mathscr{F},(\mathscr{F}_{s})_{s\geq t},\mathbb{P}^{t,x},W^{(t)},X^{(t)})$ to (\ref{eq:MainDegenSDE1})-(\ref{eq:MainDegenSDE2}) with the initial condition (\ref{eq:MainDegenSDEIniCond}), and any $\theta_{1},\,\theta_{2}\in\mathscr{T}_{t,T}^{x,X}$, the function $\mathcal{J}^{\theta_{1},\theta_{2}}_{X}$, given by (\ref{eq:ParaFunJTheta12}), satisfies
\begin{align}\label{eq:ParaEstJ12}
\left|\mathcal{J}^{\theta_{1},\theta_{2}}_{X}(t,x)\right|\leq C(1+\|x\|),
\end{align}
where $C$ is a universal positive constant, depending only on $K$ as in (\ref{eq:GrowthPara}), $c_{0}$ as in Assumption \ref{assupt:UnifPostc}, and $M$ as in (\ref{eq:SupMartM}).
\end{lemma}
By Lemma \ref{lem:BoundClass}, in the sequel, we will prove the uniqueness and existence theorems for the parabolic terminal/boundary-value and obstacle problems in scenarios (A) and (B), which was stated at the beginning of Section \ref{sec:EllipBoundValProb}. Define
\begin{align}\label{eq:DefPartialQstar}
\eth_{*}Q:=\left\{\begin{array}{ll}\eth^{1}Q, &\text{if the scenario (A) occurs},\\
\eth Q, &\text{if the scenario (B) occurs},\end{array}\right.
\end{align}
and treat the previous mentioned terminal/boundary-value problems together as
\begin{align}\label{eq:CombinedParaBound}
\left\{\begin{array}{ll} -u_{t}+\mathscr{A}u=f &\text{in }\,Q,\\ u=g &\text{on }\,\eth_{*}Q.\end{array}\right.
\end{align}

\subsection{Feynman-Kac Formulas for Parabolic Terminal/Boundary-Value Problem}

We first establish the uniqueness of Feynman-Kac formulas for solutions to the parabolic terminal/boundary value problem (\ref{eq:ParaTermBound}) with either the partial terminal/boundary condition (\ref{eq:TermBoundPart}), or the full terminal/boundary condition (\ref{eq:TermBoundFull}).
\begin{theorem}\label{thm:UniqueParaBound}
Assume that $b$, $\sigma$ and $c$ satisfy (\ref{eq:SupMartCon}), that $b$ and $\sigma$ obey the linear growth condition (\ref{eq:lineargrowth}), and that $f\in C(Q)$ which obeys the linear growth condition (\ref{eq:GrowthPara}) on $Q$,
\begin{itemize}
\item [(1)] Suppose that either case \emph{(a)}, \emph{(b)} or \emph{(c)} in Theorem \ref{thm:UniqueEllipBound} occurs. Assume that $g\in C_{\emph{loc}}(\eth^{1}Q)$ which obeys (\ref{eq:GrowthPara}) on $\eth^{1}Q$. Let
    \begin{align*}
    u\in C_{\emph{loc}}(Q\cup\eth^{1}Q)\cap C^{2}(Q)
    \end{align*}
    be a solution to the parabolic terminal/boundary value problem (\ref{eq:ParaTermBound}) and (\ref{eq:TermBoundPart}), and which obeys (\ref{eq:GrowthPara}) on $Q$. Then, { for any $(t,x)\in Q\cup\eth^{1}Q$, $u(t,x)=u_{*}^{(X)}(t,x)$}, for any weak solution { $(\Omega,\mathscr{F},(\mathscr{F}_{s})_{s\geq t},\mathbb{P}^{t,x},W^{(t)},X^{(t)})$} to (\ref{eq:MainDegenSDE1})-(\ref{eq:MainDegenSDE2}) with the initial condition (\ref{eq:MainDegenSDEIniCond}), where $u_{*}^{(X)}$ is given by (\ref{eq:ParaBoundFeynmanKacTau}).
\item [(2)] Suppose that either case \emph{(d)} or \emph{(e)} in Theorem \ref{thm:UniqueEllipBound} occurs. Assume that $g\in C_{\emph{loc}}(\eth Q)$ which obeys (\ref{eq:GrowthPara}) on $\eth Q$. Let
    \begin{align*}
    u\in C_{\emph{loc}}(Q\cup \eth Q)\cap C^{2}(Q)
    \end{align*}
    be a solution to the parabolic terminal/boundary value problem (\ref{eq:ParaTermBound}) and (\ref{eq:TermBoundFull}), and which obeys (\ref{eq:GrowthPara}) on $Q$. Then, { for any $(t,x)\in Q\cup\eth Q$, $u(t,x)=u_{*}^{(X)}(t,x)$}, for any weak solution { $(\Omega,\mathscr{F},(\mathscr{F}_{s})_{s\geq t},\mathbb{P}^{t,x},W^{(t)},X^{(t)})$} to (\ref{eq:MainDegenSDE1})-(\ref{eq:MainDegenSDE2}) with the initial condition (\ref{eq:MainDegenSDEIniCond}), where $u_{*}^{(X)}$ is given by (\ref{eq:ParaBoundFeynmanKacTau}).
\end{itemize}
\end{theorem}
\noindent
\textbf{Proof:} We need to show that if $u\in C_{\text{loc}}(Q\cup\eth_{*}Q)\cap C^{2}(Q)$ is a solution to (\ref{eq:CombinedParaBound}), which satisfies the linear growth condition (\ref{eq:GrowthPara}), then it admits the stochastic representation (\ref{eq:ParaBoundFeynmanKacTau}), for any $(t,x)\in Q\cup\eth_{*}Q$, and any weak solution $(\Omega,\mathscr{F},(\mathscr{F}_{s})_{s\geq t},\mathbb{P}^{t,x},W^{(t)},X^{(t)})$ to (\ref{eq:MainDegenSDE1})-(\ref{eq:MainDegenSDE2}) with the initial condition (\ref{eq:MainDegenSDEIniCond}). From the expression of $u_{*}^{(X)}$ and the boundary condition of $u$ in (\ref{eq:CombinedParaBound}), it is easy to see that $u(t,x)=u_{*}^{(X)}(t,x)=g(t,x)$ on $\eth_{*}Q$ for any weak solution. It remains to justify $u=u_{*}^{(X)}$ on $Q$ for any weak solution. The proof is similar as that of Theorem \ref{thm:UniqueEllipBound}, and we just outline the sketch here.

Let $(\mathscr{O}_{k})_{k\in\mathbb{N}}$ be an increasing sequence of $C^{2,\alpha}$ open subdomains of $\mathscr{O}$ as in the proof of Theorem \ref{thm:UniqueEllipBound}, with $\alpha\in(0,1)$, such that $\bar{\mathscr{O}}_{k}\subseteq\mathscr{O}$ for each $k\in\mathbb{N}$, and $\cup_{k\in\mathbb{N}}\mathscr{O}_{k}=\mathscr{O}$. For any $(t,x)\in Q$, we fix any arbitrary weak solution $(\Omega,\mathscr{F},(\mathscr{F}_{s})_{s\geq t},\mathbb{P}^{t,x},W^{(t)},X^{(t)})$ to (\ref{eq:MainDegenSDE1})-(\ref{eq:MainDegenSDE2}) with the initial condition (\ref{eq:MainDegenSDEIniCond}), and obviously $(t,x)\in\mathscr{O}_{k}$ for $k$ large enough. By It\^{o}'s formula,
\begin{align}
&\quad\,\,\mathbb{E}^{t,x}\left[\exp\left(-\int_{t}^{T\wedge\tau_{\mathscr{O}_{k}}}c(X_{s})\,ds\right)u\left(\tau_{\mathscr{O}_{k}}\wedge T,X_{T\wedge\tau_{\mathscr{O}_{k}}}\right)\right]\nonumber\\
\label{eq:ItoPara} &=u(t,x)-\mathbb{E}^{t,x}\left[\int_{t}^{T\wedge\tau_{\mathscr{O}_{k}}}\exp\left(-\int_{t}^{s}c(X_{v})\,dv\right)f(s,X_{s})\,ds\right].
\end{align}
We need to take the limit in (\ref{eq:ItoPara}) as $k\rightarrow\infty$. By a parabolic version of the growth estimate (\ref{eq:EstFirstTermJ12}) (the condition in Lemma \ref{lem:EstFunJ} can be much weaker for the stopping times bounded by $T$, due to~\cite[Problem 5.3.15]{KaratzasShreve:1991}), we can apply the dominated convergence theorem to obtain that
\begin{align*}
\lim_{k\rightarrow\infty}\mathbb{E}^{t,x}\!\left[\int_{t}^{T\wedge\tau_{\mathscr{O}_{k}}}\!\exp\!\left(\!-\!\int_{t}^{s}c(X_{v})\,dv\!\right)\!f(s,X_{s})\,ds\right]\!=\mathbb{E}^{t,x}\!\left[\int_{t}^{T\wedge\tau_{\mathscr{O}}}\!\exp\!\left(\!-\!\int_{t}^{s}c(X_{v})\,dv\!\right)\!f(s,X_{s})\,ds\right].
\end{align*}
For the non-integral term on the left-hand side of (\ref{eq:ItoPara}), by the continuity of $u$ and of sample paths of $(X_{s})_{s\geq t}$, we first have
\begin{align*}
\lim_{k\rightarrow\infty}\exp\!\left(\!-\!\!\int_{t}^{T\wedge\tau_{\mathscr{O}_{k}}}\!\!c(X_{s})ds\!\right)\!u\!\left(T\!\wedge\!\tau_{\mathscr{O}_{k}},X_{T\wedge\tau_{\mathscr{O}_{k}}}\right)\!=\exp\!\left(\!-\!\!\int_{t}^{T\wedge\tau_{\mathscr{O}}}\!\!c(X_{s})ds\!\right)\!u\!\left(T\!\wedge\!\tau_{\mathscr{O}},X_{T\wedge\tau_{\mathscr{O}}}\right),\quad\text{a.}\,\text{s.}\,.
\end{align*}
By~\cite[Theorem 4.5.4]{KailaiChung:2001}, in order to show the convergence of corresponding expectations, we only need to show that
\begin{align*}
\left\{\exp\left(-\int_{t}^{T\wedge\tau_{\mathscr{O}_{k}}}c(X_{s})\,ds\right)u\left(T\wedge\tau_{\mathscr{O}_{k}},X_{T\wedge\tau_{\mathscr{O}_{k}}}\right):\,\,k\in\mathbb{N}\right\}
\end{align*}
is a collection of uniformly integrable random variables. To do this, it suffices to show that their second moments are uniformly bounded, which follows from~\cite[Problem 3.15]{KaratzasShreve:1991} and the linear growth condition (\ref{eq:GrowthPara}) on $u$. The proof is now complete.\hfill $\Box$

\medskip
Let $C_{s,\text{loc}}^{1,1,\beta}((0,T)\times(\mathscr{O}\cup\Gamma_{0}))$ denote the subspace of $C^{1,1}((0,T)\times(\mathscr{O}\cup\Gamma_{0}))\cap C^{2}_{\text{loc}}((0,T)\times(\mathscr{O}\cup\Gamma_{0}))$ consisting of functions $\varphi$ such that, for any pre-compact open subset $V\Subset [0,T]\times(\mathscr{O}\cup\Gamma_{0})$,
\begin{align*}
\sup_{(t,x)\in V}\left(|\varphi(t,x)|+\|D\varphi(t,x)\|+\left\|x_{d}^{\beta}D^{2}\varphi(t,x)\right\|\right)<\infty.
\end{align*}
We then have the following alternative uniqueness result.
\begin{theorem}\label{thm:UniqueParaBoundPartial}
Suppose that either case \emph{(d)} or \emph{(e)} in Theorem \ref{thm:UniqueEllipBound} occurs. Let $b$, $\sigma$ and $c$ be as in Theorem \ref{thm:UniqueParaBound}, and let $c\in C_{\emph{loc}}(\mathscr{O}\cup\Gamma_{0})$. Assume that $f,g\in C_{\emph{loc}}(\eth^{1}Q)$ which obey the linear growth condition (\ref{eq:GrowthPara}) on $\eth^{1}Q$. Let
\begin{align*}
u\in C_{\emph{loc}}(Q\cup\eth^{1}Q)\cap C^{2}(Q)\cap C_{s,\emph{loc}}^{1,1,\beta}((0,T)\times(\mathscr{O}\cup\Gamma_{0}))
\end{align*}
be a solution to the parabolic terminal/boundary value problem (\ref{eq:ParaTermBound}) and (\ref{eq:TermBoundFull}), and which obeys (\ref{eq:GrowthPara}) on $Q$. Then, for any $(t,x)\in Q\cup\eth^{1}Q$, we have $u(t,x)=u_{**}^{(X)}(t,x)$, for any weak solution $(\Omega,\mathscr{F},(\mathscr{F}_{s})_{s\geq t},\mathbb{P}^{t,x},W^{(t)},X^{(t)})$ to (\ref{eq:MainDegenSDE1})-(\ref{eq:MainDegenSDE2}) with the initial condition (\ref{eq:MainDegenSDEIniCond}), where $u_{**}^{(X)}$ is given by (\ref{eq:ParaBoundFeynmanKacLambda}).
\end{theorem}
\noindent
\textbf{Proof:} For $\varepsilon>0$, by It\^{o}'s formula, we have
\begin{align*}
&\quad\,\,\mathbb{E}^{t,x}\left[\exp\left(-\int_{0}^{t\wedge\lambda_{\mathscr{U}_{k}}}c(X^{(\varepsilon)}_{s})\,ds\right)u\left(t\wedge\lambda_{\mathscr{U}_{k}},X^{(\varepsilon)}_{t\wedge\lambda_{\mathscr{U}_{k}}}\right)\right]\\ &=u(t,x)-\mathbb{E}^{t,x}\left[\int_{0}^{t\wedge\lambda_{\mathscr{U}_{k}}}\exp\left(-\int_{0}^{s}c(X^{(\varepsilon)}_{v})\,dv\right)\mathscr{A}^{\varepsilon}u(s,X^{(\varepsilon)}_{s})\,ds\right],
\end{align*}
where $\mathscr{U}_{k}$, $\mathscr{A}^{\varepsilon}$ and $X^{(\varepsilon)}$ were defined as in the proof of Theorem \ref{thm:UniqueEllipBoundPartial}. Now the proof follows from the same argument as that of Theorem \ref{thm:UniqueEllipBoundPartial}. \hfill $\Box$

\medskip
Similar to the elliptic case, we have the following two results on the existence of solutions to the parabolic terminal/boundary value problem with continuous and H\"{o}lder continuous terminal/boundary data, respectively.
\begin{theorem}\label{thm:ExistenceParaBound}
In addition to the hypothesis of Theorem \ref{thm:UniqueParaBound}, assume that $\Gamma_{1}$ is of class $C^{2,\alpha}$, that $b,\sigma\in C^{0,\alpha}(\mathscr{O})$, and that $f\in C^{0,\alpha}(Q)$, for some $\alpha\in(0,1)$.
\begin{itemize}
\item [(1)] Suppose that either case \emph{(a)}, \emph{(b)} or \emph{(c)} in Theorem \ref{thm:UniqueEllipBound} occurs. Assume that $g\in C_{\emph{loc}}(\overline{\eth^{1}Q})$ which obeys (\ref{eq:GrowthPara}) on $\eth^{1}Q$. For any $(t,x)\in Q\cup\eth^{1}Q$, let $(\Omega,\mathscr{F},(\mathscr{F}_{s})_{s\geq t},\mathbb{P}^{t,x},W^{(t)},X^{(t)})$ be a weak solution to (\ref{eq:MainDegenSDE1})-(\ref{eq:MainDegenSDE2}) with the initial condition (\ref{eq:MainDegenSDEIniCond}), and let $u_{*}^{(X)}$ be defined as in (\ref{eq:ParaBoundFeynmanKacTau}). Then, $u_{*}^{(X)}$ is a solution to (\ref{eq:ParaTermBound}) with the terminal/boundary condition (\ref{eq:TermBoundPart}) along $\eth^{1}Q$. In particular, $u_{*}^{(X)}\in C_{\emph{loc}}(Q\cup\eth^{1}Q)\cap C^{2,\alpha}(Q)$ which obeys (\ref{eq:GrowthPara}) on $Q\cup\eth^{1}Q$.
\item [(2)] Suppose that either case \emph{(d)} or \emph{(e)} in Theorem \ref{thm:UniqueEllipBound} occurs. Assume that $g\in C_{\emph{loc}}(\overline{\eth Q})$ which obeys (\ref{eq:GrowthPara}) on $\eth Q$. For any $(t,x)\in Q\cup\eth Q$, let $(\Omega,\mathscr{F},(\mathscr{F}_{s})_{s\geq t},\mathbb{P}^{t,x},W^{(t)},X^{(t)})$ be a weak solution to (\ref{eq:MainDegenSDE1})-(\ref{eq:MainDegenSDE2}) with the initial condition (\ref{eq:MainDegenSDEIniCond}), and let $u_{*}^{(X)}$ be defined as in (\ref{eq:ParaBoundFeynmanKacTau}). Then, $u_{*}^{(X)}$ is a solution to (\ref{eq:ParaTermBound}) with the terminal/boundary condition (\ref{eq:TermBoundFull}) along $\eth Q$. In particular, $u_{*}^{(X)}\in C_{\emph{loc}}(Q\cup\eth Q)\cap C^{2,\alpha}(Q)$ which obeys (\ref{eq:GrowthPara}) on $Q\cup\eth Q$.
\end{itemize}
\end{theorem}
\noindent
\textbf{Proof:} The proof is similar to that of Theorem \ref{thm:ExistenceEllipBound}, and we will just sketch the outline. Again for simplicity, we will omit the superscript $X$ of $u_{*}$ since, for each $(t,x)\in Q\cup\eth_{*}Q$, we fix an arbitrary weak solution $(\Omega,\mathscr{F},(\mathscr{F}_{s})_{s\geq t},\mathbb{P}^{t,x},W^{(t)},X^{(t)})$. From the expression of (\ref{eq:ParaBoundFeynmanKacLambda}), clearly we have $u_{*}=g$ on $\eth_{*}Q$. By~\cite[Theorem 3.1.2]{Friedman:1964}, we may extend $g\in C_{\text{loc}}(\overline{\eth_{*}Q})$ to a function (called $g$ again) on $[0,T]\times\mathbb{R}^{d}$, such that its extension belongs to $C_{\text{loc}}([0,T]\times\mathbb{R}^{d})$. Let $(\mathscr{O}_{k})_{k\in\mathbb{N}}$ be an increasing sequence of $C^{2,\alpha}$ open subdomains of $\mathscr{O}$ as in the proof of Theorem \ref{thm:UniqueEllipBound}, with $\alpha\in(0,1)$, such that $\bar{\mathscr{O}}_{k}\subseteq\mathscr{O}$ for each $k\in\mathbb{N}$, and $\cup_{k\in\mathbb{N}}\mathscr{O}_{k}=\mathscr{O}$. Let $Q_{k}:=(0,T)\times\mathscr{O}_{k}$ for each $k\in\mathbb{N}$. On $Q_{k}$, by~\cite[Theorem 3.4.9]{Friedman:1964}, the terminal/boundary value problem
\begin{align*}
\left\{\begin{array}{ll} -u_{t}+\mathscr{A}u=f &\text{in }\,Q_{k}, \\ u=g &\text{on }\,\left((0,T)\times\partial\mathscr{O}_{k}\right)\cup\left(\{T\}\times\bar{\mathscr{O}}_{k}\right),
\end{array}\right.
\end{align*}
has a unique solution $u_{k}\in C(\bar{Q}_k)\cap C^{2,\alpha}(Q_{k})$, and by~\cite[Theorem 6.5.2]{Friedman:1976} it admits the stochastic representation: for any $(t,x)\in\left((0,T)\times\partial\mathscr{O}_{k}\right)\cup\left(\{T\}\times\bar{\mathscr{O}}_{k}\right)$,
\begin{align*}
u_{k}(t,x)&=\mathbb{E}^{t,x}\left[\exp\left(-\int_{t}^{\tau_{\mathscr{O}_{k}}\wedge T}c(X_{s})\,ds\right)g\left(\tau_{\mathscr{O}_{k}}\wedge T,X_{\tau_{\mathscr{O}_{k}}\wedge T}\right)\right]\\
&\quad\,+\mathbb{E}^{t,x}\left[\int_{t}^{\tau_{\mathscr{O}_{k}}\wedge T}\exp\left(-\int_{t}^{s}c(X_{v})\,dv\right)f(s,X_{s})\,ds\right].
\end{align*}
Here, we note that for any weak solution to (\ref{eq:MainDegenSDE1})-(\ref{eq:MainDegenSDE2}), with the same initial data $(t,x)\in Q$, has the unique law up to $\tau_{\mathscr{O}}^{t,x,X}$. Since $\tau_{\mathscr{O}_{k}}\rightarrow\tau_{\mathscr{O}}$ a.$\,$s., as $k\rightarrow\infty$, and using the same argument as in the proof of Theorem \ref{thm:UniqueParaBound}, we can show the convergence of the right-hand side of the above equation. Therefore, as $k\rightarrow\infty$, $u_{k}$ converges to $u_{*}$ pointwisely in $Q$.

With the help of interior Schauder estimate for parabolic equations~\cite[Exercise 10.4.2]{Krylov:1996}, and using the same argument as in the first step of proof for Theorem \ref{thm:ExistenceHolderCont}, we can obtain that $u_{*}\in C^{2,\alpha}(Q)$. To get the continuity of $u_{*}$ up to the boundary $\eth_{*}Q$, we may use the same argument as in the second step of the proof of Theorem \ref{thm:ExistenceHolderCont}. \hfill $\Box$
\begin{theorem}\label{thm:ExistenceParaHolderCont}
In addition to the hypotheses of Theorem \ref{thm:UniqueParaBound}, assume that $b,\sigma\in C^{0,\alpha}(\mathscr{O})$, and that $f\in C^{0,\alpha}_{\emph{loc}}(Q\cup\eth^{1}Q)$, for some $\alpha\in(0,1)$.
\begin{itemize}
\item [(1)] Suppose that either case \emph{(a)}, \emph{(b)} or \emph{(c)} in Theorem \ref{thm:UniqueEllipBound} occurs. Assume that the boundary portion $\Gamma_{1}$ is of class $C^{2,\alpha}$, and that $g\in C_{\emph{loc}}^{2,\alpha}(Q\cup\eth^{1}Q)$ which obeys the linear growth condition (\ref{eq:GrowthPara}) on $Q\cup\eth^{1}Q$, and
    \begin{align}\label{compatibility1}
    -g_{t}+\mathscr{A}g=f\quad\text{ on }\,\{T\}\times\Gamma_{1}.
    \end{align}
    For any $(t,x)\in Q\cup\eth^{1}Q$, let $(\Omega,\mathscr{F},(\mathscr{F}_{s})_{s\geq t},\mathbb{P}^{t,x},W^{(t)},X^{(t)})$ be a weak solution to (\ref{eq:MainDegenSDE1})-(\ref{eq:MainDegenSDE2}) with the initial condition (\ref{eq:MainDegenSDEIniCond}), and let $u_{*}^{(X)}$ be defined as in (\ref{eq:ParaBoundFeynmanKacTau}). Then, $u_{*}^{(X)}$ is a solution to (\ref{eq:ParaTermBound}) with the terminal/boundary condition (\ref{eq:TermBoundPart}) along $\eth^{1}Q$. In particular, $u_{*}^{(X)}\in C_{\emph{loc}}^{2,\alpha}(Q\cup\eth^{1}Q)$ which obeys (\ref{eq:GrowthPara}) on $Q\cup\eth^{1}Q$.
\item [(2)] Suppose that either case \emph{(d)} or \emph{(e)} in Theorem \ref{thm:UniqueEllipBound} occurs. Assume that the boundary $\partial\mathscr{O}$ is of class $C^{2,\alpha}$, and that $g\in C_{\emph{loc}}^{2,\alpha}(Q\cup\eth^{1}Q)\cap C_{\emph{loc}}(Q\cup\eth Q)$ which obeys the linear growth condition (\ref{eq:GrowthPara}) on $Q\cup\eth Q$, and
    \begin{align}\label{compatibility2}
    -g_{t}+\mathscr{A}g=f\quad\text{ on }\,\{T\}\times\partial\mathscr{O}.
    \end{align}
    For any $(t,x)\in Q\cup\eth Q$, let $(\Omega,\mathscr{F},(\mathscr{F}_{s})_{s\geq t},\mathbb{P}^{t,x},W^{(t)},X^{(t)})$ be a weak solution to (\ref{eq:MainDegenSDE1})-(\ref{eq:MainDegenSDE2}) with the initial condition (\ref{eq:MainDegenSDEIniCond}), and let $u_{*}^{(X)}$ be defined as in (\ref{eq:ParaBoundFeynmanKacTau}). Then, $u_{*}^{(X)}$ is a solution to (\ref{eq:ParaTermBound}) with the terminal/boundary condition (\ref{eq:TermBoundPart}). In particular, $u_{*}^{(X)}\in C_{\emph{loc}}^{2,\alpha}(Q\cup\eth^{1}Q)$ which obeys (\ref{eq:GrowthPara}) on $Q\cup\eth Q$.
\end{itemize}
\end{theorem}
\noindent
\textbf{Proof:} We may use the same strategy as in the proof of Theorem \ref{thm:ExistenceHolderCont} for this parabolic case. The only difference is that, when we prove $u_{*}\in C^{2,\alpha}(Q)$ (or $u_{*}\in C^{2,\alpha}(\eth_{*}Q)$), where $u_{*}$ is a solution that we obtain by a similar limiting argument as in the proof of Theorem \ref{thm:ExistenceParaBound}, we use interior Schauder estimate~\cite[Exercise 10.4.2]{Krylov:1996} (or boundary Schauder estimate~\cite[Proposition A.1]{FeehanPop:2013(2)}) for parabolic equations, instead of~\cite[Corollary 6.3]{GilbargTrudinger:1983} (or~\cite[Corollary 6.7]{GilbargTrudinger:1983}) for elliptic case. \hfill $\Box$
\begin{remark}\label{rem:ExistPara}
In contrast to Theorem \ref{thm:ExistenceHolderCont} for the existence of solutions to the elliptic terminal/boundary-value problem with H\"{o}lder continuous terminal/boundary conditions, the parabolic case requires the compatibility conditions (\ref{compatibility1}) and (\ref{compatibility2}) (cf.~\cite[Section 10.4]{Krylov:1996}).
\end{remark}

\subsection{Feynman-Kac Formulas for Parabolic Obstacle Problems}

In this final subsection, we briefly investigate the Feynman-Kac Formulas for parabolic obstacle problem (\ref{eq:ParaTermObs}) with partial/full terminal/boundary conditions. Similar to the elliptic case, those two scenarios (see Section \ref{sec:EllipBoundValProb}) depending on whether the process $X$ reaches $\Gamma_{0}$, can be united as one obstacle problem
\begin{align*}
\left\{\begin{array}{ll}\min\{-u_{t}+\mathscr{A}u-f,u-\psi\}=0 &\text{in }\,Q, \\ u=g &\text{on }\,\eth_{*}Q,\end{array}\right.
\end{align*}
where $\eth_{*}Q $ is defined in (\ref{eq:DefPartialQstar}). The proofs of the following two theorems are then similar to those of Theorem \ref{thm:UniqueEllipObs} and Theorem \ref{thm:UniqueEllipObsPartial}, respectively, and is thus omitted.
\begin{theorem}\label{thm:UniqueParaObs}
Let $f, b, \sigma$ and $c$ be as in Theorem \ref{thm:UniqueParaBound}, and let $\psi\in C(Q)$ which obeys the linear growth condition (\ref{eq:GrowthPara}).
\begin{itemize}
\item [(1)] Suppose that either case \emph{(a)}, \emph{(b)} or \emph{(c)} in Theorem \ref{thm:UniqueEllipBound} occurs. Assume that $\psi\in C_{\emph{loc}}(Q\cup\eth^{1}Q)$, and that $g\in C_{\emph{loc}}(\eth^{1}Q)$ which obeys (\ref{eq:CompatObsPartPara}) and (\ref{eq:GrowthPara}) on $\eth^{1}Q$. Let
    \begin{align*}
    u\in C_{\emph{loc}}(Q\cup\eth^{1}Q)\cap C^{2}(Q)
    \end{align*}
    be a solution to the parabolic obstacle problem (\ref{eq:ParaTermObs}) and (\ref{eq:TermBoundPart}), such that both $u$ and $\mathscr{A}u$ obey (\ref{eq:GrowthPara}) on $Q$. Then, for any $(t,x)\in Q\cup\eth^{1}Q$, $u(t,x)=v_{*}^{(X)}(t,x)$, for any weak solution $(\Omega,\mathscr{F},(\mathscr{F}_{s})_{s\geq t},\mathbb{P}^{t,x},W^{(t)},X^{(t)})$ to (\ref{eq:MainDegenSDE1})-(\ref{eq:MainDegenSDE2}) with the initial condition (\ref{eq:MainDegenSDEIniCond}), where $v_{*}^{(X)}$ is given by (\ref{eq:ParaObsTau}).
\item [(2)] Suppose that either case \emph{(d)} or \emph{(e)} in Theorem \ref{thm:UniqueEllipBound} occurs. Assume that $\psi\in C_{\emph{loc}}(Q\cup\eth Q)$, and that $g\in C_{\emph{loc}}(\eth Q)$ which obeys (\ref{eq:CompatObsFullPara}) and (\ref{eq:GrowthPara}) on $\eth Q$. Let
    \begin{align*}
    u\in C_{\emph{loc}}(Q\cup\eth Q)\cap C^{2}({Q})
    \end{align*}
    be a solution to the elliptic obstacle problem (\ref{eq:ParaTermObs}) and (\ref{eq:TermBoundFull}), such that  both $u$ and $\mathscr{A}u$ obey (\ref{eq:GrowthPara}) on $Q$. Then, for any $(t,x)\in Q\cup\eth Q$, $u(t,x)=v_{*}^{(X)}(t,x)$, for any weak solution $(\Omega,\mathscr{F},(\mathscr{F}_{s})_{s\geq t},\mathbb{P}^{t,x},W^{(t)},X^{(t)})$ to (\ref{eq:MainDegenSDE1})-(\ref{eq:MainDegenSDE2}) with the initial condition (\ref{eq:MainDegenSDEIniCond}), where $v_{*}^{(X)}$ is given by (\ref{eq:ParaObsTau}).
\end{itemize}
\end{theorem}
\begin{theorem}\label{thm:UniqueParaObsPartial}
Suppose that either case \emph{(d)} or \emph{(e)} in Theorem \ref{thm:UniqueEllipBound} occurs. Let $f$, $b$, $\sigma$ and $c$ be as in Theorem \ref{thm:UniqueParaBoundPartial}. Let $\psi\in C_{\emph{loc}}(Q\cup\eth^{1}Q)$ which obeys the linear growth condition (\ref{eq:GrowthPara}) on $Q$, and let $g\in C_{\emph{loc}}(\eth^{1}Q)$ which obeys (\ref{eq:CompatObsPartPara}) and (\ref{eq:GrowthPara}) on $\eth^{1}Q$. Let
\begin{align*}
u\in C_{\emph{loc}}(Q\cup\eth^{1}Q)\cap C^{2}(Q)\cap C_{\emph{s,loc}}^{1,1,\beta}(Q\cup (0,T)\times(\mathscr{O}\cup\Gamma_{0}))
\end{align*}
is a solution to the parabolic obstacle problem (\ref{eq:ParaTermObs}) and (\ref{eq:TermBoundPart}), such that both $u$ and $\mathscr{A}u$ obey (\ref{eq:GrowthPara}). Then, for any $(t,x)\in Q\cup\eth^{1}Q$, $u(t,x)=v_{**}^{(X)}(t,x)$, for any weak solution $(\Omega,\mathscr{F},(\mathscr{F}_{s})_{s\geq t},\mathbb{P}^{t,x},W^{(t)},X^{(t)})$ to (\ref{eq:MainDegenSDE1})-(\ref{eq:MainDegenSDE2}) with the initial condition (\ref{eq:MainDegenSDEIniCond}), where $v_{**}^{(X)}$ is given by (\ref{eq:ParaObsLambda}).
\end{theorem}

\end{document}